\documentclass[11pt]{article}
\usepackage{amssymb}
\usepackage{amsmath}
\usepackage[all]{xy}
\usepackage{times}
\usepackage{graphicx}

\title{On small analytic relations\indent}
\author{Dominique LECOMTE}
\date{\today}

\def\ufootnote#1{\let\savedthfn\thefootnote\let\thefootnote\relax
\footnote{#1}\let\thefootnote\savedthfn\addtocounter{footnote}{-1}}

\newcommand{\Ana}{{\it\Sigma}^{1}_{1}}

\newcommand{\Ca}{{\it\Pi}^{1}_{1}}

\newcommand{\Borel}{{\it\Delta}^{1}_{1}}
\newcommand{\ana}{{\bf\Sigma}^{1}_{1}}

\newcommand{\boraone}{{\bf\Sigma}^{0}_{1}}
\newcommand{\boratwo}{{\bf\Sigma}^{0}_{2}}

\newcommand{\boraxi}{{\bf\Sigma}^{0}_{\xi}}

\newcommand{\borone}{{\bf\Delta}^{0}_{1}}
\newcommand{\bortwo}{{\bf\Delta}^{0}_{2}}

\newcommand{\bormone}{{\bf\Pi}^{0}_{1}}

\newcommand{\bormtwo}{{\bf\Pi}^{0}_{2}}

\newcommand{\bormxi}{{\bf\Pi}^{0}_{\xi}}

\newcommand{\borxi}{{\bf\Delta}^{0}_{\xi}}
\newcommand{\borme}{{\bf\Pi}^{0}_{\eta}}

\newtheorem{thm} {Theorem} [section]
\newtheorem{defi} [thm] {Definition}
\newtheorem{cor} [thm] {Corollary}

\newtheorem{lem} [thm] {Lemma}

\begin{document}

\maketitle

\centerline{$\bullet$ 1) Sorbonne Universit\' e, CNRS, Institut de Math\'ematiques de Jussieu-Paris Rive Gauche,}

\centerline{IMJ-PRG, F-75005 Paris, France}

\centerline{2) Universit\'e de Paris, IMJ-PRG, F-75013 Paris, France}

\centerline{dominique.lecomte@upmc.fr}\medskip

\centerline{$\bullet$ Universit\'e de Picardie, I.U.T. de l'Oise, site de Creil,}

\centerline{13, all\'ee de la fa\"\i encerie, 60 100 Creil, France}\medskip\medskip\medskip\medskip\medskip

\ufootnote{{\it 2010 Mathematics Subject Classification.}~Primary: 03E15, Secondary: 28A05, 54H05}

\ufootnote{{\it Keywords and phrases.}~acyclic, analytic relation, antichain basis, Borel class, countable, descriptive complexity, continuous reducibility, graph, Ramsey}

\noindent {\bf Abstract.} We study the class of analytic binary relations on Polish spaces, compared with the notions of continuous reducibility or injective continuous reducibility. In particular, we characterize when a locally countable Borel relation is $\boraxi$ (or 
$\bormxi$), when $\xi\!\geq\! 3$, by providing a concrete finite antichain basis. We give a similar characterization for arbitrary relations when $\xi\! =\! 1$. When $\xi\! =\! 2$, we provide a concrete antichain of size continuum made of locally countable Borel relations minimal among non-$\boratwo$ (or non-$\bormtwo$) relations. The proof of this last result allows us to strengthen a result due to Baumgartner in topological Ramsey theory on the space of rational numbers. We prove that positive results hold when $\xi\! =\! 2$ in the acyclic case. We give a general positive result in the non-necessarily locally countable case, with another suitable acyclicity assumption. We provide a concrete finite antichain basis for the class of uncountable analytic relations. Finally, we deduce from our positive results some antichain basis for graphs, of small cardinality (most of the time 1 or 2). 

\vfill\eject

\section{$\!\!\!\!\!\!$ Introduction}\indent

 This article presents a continuation of the work in [L5], in which the descriptive complexity of Borel equivalence relations on Polish spaces was studied (recall that a topological space is {\bf Polish} if it is separable and completely metrizable). These relations are compared using the notion of {\bf continuous reducibility}, which is as follows. Recall that if $X,Y$ are topological spaces and $A\!\subseteq\! X^2$, 
$B\!\subseteq\! Y^2$, 
$$(X,A)\leq_c(Y,B)\Leftrightarrow\exists f\! :\! X\!\rightarrow\! Y\mbox{ continuous with }A\! =\! (f\!\times\! f)^{-1}(B)$$
(we say that $f$ {\bf reduces} $A$ to $B$). When the function $f$ can be injective, we write $\sqsubseteq_c$ instead of 
$\leq_c$. Sometimes, when the space is clear for instance, we will talk about $A$ instead of $(X,A)$. The motivation for considering these quasi-orders is as follows (recall that a {\bf quasi-order} is a reflexive and transitive relation). A  standard way of comparing the descriptive complexity of subsets of zero-dimensional Polish spaces is the {\bf Wadge quasi-order} (see [W]; recall that a topological space is 
{\bf zero-dimensional} if there is a basis for its topology made of clopen, i.e., closed and open, sets). If $S,Z$ are zero-dimensional Polish spaces and $C\!\subseteq\! S$, $D\!\subseteq\! Z$, 
$$(S,C)\leq_W(Z,D)\Leftrightarrow\exists g\! :\! S\!\rightarrow\! Z\mbox{ continuous with }C\! =\! g^{-1}(D).$$
However, the pre-image of a graph by an arbitrary continuous map is not in general a graph, for instance. Note that the classes of reflexive relations, irreflexive relations, symmetric relations, transitive relations are closed under pre-images by a square map. Moreover, the class of antisymmetric relations is closed under pre-images by the square of an injective map. This is the reason why square maps are considered to compare graphs,  equivalence relations...  The most common way of comparing Borel equivalence relations is the notion of Borel reducibility (see, for example, [G], [Ka]). However, very early in the theory, injective continuous reducibility was considered, for instance in Silver's theorem (see [S]).\medskip

 The most classical hierarchy of topological complexity in descriptive set theory is the one given by the Borel classes. If $\bf\Gamma$ is a class of subsets of the metrizable spaces, then $\check {\bf\Gamma}\! :=\!\{\neg S\mid S\!\in\! {\bf\Gamma}\}$ is its {\bf dual class}, and $\Delta ({\bf\Gamma})\! :=\! {\bf\Gamma}\cap\check {\bf\Gamma}$. Recall that the Borel hierarchy is the inclusion from left to right in the following picture:\medskip
 
\scalebox{0.75}{$$\!\!\!\!\!\!\!\!\!\!\!\!\!\!\!\!\!\!\!\!\!\!\!\xymatrix@1{ 
& & \boraone\! =\!\mbox{open} & & \boratwo\! =\! F_\sigma & & & 
\boraxi\! =\! (\bigcup_{\eta <\xi}~\borme )_\sigma & \\ 
& \borone\! =\!\mbox{clopen} & & \bortwo\! =\!\boratwo\cap\bormtwo & & \cdots & \borxi\! =\!\boraxi\cap\bormxi & & \cdots\\ 
& & \bormone\! =\!\mbox{closed} & & \bormtwo\! =\! G_\delta & & & \bormxi\! =\!\check\boraxi & }$$}\medskip
 
\noindent This hierarchy is strict in uncountable Polish spaces, in which the non self-dual classes are those of the form $\boraxi$ or 
$\bormxi$. In the sequel, by non self-dual Borel class, we mean exactly those classes. A class $\bf\Gamma$ of subsets of the zero-dimensional Polish spaces is a {\bf Wadge class} if there is a zero-dimensional space $Z$ and a subset $D$ of $Z$ in $\bf\Gamma$ such that a subset $C$ of a zero-dimensional space $S$ is in $\bf\Gamma$ exactly when $(S,C)\leq_W(Z,D)$. The hierarchy of the Wadge classes of Borel sets, compared with the inclusion, refines greatly the hierarchy of the non self-dual Borel classes, and is the finest hierarchy of topological complexity considered in descriptive set theory (see  [Lo-SR2]).

\vfill\eject

 We are interested in the descriptive complexity of Borel relations on Polish spaces. In order to approach this problem, it is useful to consider invariants for the considered quasi-order. A natural invariant for Borel reducibility has been studied, the notion of potential complexity (see, for example, [L2], [L3], and [Lo2] for the definition). A Borel relation $R$ on a Polish space $X$ is {\bf potentially} in a Wadge class $\bf\Gamma$ if we can find a finer Polish topology $\tau$ on $X$ such that $R$ is in $\bf\Gamma$ in the product 
$(X,\tau )^2$. This is an invariant in the sense that any relation which is Borel reducible to a relation potentially in 
$\bf\Gamma$ has also to be potentially in $\bf\Gamma$. Along similar lines, any relation which is continuously reducible to a relation in 
$\bf\Gamma$ has also to be in $\bf\Gamma$.\medskip

 We already mentioned the equivalence relations. A number of other interesting relations can be considered on a Polish space $X$, in the descriptive set theoretic context. Let us mention\medskip
 
\noindent - the {\bf digraphs} (which do not meet the {\bf diagonal} $\Delta (X)\! :=\!\{ (x,x)\mid x\!\in\! X\}$ of $X$),\smallskip

\noindent - the {\bf graphs} (i.e., the symmetric digraphs),\smallskip

\noindent - the {\bf oriented graphs} (i.e., the antisymmetric digraphs),\smallskip

\noindent - the {\bf quasi-orders}, strict or not,\smallskip

\noindent - the {\bf partial orders} (i.e., the antisymmetric quasi-orders), strict or not.\medskip

 For instance, we refer to [Lo3], [L1], [K-Ma]. For {\bf locally countable} relations (i.e., relations with countable horizontal and vertical sections), we refer to [K2] in the case of equivalence relations. An important subclass of the class of locally countable Borel equivalence relations is the class of treeable locally countable Borel equivalence relations, generated by an acyclic locally countable Borel graph. More generally, the locally countable digraphs have been widely considered, not necessarily to study equivalence relations (see [K-Ma]). All this motivates the  work in the present paper, mostly devoted to the study of the descriptive complexity of arbitrary locally countable or/and acyclic Borel relations on Polish spaces.\medskip 
 
 We are looking for characterizations of the relations in a fixed Borel class $\bf\Gamma$. So we will consider the continuous and injective continuous reducibilities. In other words, we want to give answers to the following very simple question: when is a relation ${\bf\Sigma}^0_\xi$ (or ${\bf\Pi}^0_\xi$)? We are looking for characterizations of the following form: a relation is either simple, or more complicated than a typical complex relation. So we need to introduce, for some Borel classes $\bf\Gamma$, examples of complex relations.\medskip

\noindent\bf Notation.\rm\ Let $\bf\Gamma$ be a non self-dual Borel class, $\mathbb{K}$ be a metrizable compact space, and $\mathbb{C}\!\in\!\check {\bf\Gamma}(\mathbb{K})\!\setminus\! {\bf\Gamma}(\mathbb{K})$.\medskip 

 If the {\bf rank} of $\bf\Gamma$ is one (i.e., if ${\bf\Gamma}\!\in\!\{\boraone ,\bormone\}$), then we set 
$\mathbb{K}\! :=\!\{ 0\}\cup\{ 2^{-k}\mid k\in\omega\}\!\subseteq\!\mathbb{R}$, $\mathbb{C}\! :=\!\{ 0\}$ if ${\bf\Gamma}\! =\!\boraone$, and 
$\mathbb{C}\! :=\!\mathbb{K}\!\setminus\!\{ 0\}$ if ${\bf\Gamma}\! =\!\bormone$, since we want some injectivity results.\medskip

 If the rank of $\bf\Gamma$ is at least two, then we set $\mathbb{K}\! :=\! 2^\omega$, and 
$\mathbb{C}\cap N_s\!\in\!\check {\bf\Gamma}(N_s)\!\setminus\! {\bf\Gamma}(N_s)$ for each $s\!\in\! 2^{<\omega}$ (this is possible, by Lemma 4.5 in [L5]). In particular, $\mathbb{C}$ is dense and co-dense in $2^\omega$. We set 
$$\mathbb{C}\! :=\!\mathbb{P}_\infty\! :=\!\{\alpha\!\in\! 2^\omega\mid\exists^\infty n\!\in\!\omega ~~\alpha (n)\! =\! 1\}$$ 
if ${\bf\Gamma}\! =\!\boratwo$, and 
$\mathbb{C}\! :=\!\mathbb{P}_f\! :=\!\{\alpha\!\in\! 2^\omega\mid\forall^\infty n\!\in\!\omega ~~\alpha (n)\! =\! 0\}$ if  
${\bf\Gamma}\! =\!\bormtwo$, for injectivity reasons again. In the sequel, we will say that $(\mathbb{K},\mathbb{C})$ is 
$\bf\Gamma$-{\bf good} if it satisfies all the properties mentioned here.\medskip

\noindent\bf Examples.\rm\ Let ${\bf\Gamma}$ be a non self-dual Borel class, and $(\mathbb{K},\mathbb{C})$ be $\bf\Gamma$-good. We define a relation on 
$\mathbb{D}\! :=\! 2\!\times\!\mathbb{K}$ as follows: 
$(\varepsilon ,x)~\mathbb{E}^{\bf\Gamma}_3~(\eta ,y)\Leftrightarrow (\varepsilon ,x)\! =\! (\eta ,y)\vee (x\! =\! y\!\in\!\mathbb{C})$. The graph 
$\mathbb{G}^{\bf\Gamma}_m\! :=\!\mathbb{E}^{\bf\Gamma}_3\!\setminus\!\Delta (\mathbb{D})$ will be very important in the sequel.\smallskip

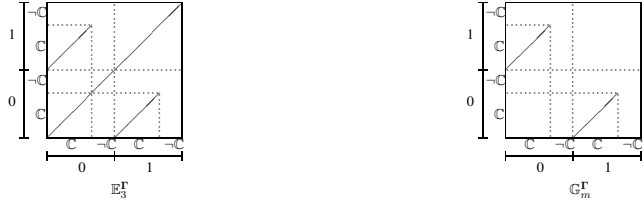
\scalebox{0.5}{\setlength{\unitlength}{0.6cm}
\begin{picture}(0,6)(-12,-6)
   
   \multiput(-0.9,-4)(0.2,0){30}{\circle*{0.05}}
   \put(-1,-7){\line(1,1){6}}
   \put(-1,-7){\line(1,0){6}}
   \put(-1,-7){\line(0,1){6}}
   \put(-1,-1){\line(1,0){6}}
   \put(5,-7){\line(0,1){6}}
   \multiput(-0.9,-5)(0.2,0){25}{\circle*{0.05}}
   \multiput(1,-6.9)(0,0.2){25}{\circle*{0.05}}
   \multiput(2,-6.9)(0,0.2){30}{\circle*{0.05}}
   \multiput(-0.9,-2)(0.2,0){10}{\circle*{0.05}}
   \multiput(4,-6.9)(0,0.2){10}{\circle*{0.05}}
   \put(-1,-4){\line(1,1){2}}
   \put(2,-7){\line(1,1){2}}
   \put(-0.15,-7.5){$\mathbb{C}$}
   \put(1.2,-7.5){$\neg\mathbb{C}$}
   \put(2.85,-7.5){$\mathbb{C}$}
   \put(4.2,-7.5){$\neg\mathbb{C}$}
   \put(-1,-7.8){\line(1,0){6}}
   \put(-1,-8){\line(0,1){0.4}}
   \put(2,-8){\line(0,1){0.4}}
   \put(5,-8){\line(0,1){0.4}}
   \put(0.4,-8.5){0}
   \put(3.4,-8.5){1}
   \put(-1.5,-6.15){$\mathbb{C}$}
   \put(-1.87,-4.65){$\neg\mathbb{C}$}
   \put(-1.5,-3.15){$\mathbb{C}$}
   \put(-1.87,-1.65){$\neg\mathbb{C}$}
   \put(-2,-7){\line(0,1){6}}
   \put(-2.2,-7){\line(1,0){0.4}}
   \put(-2.2,-4){\line(1,0){0.4}}
   \put(-2.2,-1){\line(1,0){0.4}}
   \put(-2.7,-5.6){0}
   \put(-2.7,-2.6){1}
   \put(1.85,-9.5){$\mathbb{E}^{\bf\Gamma}_3$}
     
\end{picture}

\setlength{\unitlength}{0.6cm}
\begin{picture}(0,6)(-32,-6)
   
   \multiput(-0.9,-4)(0.2,0){30}{\circle*{0.05}}
   \put(-1,-7){\line(1,0){6}}
   \put(-1,-7){\line(0,1){6}}
   \put(-1,-1){\line(1,0){6}}
   \put(5,-7){\line(0,1){6}}
   \multiput(-0.9,-5)(0.2,0){25}{\circle*{0.05}}
   \multiput(1,-6.9)(0,0.2){25}{\circle*{0.05}}
   \multiput(2,-6.9)(0,0.2){30}{\circle*{0.05}}
   \multiput(-0.9,-2)(0.2,0){10}{\circle*{0.05}}
   \multiput(4,-6.9)(0,0.2){10}{\circle*{0.05}}
   \put(-1,-4){\line(1,1){2}}
   \put(2,-7){\line(1,1){2}}
   \put(-0.15,-7.5){$\mathbb{C}$}
   \put(1.2,-7.5){$\neg\mathbb{C}$}
   \put(2.85,-7.5){$\mathbb{C}$}
   \put(4.2,-7.5){$\neg\mathbb{C}$}
   \put(-1,-7.8){\line(1,0){6}}
   \put(-1,-8){\line(0,1){0.4}}
   \put(2,-8){\line(0,1){0.4}}
   \put(5,-8){\line(0,1){0.4}}
   \put(0.4,-8.5){0}
   \put(3.4,-8.5){1}
   \put(-1.5,-6.15){$\mathbb{C}$}
   \put(-1.87,-4.65){$\neg\mathbb{C}$}
   \put(-1.5,-3.15){$\mathbb{C}$}
   \put(-1.87,-1.65){$\neg\mathbb{C}$}
   \put(-2,-7){\line(0,1){6}}
   \put(-2.2,-7){\line(1,0){0.4}}
   \put(-2.2,-4){\line(1,0){0.4}}
   \put(-2.2,-1){\line(1,0){0.4}}
   \put(-2.7,-5.6){0}
   \put(-2.7,-2.6){1}
   \put(1.85,-9.5){$\mathbb{G}^{\bf\Gamma}_m$}
     
\end{picture}}\medskip\medskip\medskip\medskip\medskip\medskip\medskip

\noindent The main result in [L5] is as follows. Most of our results will hold in analytic spaces and not only in Polish spaces. Recall that a separable metrizable space is an {\bf analytic space} if it is homeomorphic to an analytic subset of a Polish space.

\begin{thm} Let ${\bf\Gamma}$ be a non self-dual Borel class of rank at least three. Then $\mathbb{E}^{\bf\Gamma}_3$ is 
$\sqsubseteq_c$-minimum among non-$\bf\Gamma$ locally countable Borel equivalence relations on an analytic space.\end{thm}

 In fact, this result is also valid for equivalence relations with $\boratwo$ classes, i.e., $\boratwo$ sections. Recall that if $(Q,\leq )$ is a quasi-ordered class, then a {\bf basis} is a subclass $B$ of $Q$ such that any element of $Q$ is $\leq$-above an element of $B$. We are looking for basis as small as possible, so in fact for antichains (an {\bf antichain} is a subclass of $Q$ made of pairwise $\leq$-incomparable elements). So we want antichain basis. As we will see, the solution of our problem heavily depends on the rank of the non self-dual Borel class considered. The next result solves our problem for the classes of rank at least three.

\begin{thm} \label{morethree} (1) Let ${\bf\Gamma}$ be a non self-dual Borel class of rank at least three. Then there is a concrete 34 elements $\sqsubseteq_c$ and $\leq_c$-antichain basis for the class of non-$\bf\Gamma$ locally countable Borel relations on an analytic space.\smallskip

\noindent (2) $\mathbb{G}^{\bf\Gamma}_m$ is $\sqsubseteq_c$-minimum among non-$\bf\Gamma$ locally countable Borel graphs on an analytic space.\smallskip

\noindent (3) These results also hold when $\bf\Gamma$ has rank two for relations whose sections are in $\Delta ({\bf\Gamma})$.\end{thm}

 The next surprising result shows that this complexity assumption on the sections is useful for the classes of rank two, for which any basis must have size continuum.
 
\begin{thm} \label{two} Let ${\bf\Gamma}$ be a non self-dual Borel class of rank two. Then there is a concrete $\leq_c$-antichain of size continuum made of locally countable Borel relations on $2^\omega$ which are $\leq_c$ and $\sqsubseteq_c$-minimal among 
non-$\bf\Gamma$ relations on an analytic space.\end{thm}

 Similar results hold for graphs (see Corollary \ref{Fsigmastruct} and Theorem \ref{Gdeltastruct}). Our analysis of the rank two also provides a basis for the class of non-$\boratwo$ locally countable Borel relations on an analytic space (see Theorem \ref{twocases}). The proof of Theorem \ref{two} strengthens Theorem 1.1 in [B] (see also Theorem 6.31 in [T], in topological Ramsey theory on the space $\mathbb{Q}$ of rational numbers).
 
\begin{thm} \label{baum} There is an onto coloring $c\! :\!\mathbb{Q}^{[2]}\!\rightarrow\!\omega$ with the property that, for any $H\!\subseteq\!\mathbb{Q}$ homeomorphic to $\mathbb{Q}$, there is $h\! :\!\mathbb{Q}\!\rightarrow\! H$, homeomorphism onto its range, for which $c(\{ x,y\} )\! =\! c\big(\{ h(x),h(y)\}\big)$ if $x,y\!\in\!\mathbb{Q}$. In particular, $c$ takes all the values from $\omega$ on $H^{[2]}$.\end{thm}
 
\noindent\bf Question.\rm\ For the classes of rank two, we saw that any basis must have size continuum. Is there an antichain basis?\medskip
 
 The next result solves our problem for the classes of rank one.

\vfill\eject

\begin{thm} \label{one} Let ${\bf\Gamma}$ be a non self-dual Borel class of rank one. Then there is a concrete 7360 elements $\leq_c$-antichain basis, made of relations on a countable metrizable compact space, for the class of non-$\bf\Gamma$ relations on a first countable topological space. A similar result holds for $\sqsubseteq_c$, with 2 more elements in the antichain basis.\end{thm}

 Note that in Theorem \ref{one}, the fact of assuming that $X$ is analytic or that $R$ is locally countable Borel does not change the result since the elements of the antichain basis satisfy these stronger assumptions. The ``first countable" assumption ensures that closures and sequential closures coincide. Here again, similar results hold for graphs, with much smaller antichain basis (of cardinality $\leq_c$-5, $\sqsubseteq_c$-6 for $\bormone$, and 10 for $\boraone$). We will not give the proof of Theorem \ref{one} since it is elementary. We simply describe the different antichain basis in Section 5.\medskip
 
\noindent\bf Remark.\rm\ Theorem \ref{one} provides a finite antichain basis for the class of non-closed Borel relations on a Polish space, for $\sqsubseteq_c$ and $\leq_c$. This situation is very different for the class $\mathcal{C}$ of non-potentially closed Borel relations on a Polish space. Indeed, [L1] provides an antichain of size continuum made of minimal elements of $\mathcal{C}$, for any of these two quasi-orders. It also follows from [L-M] that in fact there is no antichain basis in $\mathcal{C}$, for any of these two quasi-orders again.\medskip

 The works in [K-S-T], [L-M], [L-Z], [L4] and also [C-L-M] show that an acyclicity assumption can give positive dichotomy results (see, for example, Theorem 1.9 in [L4]). This is a way to fix our problem with the rank two. If $A$ is a binary relation on a set $X$, then 
${A^{-1}\! :=\!\{ (x,y)\!\in\! X^2\mid (y,x)\!\in\! A\}}$. The {\bf symmetrization} of $A$ is $s(A)\! :=\! A\cup A^{-1}$. An $A$-{\bf path} is a finite sequence $(x_i)_{i\leq n}$ of points of $X$ such that $(x_i,x_{i+1})\!\in\! A$ if $i\! <\! n$. We say that $A$ is {\bf acyclic} if there is no injective $A$-path $(x_i)_{i\leq n}$ with $n\!\geq\! 2$ and $(x_n,x_0)\!\in\! A$.  In practice, we will consider acyclicity only for symmetric relations. We will say that $A$ is {\bf s-acyclic} if $s(A)$ is acyclic.

\begin{thm} \label{Sigma02ac} (1) There is a concrete 34 elements $\sqsubseteq_c$ and $\leq_c$-antichain basis for the class of non-$\boratwo$ s-acyclic locally countable Borel relations on an analytic space.\smallskip

\noindent (2) $\mathbb{G}^{\boratwo}_m$ is $\sqsubseteq_c$-minimum among non-$\boratwo$ acyclic locally countable Borel graphs on an analytic space.\end{thm}

\noindent\bf Notation.\rm\ For the class $\bormtwo$, we need some more examples since the complexity of a locally countable relation can come from the complexity of one of its sections in this case. Let ${\bf\Gamma}$ be a non self-dual Borel class of rank at least two, and $(2^\omega ,\mathbb{C})$ be $\bf\Gamma$-good. We set 
${\mathbb S}\! :=\!\{ 0^\infty\}\cup N_1$ (where $N_1$ is the basic clopen set $\{\alpha\!\in\! 2^\omega\mid\alpha (0)\! =\! 1\}$), and 
$\mathbb{G}^{{\bf\Gamma},a}_m\! :=\! s(\{ (0^\infty ,1\alpha )\mid\alpha\!\in\!\mathbb{C}\} )$.\smallskip

\scalebox{0.5}{\setlength{\unitlength}{0.6cm}
\begin{picture}(0,6)(-21,-6)
   
   \put(2,-5.6){\line(1,0){2}}
   \put(0.1,-4){\line(0,1){2}}
   \put(2.85,-7.5){$\mathbb{C}$}
   \put(4.2,-7.5){$\neg\mathbb{C}$}
   \multiput(2,-7.6)(0,0.2){10}{\circle*{0.05}}
   \multiput(4,-7.6)(0,0.2){10}{\circle*{0.05}}
   \multiput(-1.9,-2)(0.2,0){10}{\circle*{0.05}}
   \multiput(-1.9,-4)(0.2,0){10}{\circle*{0.05}}
   \put(2,-7.8){\line(1,0){3}}
   \put(2,-8){\line(0,1){0.4}}
   \put(4,-8){\line(0,1){0.4}}
   \put(5,-8){\line(0,1){0.4}}
   \put(0,-8.5){$0^\infty$}
   \put(3.4,-8.5){1}
   \put(-1.5,-3.15){$\mathbb{C}$}
   \put(-1.8,-1.65){$\neg\mathbb{C}$}
   \put(-2.1,-5.6){$\bullet$}
   \put(0,-7.8){$\bullet$}
   \put(-2,-4){\line(0,1){3}}
   \put(-2.2,-4){\line(1,0){0.4}}
   \put(-2.2,-2){\line(1,0){0.4}}
   \put(-2.2,-1){\line(1,0){0.4}}
   \put(-2.8,-5.6){$0^\infty$}
   \put(-2.7,-2.6){1}
   \put(1.85,-9.5){$\mathbb{G}^{{\bf\Gamma},a}_m$}
     
\end{picture}}\medskip\medskip\medskip\medskip\medskip\medskip\medskip

\begin{thm} \label{Pi02ac} (1) There is a concrete 52 elements $\sqsubseteq_c$ and $\leq_c$-antichain basis for the class of non-$\bormtwo$ s-acyclic locally countable Borel relations on an analytic space.\smallskip

\noindent (2) $\{ (\mathbb{D},\mathbb{G}^{\bormtwo}_m),(\mathbb{S},\mathbb{G}^{{\bormtwo},a}_m)\}$ is a 
$\sqsubseteq_c$ and $\leq_c$-antichain basis for the class of non-$\bormtwo$ acyclic locally countable Borel graphs on an analytic space.\end{thm}

 Theorems \ref{Sigma02ac} and \ref{Pi02ac} are consequences of the following, not involving local countability.
 
\begin{thm} \label{2ac} (1) Let ${\bf\Gamma}$ be a non self-dual Borel class of rank two.\smallskip

\noindent (1) There is a concrete 76 elements $\sqsubseteq_c$ and $\leq_c$-antichain basis for the class of non-$\bf\Gamma$ s-acyclic Borel relations on an analytic space.\smallskip

\noindent (2) $\{ (\mathbb{D},\mathbb{G}^{\bf\Gamma}_m),(\mathbb{S},\mathbb{G}^{{\bf\Gamma},a}_m)\}$ is a 
$\sqsubseteq_c$ and $\leq_c$-antichain basis for the class of non-$\bf\Gamma$ acyclic Borel graphs on an analytic space.\end{thm}

 This result can be extended, with a suitable acyclicity assumption. In [L4], it is shown that the containment in a s-acyclic 
$\boratwo$ relation allows some positive reducibility results (see, for example Theorem 4.1 in [L4]). A natural way to ensure this is to have a s-acyclic closure (recall Theorem 1.9 in [L4]). In sequential spaces like the Baire space $\omega^\omega$, having s-acyclic levels is sufficient to ensure this (see Proposition 2.7 in [L4]). Moreover, there is a s-acyclic closed relation on $2^\omega$ containing Borel relations of arbitrarily high complexity (even potential complexity), by Proposition 3.17 in [L4]. The next result unifies the classes of rank at least two. 

\begin{thm} \label{genac} Let ${\bf\Gamma}$ be a non self-dual Borel class of rank at least two.\smallskip

\noindent (1) There is a concrete 76 elements $\sqsubseteq_c$ and $\leq_c$-antichain basis for the class of non-$\bf\Gamma$ Borel relations on an analytic space contained in a s-acyclic $\boratwo$ relation.\smallskip

\noindent (2) $\{ (\mathbb{D},\mathbb{G}^{\bf\Gamma}_m),(\mathbb{S},\mathbb{G}^{{\bf\Gamma},a}_m)\}$ is a $\sqsubseteq_c$ and $\leq_c$-antichain basis for the class of 
non-$\bf\Gamma$ Borel graphs on an analytic space contained in an acyclic $\boratwo$ graph.\end{thm}

\noindent\bf Questions.\rm\ For the classes of rank at least three, we gave finite antichain basis for small relations. Is there an antichain basis if we do not assume smallness? If yes, is it finite? Countable? Is it true that any basis must have size continuum? The graph 
$\big( 2^\omega ,\mathbb{C}^2\!\setminus\!\Delta (2^\omega )\big)$ shows that we cannot simply erase the acyclicity assumptions in Theorems \ref{2ac}.(2) and \ref{genac}.(2).\medskip 

 Theorems \ref{morethree}, \ref{two}-\ref{baum}, \ref{Sigma02ac}-\ref{genac} are proved in Sections 2, 3, 4 respectively. In Section 6, we close this study of $\sqsubseteq_c$ by providing an antichain basis for the class of uncountable analytic relations on a Hausdorff topological space, which gives a perfect set theorem for binary relations. We extend the notation $\mathbb{G}^{\bf\Gamma}_m,\mathbb{G}^{{\bf\Gamma},a}_m$ to the class ${\bf\Gamma}\! =\!\{\emptyset\}$.
  
\begin{thm} \label{cou} (1) There is a concrete 13 elements $\sqsubseteq_c$-antichain basis for the class of uncountable analytic relations on a Hausdorff topological space.\smallskip

\noindent (2) $\{ (\mathbb{D},\mathbb{G}^{\{\emptyset\}}_m),({\mathbb S},\mathbb{G}^{\{\emptyset\} ,a}_m),(2^\omega ,\not= )\}$ is a $\sqsubseteq_c$-antichain basis for the class of uncountable analytic graphs on a Hausdorff topological space.\end{thm}

 Note that $(2^\omega ,\not= )$ is not acyclic, so that we recover the basis met in Theorems \ref{2ac}.(2) and \ref{genac}.(2) in the acyclic case. Also, $({\mathbb S},\mathbb{G}^{\{\emptyset\} ,a}_m)$ and $(2^\omega ,\not= )$ are not locally countable. In conclusion, [L5] and the present study show that, when our finite antichain basis exist, they are small in the cases of equivalence relations and graphs.\medskip
 
 We saw that there is no antichain basis in the class of non-potentially closed Borel relations on a Polish space. However, it follows from the main results in [L2] and [L3] that this problem can be fixed if we allow partial reductions, on a closed relation (which in fact is suitable for any non self-dual Borel class).
 
\vfill\eject 
 
 This solution involves the following quasi-order, less considered than $\sqsubseteq_c$ and $\leq_c$. Let $X$, $Y$ be topological spaces, and $A_0,A_1\!\subseteq\! X^2$ (resp., $B_0,B_1\!\subseteq\! Y^2$) be disjoint. Then we set 
$$(X, A_0, A_1)\leq (Y, B_0, B_1)\Leftrightarrow
\exists f\! :\! X\!\rightarrow\! Y\mbox{ continuous with }\forall\varepsilon\!\in\! 2\ \ A_\varepsilon\!\subseteq\! (f\!\times\! f)^{-1}(B_\varepsilon ).$$
A similar result holds here, for the Borel classes instead of the potential Borel classes. We define  
$\mathbb{O}^{\bf\Gamma}_m\! :=\!\big\{\big( (0,x),(1,x)\big)\mid x\!\in\!\mathbb{C}\big\}$, so that $\mathbb{G}^{\bf\Gamma}_m\! =\! s(\mathbb{O}^{\bf\Gamma}_m)$. 
 
\begin{thm} \label{partial} Let ${\bf\Gamma}$ be a non self-dual Borel class.\smallskip

\noindent (1) Let $X$ be an analytic space, and $R$ be a Borel relation on $X$. Exactly one of the following holds:\smallskip

(a) the relation $R$ is a $\bf\Gamma$ subset of $X^2$,\smallskip

(b) $(\mathbb{D},\mathbb{O}^{\bf\Gamma}_m,\mathbb{O}^{\{\emptyset\}}_m\!\setminus\!\mathbb{O}^{\bf\Gamma}_m)\leq (X,R,X^2\!\setminus\! R)$.\smallskip

\noindent (2) A similar statement holds for graphs, with 
$(\mathbb{G}^{\bf\Gamma}_m,\mathbb{G}^{\{\emptyset\}}_m\!\setminus\!\mathbb{G}^{\bf\Gamma}_m)$ instead of 
$(\mathbb{O}^{\bf\Gamma}_m,\mathbb{O}^{\{\emptyset\}}_m\!\setminus\!\mathbb{O}^{\bf\Gamma}_m)$.\end{thm}

 This last result can be extended to any non self-dual Wadge class of Borel sets. Note that there is no injectivity in Theorem \ref{partial}, because of the examples $\mathbb{G}^{\bf\Gamma}_m$ and $\mathbb{G}^{\bf\Gamma ,a}_m$ for instance. 
  
\section{$\!\!\!\!\!\!$ The general case}

\label{SectionGeneral}

\subsection{$\!\!\!\!\!\!$ Preliminary results}\indent 

 We first extend Lemma 4.1 in [L5].

\begin{lem} \label{basicr} Let $\bf\Gamma$ be a class of sets closed under continuous pre-images, 
$Y,Z$ be topological spaces, and $R,S$ be a relation on $Y,Z$ respectively.\smallskip

(a) If $R$ is in $\bf\Gamma$, then the sections of $R$ are also in $\bf\Gamma$.\smallskip

(b) If $S$ has vertical (resp., horizontal) sections in $\bf\Gamma$ and $(Y,R)\leq_c(Z,S)$, then the vertical (resp., horizontal) sections of $R$ are also in $\bf\Gamma$.\smallskip

(c) If $S$ has countable vertical (resp., horizontal) sections and $(Y,R)\sqsubseteq_c(Z,S)$, then the vertical (resp., horizontal) sections of $R$ are also countable.\end{lem}

\noindent\bf Proof.\rm\ (a) comes from the fact that if $y\!\in\! Y$, then the maps $i_y\! :\! y'\!\mapsto\! (y,y')$, 
$j_y\! :\! y'\!\mapsto\! (y',y)$ are continuous and satisfy $R_y\! =\! i_y^{-1}(R)$, $R^y\! =\! j_y^{-1}(R)$. The statements (b), (c) come from the facts that $R_y\! =\! f^{-1}(S_{f(y)})$, $R^y\! =\! f^{-1}(S^{f(y)})$.\hfill{$\square$}\medskip

 We now extend Theorem 4.3 in [L5].\medskip
 
\noindent\bf Notation.\rm\ Let $R$ be a relation on $\mathbb{D}$. We set, for $\varepsilon ,\eta\!\in\! 2$, 
$R_{\varepsilon ,\eta}\! :=\!
\big\{ (\alpha ,\beta )\!\in\!\mathbb{K}^2\mid\big( (\varepsilon ,\alpha ),(\eta ,\beta )\big)\!\in\! R\big\}$.

\begin{thm} \label{contBr} Let ${\bf\Gamma}$ be a non self-dual Borel class of rank at least two, $(2^\omega ,\mathbb{C})$ be 
$\bf\Gamma$-good, $X$ be an analytic space, and $R$ be a Borel relation on $X$. Exactly one of the following holds:\smallskip  

(a) the relation $R$ is a $\bf\Gamma$ subset of $X^2$,\smallskip  

(b) one of the following holds:\smallskip

\noindent (1) the relation $R$ has at least one section not in $\bf\Gamma$,\smallskip

\noindent (2) there is a relation $\mathbb{R}$ on $2^\omega$ such that 
$\mathbb{R}\cap\Delta (2^\omega )\! =\!\Delta (\mathbb{C})$ and $(2^\omega ,\mathbb{R})\sqsubseteq_c(X,R)$,\smallskip

\noindent (3) there is a relation $\mathbb{R}$ on $\mathbb{D}$ such that 
$\mathbb{R}_{0,1}\cap\Delta (2^\omega )\! =\!\Delta (\mathbb{C})$ and $(\mathbb{D},\mathbb{R})\sqsubseteq_c(X,R)$.\end{thm}

\noindent\bf Proof.\rm\ We first note that (a) and (b) cannot hold simultaneously. Indeed, we argue by contradiction, so that $R$ has sections in $\bf\Gamma$ by Lemma \ref{basicr}.(a), $\mathbb{R}\!\in\! {\bf\Gamma}\big( (2^\omega )^2\big)$, 
$\mathbb{R}\!\in\! {\bf\Gamma}(\mathbb{D}^2)$, and 
$\mathbb{R}\cap\Delta (2^\omega )\!\in\! {\bf\Gamma}\big(\Delta (2^\omega )\big)$, 
$\mathbb{R}_{0,1}\cap\Delta (2^\omega )\!\in\! {\bf\Gamma}\big(\Delta (2^\omega )\big)$ respectively. This implies that 
$\mathbb{C}\!\in\!\ {\bf\Gamma}(2^\omega )$, which is absurd.\medskip

 Assume now that (a) and (b).(1) do not hold. Theorem 1.9 in [L5] gives $f\! :=\! (f_0,f_1)\! :\! 2^\omega\!\rightarrow\! X^2$ continuous with injective coordinates such that $\mathbb{C}\! =\! f^{-1}(R)$.\medskip
 
\noindent\bf Case 1.\rm\ $f[\neg\mathbb{C}]\!\subseteq\!\Delta (X)$.\medskip

 Note that $f[2^\omega ]\!\subseteq\!\Delta (X)$, by the choice of $\mathbb{C}$, so that $f_0\! =\! f_1$. We set 
$\mathbb{R}\! :=\! (f_0\!\times \! f_0)^{-1}(R)$, so that $\mathbb{R}\cap\Delta (2^\omega )\! =\!\Delta (\mathbb{C})$. Note that $(2^\omega ,\mathbb{R})\sqsubseteq_c(X,R)$, with witness $f_0$.\medskip

\noindent\bf Case 2.\rm\ $f[\neg\mathbb{C}]\!\not\subseteq\!\Delta (X)$.\medskip

 We may assume that $f_0$ and $f_1$ have disjoint ranges, by the choice of 
$\mathbb{C}$. We define $g\! :\!\mathbb{D}\!\rightarrow\! X$ by $g(\varepsilon ,\alpha )\! :=\! f_\varepsilon (\alpha )$. Note that $g$ is injective continuous, 
$\big\{\big( (0,\alpha ),(1,\alpha )\big)\mid\alpha\!\in\!\mathbb{C}\big\}\!\subseteq\! (g\!\times\! g)^{-1}(R)$ and 
$\big\{\big( (0,\alpha ),(1,\alpha )\big)\mid\alpha\!\notin\!\mathbb{C}\big\}\!\subseteq\! (g\!\times\! g)^{-1}(\neg R)$. It remains to set $\mathbb{R}\! :=\! (g\!\times\! g)^{-1}(R)$.\hfill{$\square$}\medskip

 The following property is crucial in the sequel, as well as in [L5].

\begin{defi} Let $f\! :\!\mathbb{K}\!\rightarrow\!\mathbb{K}$ be a function, and $\mathbb{C}\!\subseteq\!\mathbb{K}$. We say that $f$ {\bf preserves} $\mathbb{C}$ if $\mathbb{C}\! =\! f^{-1}(\mathbb{C})$.\end{defi}

 It is strongly related to condition (2) in Theorem \ref{contBr}. 

\begin{lem} \label{pre} Let $f\! :\!\mathbb{K}\!\rightarrow\!\mathbb{K}$ be a function, $\mathbb{C}\!\subseteq\!\mathbb{K}$, and $R$ be a relation on $\mathbb{K}$ with $R\cap\Delta (\mathbb{K})\! =\!\Delta (\mathbb{C})$. Then the following are equivalent:\smallskip

(1) $f$ preserves $\mathbb{C}$.\smallskip

(2) $(f\!\times\! f)^{-1}(R)\cap\Delta (\mathbb{K})\! =\!\Delta (\mathbb{C})$.\end{lem}

\noindent\bf Proof.\rm\ We just apply the definitions.\hfill{$\square$}\medskip

 Finally, we will make a strong use of the following result (see page 433 in [Lo-SR1]).

\begin{thm} \label{Lo-SR} (Louveau-Saint Raymond) Let $\xi\!\geq\! 1$ be a countable ordinal, $(\mathbb{K},\mathbb{C})$ be $\boraxi$-good (or simply $\mathbb{C}\!\in\!\bormxi (2^\omega )\!\setminus\!\boraxi (2^\omega )$ if $\xi\!\geq\! 3$), $X$ be a Polish space, and $A,B$ be disjoint analytic subsets of $X$. Then exactly one of the following holds:\smallskip

(a) the set $A$ is separable from $B$ by a $\boraxi$ set,\smallskip

(b) we can find $f\! :\!\mathbb{K}\!\rightarrow\! X$ injective continuous such that 
$\mathbb{C}\!\subseteq\! f^{-1}(A)$ and $\neg\mathbb{C}\!\subseteq\! f^{-1}(B)$.\end{thm}

 A first consequence of this is Theorem \ref{partial} (these two results extend to the non self-dual Wadge classes of Borel sets, see Theorem 5.2 in [Lo-SR2]).\medskip
 
\noindent\bf Proof of Theorem \ref{partial}.\rm\ (1) As $\mathbb{C}\!\notin\! {\bf\Gamma}$, $\mathbb{O}^{\bf\Gamma}_m$ is not separable from $\mathbb{O}^{\{\emptyset\}}_m\!\setminus\!\mathbb{O}^{\bf\Gamma}_m$ by a set in $\bf\Gamma$, and (a), (b) cannot hold simultaneously. So assume that (a) does not hold. As $X$ is separable metrizable, we may assume that $X$ is a subset of the Polish space $[0,1]^\omega$ (see 4.14 in [K1]). Note that the analytic set $R$ is not separable from the analytic set $X^2\!\setminus\! R$ by a $\bf\Gamma$ subset of $([0,1]^\omega )^2$. Theorem \ref{Lo-SR} provides 
$h\! :\!\mathbb{K}\!\rightarrow\! ([0,1]^\omega )^2$ continuous such that $\mathbb{C}\!\subseteq\! h^{-1}(R)$ and 
$\neg\mathbb{C}\!\subseteq\! h^{-1}(X^2\!\setminus\! R)$. Note that $h$ takes values in $X^2$. We define 
$f\! :\!\mathbb{D}\!\rightarrow\! X$ by $f(\varepsilon ,x)\! :=\! h_\varepsilon (x)$, so that $f$ is continuous. Moreover,
$$\big( (0,x),(1,x)\big)\!\in\!\mathbb{O}^{\bf\Gamma}_m\Leftrightarrow x\!\in\!\mathbb{C}\Leftrightarrow h(x)\! =\!\big( h_0(x),h_1(x)\big)\!\in\! R\Leftrightarrow
\big( f(0,x),f(1,x)\big)\!\in\! R\mbox{,}$$
so that (b) holds.\medskip

\noindent (2) We just have to consider the symmetrizations of the relations appearing in (1).\hfill{$\square$}

\subsection{$\!\!\!\!\!\!$ Simplifications for the rank two}\indent 

 We will see that, for classes of rank two, the basic examples are contained in $\Delta (2^\omega )\cup\mathbb{P}_f^2$. The next proof is the first of our two proofs using effective descriptive set theory (see also Theorem \ref{morethreeacy}).
 
\begin{lem} \label{bot} Let $R$ be a Borel relation on $\mathbb{P}_\infty$ whose sections are separable from $\mathbb{P}_f$ by a 
$\boratwo$ set. Then we can find a sequence $(R_n)_{n\in\omega}$ of relations closed in $\mathbb{P}_\infty\!\times\! 2^\omega$ and 
$2^\omega\!\times\!\mathbb{P}_\infty$, as well as $f\! :\! 2^\omega\!\rightarrow\! 2^\omega$ injective continuous preserving 
$\mathbb{P}_f$ such that $(f\!\times\! f)^{-1}(R)\!\subseteq\!\bigcup_{n\in\omega}~R_n$.\end{lem}

\noindent\bf Proof.\rm\ In order to simplify the notation, we assume that $R$ is a $\Borel$ relation on $2^\omega$. Recall that we can find $\Ca$ sets $W\!\subseteq\! 2^\omega\!\times\!\omega$ and $C\!\subseteq\!2^\omega\!\times\!\omega\!\times\! 2^\omega$ such that 
$\Borel (\alpha )(2^\omega )\! =\!\{ C_{\alpha ,n}\mid (\alpha ,n)\!\in\! W\}$ for each $\alpha\!\in\! 2^\omega$ and 
$\{ (\alpha ,n,\beta )\!\in\! 2^\omega\!\times\!\omega\!\times\! 2^\omega\mid (\alpha ,n)\!\in\! W~\wedge ~(\alpha ,n,\beta )\!\notin\! C\}$ 
is a $\Ca$ subset of $2^\omega\!\times\!\omega\!\times\! 2^\omega$ (see Section 2 in [Lo1]). Intuitively, $W_\alpha$ is the set of codes for the $\Borel (\alpha )$ subsets of $2^\omega$. We set 
$W_2\! :=\!\{ (\alpha ,n)\!\in\! W\mid C_{\alpha ,n}\mbox{ is a }\bormtwo\cap\Borel (\alpha )\mbox{ subset of }2^\omega\}$. Intuitively, 
$(W_2)_\alpha$ is the set of codes for the $\bormtwo\cap\Borel (\alpha )$ subsets of $2^\omega$. By Section 2 in [Lo1], the set $W_2$ is $\Ca$. We set 
$$P\! :=\{ (\alpha ,n)\!\in\! 2^\omega\!\times\!\omega\mid 
(\alpha ,n)\!\in\! W_2\wedge R_\alpha\!\subseteq\!\neg C_{\alpha ,n}\!\subseteq\!\mathbb{P}_\infty\} .$$ 
Note that $P$ is $\Ca$. Moreover, for each $\alpha\!\in\!\mathbb{P}_\infty$, there is $n\!\in\!\omega$ such that $(\alpha ,n)\!\in\! P$, by Theorem 2.B' in [Lo1]. The $\it\Delta$-selection principle provides a $\Borel$-recursive map $f\! :\! 2^\omega\!\rightarrow\!\omega$ such that $\big(\alpha ,f(\alpha )\big)\!\in\! P$ if $\alpha\!\in\!\mathbb{P}_\infty$ (see 4B.5 in [Mo]). We set 
$B\! :=\!\big\{ (\alpha ,\beta )\!\in\!\mathbb{P}_\infty\!\times\! 2^\omega\mid\big(\alpha ,f(\alpha ),\beta \big)\!\notin\! C\big\}$. Note that $B$ is a $\Borel$ set with vertical sections in $\boratwo$, and $R_\alpha\!\subseteq\! B_\alpha\!\subseteq\!\mathbb{P}_\infty$ for each 
$\alpha\!\in\!\mathbb{P}_\infty$. Theorem 3.6 in [Lo1] provides a finer Polish topology $\tau$ on $2^\omega$ such that 
$B\!\in\!\boratwo\big( (2^\omega ,\tau)\!\times\! 2^\omega\big)$. Note that the identity map from $(2^\omega ,\tau )$ into $2^\omega$ is a continuous bijection. By 15.2 in [K1], it is a Borel isomorphism. By 11.5 in [K1], its inverse is Baire measurable. By 8.38 in [K1], there is a dense $G_\delta$ subset $G$ of $2^\omega$ on which $\tau$ coincides with the usual topology on $2^\omega$. In particular, 
$B\cap (G\!\times\! 2^\omega )\!\in\!\boratwo (G\!\times\! 2^\omega )$ and we may assume that $G\!\subseteq\!\mathbb{P}_\infty$. Note that $G$ is not separable from $\mathbb{P}_f$ by a set in $\boratwo$, by Baire's theorem. Theorem \ref{Lo-SR} provides 
$g\! :\! 2^\omega\!\rightarrow\! 2^\omega$ injective continuous such that $\mathbb{P}_\infty\!\subseteq\! g^{-1}(G)$ and 
$\mathbb{P}_f\!\subseteq\! g^{-1} (\mathbb{P}_f)$. Note that $(g\!\times\! g)^{-1}(B)$ is $\boratwo$ in 
$\mathbb{P}_\infty\!\times\! 2^\omega$ and contained in $\mathbb{P}_\infty^2$. So, replacing $B$ with $(g\!\times\! g)^{-1}(B)$ if necessary, we may assume that $R$ is contained in a Borel set $B$ which is $\boratwo$ in $\mathbb{P}_\infty\!\times\! 2^\omega$ and contained in $\mathbb{P}_\infty^2$. Similarly, we may assume that $R$ is contained in a Borel set $D$ which is $\boratwo$ in $2^\omega\!\times\!\mathbb{P}_\infty$ and contained in $\mathbb{P}_\infty^2$. Let $(B_p)_{p\in\omega},(D_q)_{q\in\omega}$ be sequences of closed relations on $2^\omega$ with ${B\! =\!\bigcup_{p\in\omega}~B_p\cap (\mathbb{P}_\infty\!\times\! 2^\omega )}$, 
$D\! =\!\bigcup_{q\in\omega}~D_q\cap (2^\omega\!\times\!\mathbb{P}_\infty )$ respectively. We set 
$R_{p,q}\! :=\! B_p\cap D_q\cap\mathbb{P}_\infty^2$, so that $R_{p,q}$ is closed in $\mathbb{P}_\infty\!\times\! 2^\omega$ and 
$2^\omega\!\times\!\mathbb{P}_\infty$ since, for example, 
$R_{p,q}\! :=\! B_p\cap D_q\cap (\mathbb{P}_\infty\!\times\! 2^\omega )\! =\! B_p\cap D_q\cap (2^\omega\!\times\!\mathbb{P}_\infty )$, and $R\!\subseteq\!\bigcup_{p,q\in\omega}~R_{p,q}$.\hfill{$\square$}\medskip

\noindent\bf Notation.\rm\ We define a well-order $\leq_l$ of order type $\omega$ on $2^{<\omega}$ by 
${s\!\leq_l\! t\Leftrightarrow\vert s\vert\! <\!\vert t\vert\vee (\vert s\vert\! =\!\vert t\vert\wedge s\!\leq_{\text{lex}}\! t)}$ and, as usual for linear orders, set $s\! <_l\! t\Leftrightarrow s\!\leq_l\! t\wedge s\!\not=\! t$. Let 
$b\! :\! (\omega ,\leq )\!\rightarrow\! (2^{<\omega},\leq_l)$ be the increasing bijection, $\alpha_{n+1}\! :=\! b(n)10^\infty$, 
$\alpha_0\! :=\! 0^\infty$, so that $\mathbb{P}_f\! =\!\{\alpha_n\mid n\!\in\!\omega\}$. We then set 
${Q\! :=\!\{\emptyset\}\cup\{ u1\mid u\!\in\! 2^{<\omega}\}}$, so that $\mathbb{P}_f\! =\!\{ t0^\infty\mid t\!\in\! Q\}$.

\begin{defi} A {\bf finitely dense Cantor set} is a copy $C$ of $2^\omega$ in $2^\omega$ such that $\mathbb{P}_f\cap C$ is dense in $C$.\end{defi}

 Note that if $C$ is a finitely dense Cantor set, then $\mathbb{P}_f\cap C$ is countable dense, and also co-dense, in $C$, which implies that $\mathbb{P}_f\cap C$ is $\boratwo$ and not $\bormtwo$ in $C$, by Baire's theorem.\medskip
 
\noindent\bf Conventions.\rm\ In the rest of Sections 2 and 3, we will perform a number of Cantor-like constructions. The following will always hold. We fix a finitely dense Cantor set $C$, and we want to construct $f\! :\! 2^\omega\!\rightarrow\! C$ injective continuous preserving $\mathbb{P}_f$. We inductively construct a sequence $(n_t)_{t\in 2^{<\omega}}$ of positive natural numbers, and a sequence $(U_t)_{t\in 2^{<\omega}}$ of basic clopen subsets of $C$, satisfying the following conditions.
$$\begin{array}{ll}
& (1)~U_{t\varepsilon}\!\subseteq\! U_t\cr
& (2)~\alpha_{n_t}\!\in\! U_t\cr
& (3)~\mbox{diam}(U_t)\!\leq\! 2^{-\vert t\vert}\cr
& (4)~U_{t0}\cap U_{t1}\! =\!\emptyset\cr
& (5)~n_{t0}\! =\! n_t\cr
& (6)~U_{t1}\cap\{\alpha_n\mid n\!\leq\!\vert t\vert\}\! =\!\emptyset
\end{array}$$
Assume that this is done. Using (1)-(3), we define $f\! :\! 2^\omega\!\rightarrow\! C$ by 
$\{ f(\beta )\}\! :=\!\bigcap_{n\in\omega}~U_{\beta\vert n}$, and $f$ is injective continuous by (4). If 
$t\!\in\! Q$ and $\alpha\! =\! t0^\infty$, then $f(\alpha )\! =\!\alpha_{n_t}$ by (5), so that 
$\mathbb{P}_f\!\subseteq\! f^{-1}(\mathbb{P}_f)$. Condition (6) ensures that 
$\mathbb{P}_\infty\!\subseteq\! f^{-1}(\mathbb{P}_\infty)$, so that $f$ preserves $\mathbb{P}_f$. For the first step of the induction, we choose $n_\emptyset\!\geq\! 1$ in such a way that $\alpha_{n_\emptyset}\!\in\! C$, and a basic clopen neighbourhood $U_\emptyset$ of $\alpha_{n_\emptyset}$. Condition (5) defines $n_{t0}$. It will also be convenient to set 
$s_t\! :=\! b(n_t\! -\! 1)1$.

\begin{lem} \label{inf} Let $(R_n)_{n\in\omega}$ be a sequence of relations on $\mathbb{P}_\infty$ which are closed in 
$\mathbb{P}_\infty\!\times\! 2^\omega$ and in $2^\omega\!\times\!\mathbb{P}_\infty$. Then there is 
$f\! :\! 2^\omega\!\rightarrow\! 2^\omega$ injective continuous preserving $\mathbb{P}_f$ with the property that 
$\big( f(\alpha ),f(\beta )\big)\!\notin\!\bigcup_{n\in\omega}~R_n$ if $\alpha\!\not=\!\beta\!\in\!\mathbb{P}_\infty$.\end{lem}

\noindent\bf Proof.\rm\ We ensure (1)-(6) with $C\! =\! 2^\omega$ and
$$(7)~(U_{s1}\!\times\! U_{t\varepsilon})\cap (\bigcup_{n\leq l}~R_n)\! =\!\emptyset\mbox{ if }s\!\not=\! t\!\in\! 2^l$$
Assume that this is done. If $\alpha\!\not=\!\beta\!\in\!\mathbb{P}_\infty$, then we can find  $l$ with $\alpha\vert l\!\not=\!\beta\vert l$, and a strictly increasing sequence 
$(l_k)_{k\in\omega}$ of natural numbers bigger than $l$ such that $\alpha (l_k)\! =\! 1$ for each $k$. Condition (7) ensures that 
$U_{\alpha\vert (l_k+1)}\!\times\! U_{\beta\vert (l_k+1)}$ does not meet $\bigcup_{n\leq l_k}~R_n$, so that $\big( f(\alpha ),f(\beta )\big)\!\notin\!\bigcup_{n\in\omega}~R_n$.\medskip
 
 So it is enough to prove that the construction is possible. We first set $n_\emptyset\! :=\! 1$ and 
$U_\emptyset\! :=\! 2^\omega$. We choose $n_1\!\geq\! 1$ such that $\alpha_{n_1}\!\not=\!\alpha_{n_0}$, and $U_0,U_1$ disjoint with diameter at most $2^{-1}$ such that 
$\alpha_{n_\varepsilon}\!\in\! U_\varepsilon\!\subseteq\! 2^\omega\!\setminus\!\{\alpha_0\}$. Assume that 
$(n_t)_{\vert t\vert\leq l}$ and $(U_t)_{\vert t\vert\leq l}$ satisfying (1)-(7) have been constructed for some $l\!\geq\! 1$, which is the case for $l\! =\! 1$. We set $F\! :=\!\{\alpha_{n_t}\mid t\!\in\! 2^l\}\cup\{\alpha_n\mid n\!\leq\! l\}$ and $L\! :=\!\bigcup_{n\leq l}~R_n$.\medskip

\noindent\bf Claim.\it\ Let $U,U_0,\cdots ,U_m$ be nonempty open subsets of $2^\omega$, 
$\gamma_i\!\in\!\mathbb{P}_f\cap U_i$, for each $i\!\leq\! m$, and $F$ be a finite subset of $\mathbb{P}_f$ containing 
$\{\gamma_i\mid i\!\leq\! m\}$. Then we can find $\gamma\!\in\!\mathbb{P}_f\cap U$ and clopen subsets 
$V,V_0,\cdots ,V_m$ of $2^\omega$ with diameter at most $2^{-l-1}$ such that 
$\gamma\!\in\! V\!\subseteq\! U\!\setminus\! (F\cup V_0)$, $\gamma_i\!\in\! V_i\!\subseteq\! U_i$ and 
$$\big( (V\!\times\! V_i)\cup (V_i\!\times\! V)\big)\cap L\! =\!\emptyset$$ 
for each $i\!\leq\! m$.\rm\medskip

 Indeed, fix $i\!\leq\! m$. Note that $\mathbb{P}_\infty\!\times\!\{\gamma_i\}\!\subseteq\!\mathbb{P}_\infty\!\times\!\mathbb{P}_f$ and $L\!\subseteq\!\mathbb{P}_\infty^2$ are disjoint closed subsets of the zero-dimensional space $\mathbb{P}_\infty\!\times\! 2^\omega$. This gives a clopen subset ${}_iC$ of $\mathbb{P}_\infty\!\times\! 2^\omega$ with 
$\mathbb{P}_\infty\!\times\!\{\gamma_i\}\!\subseteq\! {}_iC\!\subseteq\!\neg L$ (see 22.16 in [K1]). So for each $\beta\!\in\!\mathbb{P}_\infty$ we can find a clopen subset 
${}_{i,\beta}O$ of $\mathbb{P}_\infty$ and a clopen subset ${}_{i,\beta}D$ of $2^\omega$ such that $(\beta ,\gamma_i)\!\in\! {}_{i,\beta}O\!\times\! {}_{i,\beta}D\!\subseteq\! {}_iC$ and 
${}_{i,\beta}D\!\subseteq\! U_i$.

\vfill\eject

 As $\mathbb{P}_\infty\! =\!\bigcup_{\beta\in\mathbb{P}_\infty}~{}_{i,\beta}O$, we can find a sequence $(\beta_p)_{p\in\omega}$ of points of 
$\mathbb{P}_\infty$ with the property that $\mathbb{P}_\infty\! =\!\bigcup_{p\in\omega}~{}_{i,\beta_p}O$. We set 
${}_{i,p}O\! :=\! {}_{i,\beta_p}O\!\setminus\! (\bigcup_{q<p}~{}_{i,\beta_q}O)$, so that $({}_{i,p}O)_{p\in\omega}$ is a partition of $\mathbb{P}_\infty$ into clopen sets. We set 
${}_{i,p}D\! :=\! {}_{i,\beta_p}D$, so that ${}_{i,p}O\!\times\! {}_{i,p}D\!\subseteq\! {}_iC$. Let ${}_{i,p}W$ be an open subset of $2^\omega$ with 
${}_{i,p}O\! =\!\mathbb{P}_\infty\cap {}_{i,p}W$. By 22.16 in [K1], there is a sequence $({}_{i,p}U)_{p\in\omega}$ of pairwise disjoint open subsets of $2^\omega$ such that 
${}_{i,p}U\!\subseteq\! {}_{i,p}W$ and $\bigcup_{p\in\omega}~{}_{i,p}U\! =\!\bigcup_{p\in\omega}~{}_{i,p}W$. Note that ${}_{i,p}O\! =\!\mathbb{P}_\infty\cap {}_{i,p}U$.\medskip

 Similarly, let $C_i$ be a clopen subset of $2^\omega\!\times\!\mathbb{P}_\infty$ with 
$\{\gamma_i\}\!\times\!\mathbb{P}_\infty\!\subseteq\! C_i\!\subseteq\!\neg L$, $(O_{i,q})_{q\in\omega}$ be a partition of 
$\mathbb{P}_\infty$ into clopen sets, $(D_{i,q})_{q\in\omega}$ be a sequence of clopen subsets of $2^\omega$, and 
$(U_{i,q})_{q\in\omega}$ be a sequence of pairwise disjoint open subsets of $2^\omega$ such that 
$D_{i,q}\!\times\! O_{i,q}\!\subseteq\! C_i$, $\gamma_i\!\in\! D_{i,q}\!\subseteq\! U_i$ and $O_{i,q}\! =\!\mathbb{P}_\infty\cap U_{i,q}$. Now pick $\beta\!\in\!\mathbb{P}_\infty\cap U$, $p_i\!\in\!\omega$ with $\beta\!\in\! {}_{i,p_i}O$, and $q_i\!\in\!\omega$ with $\beta\!\in\! O_{i,q_i}$. We set 
$$U'\! :=\! U\cap\bigcap_{i\leq m}~({}_{i,p_i}U\cap U_{i,q_i})\!\setminus\! F.$$
As $U'$ is an open subset of $2^\omega$ containing $\beta$, we can choose $\gamma\!\in\!\mathbb{P}_f\cap U'$, and a clopen subset $V$ of $2^\omega$ with diameter at most $2^{-l-1}$ such that $\gamma\!\in\! V\!\subseteq\! U'$. If $i\!\leq\! m$, then we choose a clopen neighbourhood $V_i$ of $\gamma_i$ with diameter at most $2^{-l-1}$ such that 
$V_i\!\subseteq\! {}_{i,p_i}D\cap D_{i,q_i}$. As 
$\gamma\!\not=\!\gamma_0$, we can ensure that $V$ does not meet $V_0$. For example, note that 
$(\mathbb{P}_\infty\cap V)\!\times\! V_i\!\subseteq\! (\mathbb{P}_\infty\cap {}_{i,p_i}U)\!\times\! {}_{i,p_i}D\!\subseteq\! 
{}_iC\!\subseteq\!\neg L$ since $\mathbb{P}_\infty\cap {}_{i,p_i}U\! =\! {}_{i,p_i}O$.\hfill{$\diamond$}\medskip

 We put a linear order on the set $2^l$ of binary sequences of length $l$, so that $2^l$ is enumerated injectively by 
$\{ t_j\mid j\! <\! 2^l\}$. Let $j\! <\! 2^l$. Note that (5) defines $n_{t_j0}$. We construct, by induction on $j\! <\! 2^l$, $n_{t_j1}$, and sequences $(U^j_{t0})_{t\in 2^l}$, $(U^j_{t_i1})_{i\leq j}$ of clopen subsets of $2^\omega$ satisfying the following:
$$\begin{array}{ll}
& (a)~U^{j+1}_{t0}\!\subseteq\! U^j_{t0}\!\subseteq\! U_t\wedge U^{j+1}_{t_i1}\!\subseteq\! U^j_{t_i1}\!\subseteq\! U_{t_i}\cr
& (b)~\alpha_{n_t}\!\in\! U^j_{t0}\wedge\alpha_{n_{t_i1}}\!\in\! U^j_{t_i1}\cr
& (c)~\mbox{diam}(U^j_{t0})\mbox{, diam}(U^j_{t_i1})\!\leq\! 2^{-l-1}\cr
& (d)~U^j_{t_j0}\cap U^j_{t_j1}\! =\!\emptyset\cr
& (e)~U^j_{t_j1}\cap F\! =\!\emptyset\cr
& (f)~\Big(\big( U^j_{t_j1}\!\times\! (U^j_{t0}\cup U^j_{t_i1})\big)\cup\big( (U^j_{t0}\cup U^j_{t_i1})\!\times\! U^j_{t_j1}\big)\Big)\cap L\! =\!\emptyset\mbox{ if }i\! <\! j
\end{array}$$
In order to do this, we first apply the claim to $U\! :=\! U_{t_0}$ and a family $(\gamma_i)_{i\leq m}$ of elements of 
$\mathbb{P}_f$  enumerating $\{\alpha_{n_t}\mid t\!\in\! 2^l\}$ in such a way that $\gamma_0\! =\!\alpha_{n_{t_0}}$. The corresponding family of open sets enumerates $\{ U_t\mid t\!\in\! 2^l\}$. The claim provides $\alpha_{n_{t_01}}\!\in\!\mathbb{P}_f\cap U_{t_0}$ and clopen subsets $U^0_{t_01},U^0_{t0}$ of $2^\omega$ with diameter at most $2^{-l-1}$ such that ${\alpha_{n_{t_01}}\!\in\! U^0_{t_01}\!\subseteq\! U_{t_0}\!\setminus\! (F\cup U^0_{t_00})}$, 
$\alpha_{n_t}\!\in\! U^0_{t0}\!\subseteq\! U_t$ and 
$$\big( (U^0_{t_01}\!\times\! U^0_{t0})\cup (U^0_{t0}\!\times\! U^0_{t_01})\big)\cap L\! =\!\emptyset$$ 
for each $t\!\in\! 2^l$. Assume then that ${j\! <\! 2^l\! -\! 1}$ and $(n_{t_k1})_{k\leq j}$, $(U^k_{t0})_{t\in 2^l,k\leq j}$ and 
$(U^k_{t_i1})_{i\leq k\leq j}$ satisfying $(a)$-$(f)$ have been constructed, which is the case for $j\! =\! 0$.\medskip

 We now apply the claim to ${U\! :=\! U_{t_{j+1}}}$ and a family $(\gamma_i)_{i\leq m}$ of elements of $\mathbb{P}_f$  enumerating $\{\alpha_{n_t}\mid t\!\in\! 2^l\}\cup\{\alpha_{n_{t_i1}}\mid i\!\leq\! j\}$ in such a way that 
 $\gamma_0\! =\!\alpha_{n_{t_{j+1}}}$. The corresponding family of open sets enumerates 
${\{ U^j_{t0}\mid t\!\in\! 2^l\}\cup\{ U^j_{t_i1}\mid i\!\leq\! j\}}$.

\vfill\eject

 The claim provides ${\alpha_{n_{t_{j+1}1}}\!\in\!\mathbb{P}_f\cap U_{t_{j+1}}}$ and clopen subsets 
$U^{j+1}_{t_{j+1}1},U^{j+1}_{t0},U^{j+1}_{t_i1}$ of $2^\omega$ of diameter at most $2^{-l-1}$ with the properties that $\alpha_{n_{t_{j+1}1}}$ is in 
$U^{j+1}_{t_{j+1}1}\!\subseteq\! U_{t_{j+1}}\!\setminus\! (F\cup U^{j+1}_{t_{j+1}0})$, $\alpha_{n_t}\!\in\! U^{j+1}_{t0}\!\subseteq\! U^j_{t0}$, 
$\alpha_{n_{t_i1}}\!\in\! U^{j+1}_{t_i1}\!\subseteq\! U^j_{t_i1}$ and, when $i\!\leq\! j$,  
$$\Big(\big( U^{j+1}_{t_{j+1}1}\!\times\! (U^{j+1}_{t0}\cup U^{j+1}_{t_i1})\big)\cup
\big( (U^{j+1}_{t0}\cup U^{j+1}_{t_i1})\!\times\! U^{j+1}_{t_{j+1}1}\big)\Big)\cap L\! =\!\emptyset .$$ 
It remains to set $U_{t\varepsilon}\! :=\! U^{2^l-1}_{t\varepsilon}$.\hfill{$\square$}
 
\begin{cor} \label{redtwo} Let $R$ be a locally countable Borel relation on $2^\omega$. Then we can find 
$f\! :\! 2^\omega\!\rightarrow\! 2^\omega$ injective continuous preserving $\mathbb{P}_f$ such that 
$R'\! :=\! (f\!\times\! f)^{-1}(R)\!\subseteq\!\Delta (2^\omega )\cup\mathbb{P}_f^2$.\end{cor}

\noindent\bf Proof.\rm\ As $\mathbb{P}_f$ is countable and $R$ is locally countable, the set 
$C\! :=\!\bigcup_{\alpha\in\mathbb{P}_f}~(R_\alpha\cup R^\alpha )$ is countable. We set $G\! :=\!\mathbb{P}_\infty\!\setminus\! C$, so that $G\!\subseteq\!\mathbb{P}_\infty$ is a non-meager subset of $2^\omega$ having the Baire property. Lemma 7.2 in [L5] provides 
$f\! :\! 2^\omega\!\rightarrow\! 2^\omega$ injective continuous such that $f[\mathbb{P}_\infty ]\!\subseteq\! G$ and 
$f[\mathbb{P}_f]\!\subseteq\!\mathbb{P}_f$. This proves that we may assume that 
$R\cap\big( (\mathbb{P}_\infty\!\times\!\mathbb{P}_f)\cup (\mathbb{P}_f\!\times\!\mathbb{P}_\infty )\big)\! =\!\emptyset$. By Lemma \ref{bot} applied to $R\cap\mathbb{P}_\infty^2$, we may assume that there is a sequence 
$(R_n)_{n\in\omega}$ of relations closed in $\mathbb{P}_\infty\!\times\! 2^\omega$ and $2^\omega\!\times\!\mathbb{P}_\infty$ such that 
$R\cap\mathbb{P}_\infty^2\!\subseteq\!\bigcup_{n\in\omega}~R_n$. It remains to apply Lemma \ref{inf}.\hfill{$\square$}\medskip
 
 Corollary \ref{redtwo} leads to the following.

\begin{defi} Let $\bf\Gamma$ be a non self-dual Borel class of rank two, and $(2^\omega ,\mathbb{C})$ be $\bf\Gamma$-good. A relation $\mathcal{R}$ on $2^\omega$ is {\bf diagonally complex} if it it satisfies the following:\smallskip

(1) $\mathcal{R}\cap\Delta (2^\omega )\! =\!\Delta (\mathbb{C})$,\smallskip

(2) $\mathcal{R}\!\subseteq\!\Delta (\mathbb{C})\cup\mathbb{P}_f^2$.\end{defi}

 Note that a diagonally complex relation is not in $\bf\Gamma$ by (1), and locally countable Borel by (2).\medskip

\noindent\bf Notation.\rm\ We set, for any digraph $\mathcal{D}$ on $Q$, 
$\mathcal{R}_\mathcal{D}\! :=\!\Delta (\mathbb{C})\cup\{ (s0^\infty ,t0^\infty )\mid (s,t)\!\in\!\mathcal{D}\}$. Note that any  diagonally complex relation is of the form $\mathcal{R}_\mathcal{D}$, for some digraph $\mathcal{D}$ on $Q$.\medskip

 In our future Cantor-like constructions, the definition of $n_{t1}$ will be by induction on $\leq_l$, except where indicated. Corollary \ref{redtwo} simplifies the locally countable Borel relations on $2^\omega$. Some further simplification is possible when the sections are nowhere dense.

\begin{lem} \label{nd} Let $R$ be a relation on $2^\omega$ with nowhere dense sections. Then there is 
$f\! :\! 2^\omega\!\rightarrow\! 2^\omega$ injective continuous preserving $\mathbb{P}_f$ such that  
$\big( f(\alpha ),f(\beta )\big)\!\notin\! R$ if $\alpha\!\not=\!\beta\!\in\!\mathbb{P}_f$.\end{lem}

\noindent\bf Proof.\rm\ We ensure (1)-(6) with $C\! =\! 2^\omega$ and
$$(7)~(\alpha_{n_s},\alpha_{n_t}),(\alpha_{n_t},\alpha_{n_s})\!\notin\! R\mbox{ if }t\!\in\! Q\wedge s\!\in\! 2^{\vert t\vert}\wedge 
s\! <_{\text{lex}}\! t$$
Assume that this is done. If $\alpha\!\not=\!\beta\!\in\!\mathbb{P}_f$, then for example 
$\alpha\! <_{\text{lex}}\!\beta$ and there are initial segments $s,t$ of $\alpha ,\beta$ respectively satisfying the assumption in (7). Condition (7) ensures that $\big( f(\alpha ),f(\beta )\big)\!\notin\! R$.\medskip

 So it is enough to prove that the construction is possible. We first set $n_\emptyset\! :=\! 1$ and 
$U_\emptyset\! :=\! 2^\omega$. Assume that $(n_t)_{\vert t\vert\leq l}$ and $(U_t)_{\vert t\vert\leq l}$ satisfying (1)-(7) have been constructed, which is the case for $l\! =\! 0$. Fix $t\!\in\! 2^l$. We choose $n_{t1}\!\geq\! 1$ such that 
$$\alpha_{n_{t1}}\!\in\! U_t\!\setminus\! (\{\alpha_{n_t}\}\cup\{\alpha_n\mid n\!\leq\! l\}\cup
\bigcup_{s\in 2^{l+1},s<_{\text{lex}}t1}~\overline{R_{\alpha_{n_s}}}\cup\overline{R^{\alpha_{n_s}}})\mbox{,}$$ 
which exists since $R$ has nowhere dense sections. We then choose a clopen neighbourhood with small diameter 
$U_{t\varepsilon}$ of $\alpha_{n_{t\varepsilon}}$ contained in $U_t$, ensuring (4) and (6).\hfill{$\square$}

\begin{cor} \label{dc} Let $\mathcal{R}$ be a diagonally complex relation.\smallskip

(1) If $\mathcal{R}$ has nowhere dense sections, then 
$(2^\omega ,\mathcal{R}_\emptyset )\sqsubseteq_c(2^\omega ,\mathcal{R})$.\smallskip

(2) If $\mathbb{P}_f^2\!\setminus\!\mathcal{R}$ has nowhere dense sections, then 
$(2^\omega ,\mathcal{R}_{\not=})\sqsubseteq_c(2^\omega ,\mathcal{R})$.\end{cor}

\noindent\bf Proof.\rm\ Lemma \ref{nd} gives $f\! :\! 2^\omega\!\rightarrow\! 2^\omega$ injective continuous preserving $\mathbb{P}_f$ such that $\big( f(\alpha ),f(\beta )\big)$ is not in $\mathcal{R}$ (resp., in $\mathcal{R}$) if 
$\alpha\!\not=\!\beta\!\in\!\mathbb{P}_f$. Note that $(f\!\times\! f)^{-1}(\mathcal{R})\! =\!\mathcal{R}_\emptyset$ (resp., 
$(f\!\times\! f)^{-1}(\mathcal{R})\! =\!\mathcal{R}_{\not=}$).\hfill{$\square$}\medskip

 Theorem \ref{morethree} provides a basis. Theorem \ref{twocases} to come provides another one, and is a consequence of Corollary \ref{redtwo}.\medskip

\noindent\bf Notation.\rm\ Let ${\bf\Gamma}$ be a non self-dual Borel class of rank at least two, and 
$(2^\omega ,\mathbb{C})$ be $\bf\Gamma$-good. We set $S_0\! :=\!\mathbb{C}$, $S_1\! :=\!\emptyset$, 
$S_2\! :=\! 2^\omega\!\setminus\!\mathbb{C}$, $S_3\! :=\! 2^\omega$. The next result motivates the introduction of these sets.

\begin{lem} \label{four} Let ${\bf\Gamma}$ be a non self-dual Borel class of rank at least two, $(2^\omega ,\mathbb{C})$ be 
$\bf\Gamma$-good, and $B$ be a Borel subset of $2^\omega$. Then we can find $j\!\in\! 4$ and 
$f\! :\! 2^\omega\!\rightarrow\! 2^\omega$ injective continuous preserving $\mathbb{C}$ such that $f^{-1}(B)\! =\! S_j$.\end{lem}

\noindent\bf Proof.\rm\ Note that since $\mathbb{C}\!\in\!\check{\bf\Gamma}(2^\omega )\!\setminus\! {\bf\Gamma}(2^\omega )$, either 
$\mathbb{C}\!\setminus\! B$ is not separable from $\neg\mathbb{C}$ by a set in $\bf\Gamma$, or $\mathbb{C}\cap B$ is not separable from $\neg\mathbb{C}$ by a set in $\bf\Gamma$, because $\bf\Gamma$ is closed under finite unions. Assume, for example, that the first case occurs. Similarly, either $\mathbb{C}\!\setminus\! B$ is not separable from $(\neg\mathbb{C})\cap (\neg B)$ by a set in 
$\bf\Gamma$, or $\mathbb{C}\!\setminus\! B$ is not separable from $(\neg\mathbb{C})\cap B$ by a set in $\bf\Gamma$, because 
$\bf\Gamma$ is closed under finite intersections. This shows that one of the following cases occurs:\medskip

\noindent - $\mathbb{C}\!\setminus\! B$ is not separable from $(\neg\mathbb{C})\cap (\neg B)$ by a set in $\bf\Gamma$,\smallskip

\noindent - $\mathbb{C}\!\setminus\! B$ is not separable from $(\neg\mathbb{C})\cap B$ by a set in $\bf\Gamma$,\smallskip

\noindent - $\mathbb{C}\cap B$ is not separable from $(\neg\mathbb{C})\cap (\neg B)$ by a set in $\bf\Gamma$,\smallskip

\noindent - $\mathbb{C}\cap B$ is not separable from $(\neg\mathbb{C})\cap B$ by a set in $\bf\Gamma$.\medskip

 Assume, for example, that we are in the first of these four cases. By Theorem \ref{Lo-SR}, there is 
$f\! :\! 2^\omega\!\rightarrow\! 2^\omega$ injective continuous such that $\mathbb{C}\!\subseteq\! f^{-1}(\mathbb{C}\!\setminus\! B)$ and 
$\neg\mathbb{C}\!\subseteq\! f^{-1}\big( (\neg\mathbb{C})\cap (\neg B)\big)$. In this case, $f^{-1}(B)\! =\!\emptyset\! =\! S_1$. In the other cases, $f^{-1}(B)\! =\! S_2,S_0,S_3$, respectively. So we proved that $f^{-1}(B)\! =\! S_j$ for some $j\!\in\! 4$.\hfill{$\square$}

\begin{cor} \label{Delta} ${\bf\Gamma}$ be a non self-dual Borel class of rank at least two, $(2^\omega ,\mathbb{C})$ be 
$\bf\Gamma$-good, and $R$ be a Borel relation on $\mathbb{D}$. Then there is 
$f\! :\! 2^\omega\!\rightarrow\! 2^\omega$ injective continuous preserving $\mathbb{C}$ such that 
$(f\!\times\! f)^{-1}(R_{\varepsilon ,\eta})\cap\Delta (2^\omega )\!\in\!\{\Delta (S_j)\mid j\!\leq\! 3\}$ for each 
$\varepsilon ,\eta\!\in\! 2$.\end{cor}

\noindent\bf Proof.\rm\ We set, for $\varepsilon ,\eta\!\in\! 2$, 
$E_{\varepsilon ,\eta}\! :=\!\{\alpha\!\in\! 2^\omega\mid (\alpha ,\alpha )\!\in\!\mathbb{R}_{\varepsilon ,\eta}\}$, so that 
$E_{\varepsilon ,\eta}$ is a Borel subset of $2^\omega$ and 
$R_{\varepsilon ,\eta}\cap\Delta (2^\omega )\! =\!\Delta (E_{\varepsilon ,\eta})$. Now fix $\varepsilon ,\eta\!\in\! 2$. Lemma \ref{four} provides $j\!\in\! 4$ and $g\! :\! 2^\omega\!\rightarrow\! 2^\omega$ injective continuous preserving $\mathbb{C}$ such that $(g\!\times\! g)^{-1}(E_{\varepsilon ,\eta})\! =\! S_j$. We just have to apply this for each $\varepsilon ,\eta\!\in\! 2$.\hfill{$\square$}

\begin{thm} \label{twocases} Let ${\bf\Gamma}$ be a non self-dual Borel class of rank two, 
$(2^\omega ,\mathbb{C})$ be $\bf\Gamma$-good, $X$ be an analytic space, and $R$ be a locally countable Borel relation on $X$ whose sections are in $\bf\Gamma$. Exactly one of the following holds.\smallskip

(a) the relation $R$ is a $\bf\Gamma$ subset of $X^2$,\smallskip

(b) one of the following holds:\smallskip

\noindent (1) there is a diagonally complex relation $\mathbb{R}$ on $2^\omega$ such that 
$(2^\omega ,\mathbb{R})\sqsubseteq_c(X,R)$,\smallskip

\noindent (2) there is a relation $\mathbb{R}$ on $\mathbb{D}$ such that, for each $\varepsilon ,\eta\!\in\! 2$, if 
$\mathbb{R}_{\varepsilon ,\eta}\cap\Delta (2^\omega )\! =\!\Delta (E_{\varepsilon ,\eta})$, then\smallskip

(i) $E_{0,1}\! =\!\mathbb{C}\! =\! S_0$, $E_{1,0}\!\in\!\{ S_j\mid j\!\leq\! 3\}$, and 
$E_{\varepsilon ,\varepsilon}\!\in\!\{ S_j\mid 1\!\leq\! j\!\leq\! 3\}$,\smallskip

(ii) $\mathbb{R}_{\varepsilon ,\eta}\!\subseteq\!\Delta (E_{\varepsilon ,\eta})\cup\mathbb{P}_f^2$ (in particular, $\mathbb{R}_{0,1}$ is diagonally complex),\smallskip

\noindent and $(\mathbb{D},\mathbb{R})\sqsubseteq_c(X,R)$.\end{thm}

\noindent\bf Proof.\rm\ By Theorem \ref{contBr}, (a) and (b) cannot hold simutaneously. Assume that (a) does not hold. By Theorem \ref{contBr} again, one of the following holds.\medskip

\noindent (1) There is a relation $\mathbb{R}$ on $2^\omega$ such that 
$\mathbb{R}\cap\Delta (2^\omega )\! =\!\Delta (\mathbb{C})$ and $(2^\omega ,\mathbb{R})\sqsubseteq_c(X,R)$. By Corollary \ref{redtwo} we may assume that $\mathbb{R}\!\subseteq\!\Delta (2^\omega )\cup\mathbb{P}_f^2$, so that $\mathbb{R}$ is a diagonally complex relation.\medskip

\noindent (2) There is a relation $\mathbb{R}'$ on $\mathbb{D}$ such that 
$\mathbb{R}'_{0,1}\cap\Delta (2^\omega )\! =\!\Delta (\mathbb{C})$ and $(\mathbb{D},\mathbb{R}')\sqsubseteq_c(X,R)$. Note that 
$\mathbb{R}'_{\varepsilon ,\eta}$ is a locally countable Borel relation on $2^\omega$, for each 
$\varepsilon, \eta\!\in\! 2$. Corollary \ref{redtwo} provides $g'\! :\! 2^\omega\!\rightarrow\! 2^\omega$ injective continuous preserving 
$\mathbb{C}$ such that $(g'\!\times\! g')^{-1}(\bigcup_{\varepsilon ,\eta\in 2}~\mathbb{R}'_{\varepsilon ,\eta})
\!\subseteq\!\Delta (2^\omega )\cup\mathbb{P}_f^2$. We define $h\! :\!\mathbb{D}\!\rightarrow\!\mathbb{D}$ by 
$h(\varepsilon ,\alpha )\! :=\!\big(\varepsilon ,g'(\alpha )\big)$ and set $\mathbb{R}\! =\! (h\!\times\! h)^{-1}(\mathbb{R}')$, so that 
$\mathbb{R}_{\varepsilon ,\eta}\!\subseteq\!\Delta (2^\omega )\cup\mathbb{P}_f^2$. By Corollary \ref{Delta}, we may assume that 
$\mathbb{R}_{\varepsilon ,\eta}\cap\Delta (2^\omega )\!\in\!\{\Delta (S_j)\mid j\!\leq\! 3\}$. We are done since we are reduced to Case (1) if $\mathbb{R}_{\varepsilon, \varepsilon}\cap\Delta (2^\omega )\! =\!\Delta (\mathbb{C})$.\hfill{$\square$}

\subsection{$\!\!\!\!\!\!$ Proof of Theorem \ref{morethree}}\indent

 We now introduce a first antichain basis.\medskip

\noindent\bf Notation.\rm\ Let ${\bf\Gamma}$ be a non self-dual Borel class of rank at least two, and 
$(2^\omega ,\mathbb{C})$ be $\bf\Gamma$-good. We set
$$P\! :=\!\big\{ t\!\in\! 4^{(2^2)}\mid t(0,0),t(1,1)\!\not=\! 0\wedge t(0,1)\! =\! 0\wedge\big( t(1,0)\! =\! 0\Rightarrow 
t(0,0)\!\leq\! t(1,1)\big)\big\}$$
and, for $t\!\in\! P$, $\mathbb{R}^{\bf\Gamma}_t\! :=\!\big\{\big( (\varepsilon ,x),(\eta ,x)\big)\!\in\!\mathbb{D}^2\mid 
x\!\in\! S_{t(\varepsilon ,\eta )}\big\}$. We order $2^2$ lexicographically, so that, for example, 
$\mathbb{E}^{\bf\Gamma}_3\! =\!\mathbb{R}^{\bf\Gamma}_{3,0,0,3}$. Note that $\mathbb{G}^{\bf\Gamma}_m\! =\!\mathbb{R}^{\bf\Gamma}_{1,0,0,1}$. Finally, $\mathcal{A}^{\bf\Gamma}\! :=\!
\big\{\big( 2^\omega ,\Delta (\mathbb{C})\big)\big\}\cup\{ (\mathbb{D},\mathbb{R}^{\bf\Gamma}_t)\mid t\!\in\! P\}$. Note that the sections of the elements of 
$\mathcal{A}^{\bf\Gamma}$ have cardinality at most two, and are in particular closed.

\begin{lem} \label{moretwoant} Let ${\bf\Gamma}$ be a non self-dual Borel class of rank at least two. Then ${\mathcal A}^{\bf\Gamma}$ is a 34 elements $\leq_c$-antichain.\end{lem} 

\noindent\bf Proof.\rm\ Let $(\mathbb{X},\mathbb{R})\!\not=\! (\mathbb{X}',\mathbb{R}')$ in ${\mathcal A}^{\bf\Gamma}$. We argue by contradiction, which gives $f\! :\!\mathbb{X}\!\rightarrow\!\mathbb{X}'$ continuous. Assume first that 
$(\mathbb{X},\mathbb{R}),(\mathbb{X}',\mathbb{R}')$ are of the form 
$(\mathbb{D},\mathbb{R}^{\bf\Gamma}_t),(\mathbb{D},\mathbb{R}^{\bf\Gamma}_{t'})$ respectively, so that $f(\varepsilon ,\alpha )$ is of the form $\big( f_0(\varepsilon ,\alpha ),f_1(\varepsilon ,\alpha )\big)\!\in\! 2\!\times\! 2^\omega$.

\vfill\eject

 Let us prove that $f_0(0,\alpha )\!\not=\! f_0(1,\alpha )$ if $\alpha\!\in\!\mathbb{C}$. We argue by contradiction, which gives 
$l\!\in\!\omega$ such that $f_0(0,\beta )\! =\! f_0(1,\beta )\! =:\!\varepsilon$ if $\beta\!\in\! N_{\alpha\vert l}$, by continuity of $f_0$. This also gives continuous maps $g_\eta\! :\! N_{\alpha\vert l}\!\rightarrow\! 2^\omega$ such that $f_1(\eta ,\beta )\! =\! g_\eta (\beta )$ if 
$\beta\!\in\! N_{\alpha\vert l}$. If $\beta\!\in\!\mathbb{C}\cap N_{\alpha\vert l}$, then $\big( (0,\beta),(1,\beta )\big)\!\in\!\mathbb{R}$, so that $\big( f(0,\beta),f(1,\beta )\big)\! =\!\Big(\big(\varepsilon ,g_0(\beta)\big) ,(\varepsilon ,g_1(\beta )\big)\Big)\!\in\!\mathbb{R}'$, 
$g_0(\beta)\! =\! g_1(\beta )$ and $g_0\! =\! g_1\! =:\! g$ by continuity of $g_0,g_1$. Note that there is $j\!\in\!\{ 1,2,3\}$ such that 
$\mathbb{C}\cap N_{\alpha\vert l}\! =\! g^{-1}(S_j)\cap N_{\alpha\vert l}$, which contradicts the choice of $\mathbb{C}$.\medskip

 Fix $\alpha\!\in\!\mathbb{C}$. Note that there is $l\!\in\!\omega$ such that $\varepsilon_0\! :=\! f_0(0,\beta )\!\not=\! f_0(1,\beta )$ if 
$\beta\!\in\! N_{\alpha\vert l}$, by continuity of $f_0$. There are $g_\eta\! :\! N_{\alpha\vert l}\!\rightarrow\! 2^\omega$ continuous such that $f_1(\eta ,\beta )\! =\! g_\eta (\beta )$ if $\beta\!\in\! N_{\alpha\vert l}$. If ${\beta\!\in\!\mathbb{C}\cap N_{\alpha\vert l}}$, then 
$\big( (0,\beta),(1,\beta )\big)\!\in\!\mathbb{R}$, so that ${\big( f(0,\beta),f(1,\beta )\big)\! =\!
\Big(\big(\varepsilon_0,g_0(\beta)\big) ,(1\! -\!\varepsilon_0,g_1(\beta )\big)\Big)\!\in\!\mathbb{R}'}$, ${g_0(\beta)\! =\! g_1(\beta )}$ and $g_0\! =\! g_1\! =:\! g$ by continuity of $g_0,g_1$. Note that 
$\mathbb{C}\cap N_{\alpha\vert l}\! =\! g^{-1}(S_{t'(\varepsilon_0,1-\varepsilon_0)})\cap N_{\alpha\vert l}$, so that 
$t'(\varepsilon_0,1\! -\!\varepsilon_0)\! =\! 0$ by the choice of $\mathbb{C}$ and 
$\mathbb{C}\cap N_{\alpha\vert l}\! =\! g^{-1}(\mathbb{C})\cap N_{\alpha\vert l}$. If $\varepsilon_0\! =\! 0$, then 
$$S_{t(\varepsilon ,\eta )}\cap N_{\alpha\vert l}\! =\! g^{-1}(S_{t'(\varepsilon ,\eta )})\cap N_{\alpha\vert l}\! =\! 
S_{t'(\varepsilon ,\eta )}\cap N_{\alpha\vert l}$$ 
for $\varepsilon ,\eta\!\in\! 2$, so that $t\! =\! t'$ by the choice of $\mathbb{C}$. Thus $\varepsilon_0\! =\! 1$ and 
$$S_{t(\varepsilon ,\eta )}\cap N_{\alpha\vert l}\! =\! g^{-1}(S_{t'(1-\varepsilon ,1-\eta )})\cap 
N_{\alpha\vert l}\! =\! S_{t'(1-\varepsilon ,1-\eta )}\cap N_{\alpha\vert l}$$ 
for $\varepsilon ,\eta\!\in\! 2$, so that $t(\varepsilon ,\eta )\! =\! t'(1\! -\!\varepsilon ,1\! -\!\eta )$ if $\varepsilon ,\eta\!\in\! 2$, by the choice of $\mathbb{C}$. In particular, note that $t'(1,0)\! =\! t(0,1)\! =\! 0\! =\! t'(0,1)\! =\! t(1,0)$, $t'(0,0),t'(1,1)\!\not=\! 0$, 
$$t(1,1)\! =\! t'(0,0)\!\leq\! t'(1,1)\! =\! t(0,0)\mbox{,}$$ 
$t(0,0),t(1,1)\!\not=\! 0$ and $t'(1,1)\! =\! t(0,0)\!\leq\! t(1,1)\! =\! t'(0,0)$, so that $t\! =\! t'$ again.\medskip

 We now have to consider $\big( 2^\omega ,\Delta (\mathbb{C})\big)$. Assume first that 
$\big( 2^\omega ,\Delta (\mathbb{C})\big)\leq_c(\mathbb{D},\mathbb{R}^{\bf\Gamma}_t)$, with 
$(\mathbb{D},\mathbb{R}_t)$ in $\mathcal{A}^{\bf\Gamma}$ and witness $f$. Let $\alpha\!\in\!\mathbb{C}$. Then we can find  
$\varepsilon_0\!\in\! 2$ with $f_0(\alpha )\! =\!\varepsilon_0$, and $l\!\in\!\omega$ such that $f_0(\beta )\! =\!\varepsilon_0$ if 
$\beta\!\in\! N_{\alpha\vert l}$, by continuity of $f_0$. Note that 
$\mathbb{C}\cap N_{\alpha\vert l}\! =\! f_1^{-1}(S_{t(\varepsilon_0,\varepsilon_0)})\cap N_{\alpha\vert l}$, so that $t(\varepsilon_0,\varepsilon_0)\! =\! 0$ by the choice of $\mathbb{C}$, which contradicts the definition of $\mathcal{A}^{\bf\Gamma}$.\medskip

 Assume now that $(\mathbb{D},\mathbb{R}^{\bf\Gamma}_t)\leq_c\big( 2^\omega ,\Delta (\mathbb{C})\big)$, with 
$(\mathbb{D},\mathbb{R}_t)$ in $\mathcal{A}^{\bf\Gamma}$ and witness $f$. Let $\alpha\!\in\!\mathbb{C}$. Then 
$\big( (0,\alpha ),(1,\alpha )\big)\!\in\! \mathbb{R}^{\bf\Gamma}_t$, so that $f(0,\alpha )\! =\! f(1,\alpha )\!\in\!\mathbb{C}$. By the choice of $\mathbb{C}$, $g(\beta )\! :=\! f(0,\beta )\! =\! f(1,\beta )$ for each $\beta\!\in\! 2^\omega$. We set, for $\varepsilon ,\eta\!\in\! 2$, 
$\mathbb{R}_{\varepsilon ,\eta}\! :=\!\big\{ (\alpha ,\beta )\!\in\! 2^\omega\!\times\! 2^\omega\mid
\big( (\varepsilon ,\alpha ),(\eta ,\beta )\big)\!\in\!\mathbb{R}^{\bf\Gamma}_t\big\}$. Note that 
$$(\alpha ,\beta )\!\in\!\mathbb{R}_{\varepsilon ,\eta}\Leftrightarrow\big( g(\alpha ),g(\beta )\big)\!\in\!\Delta (\mathbb{C})\mbox{,}$$
so that $\mathbb{R}_{0,0}\! =\!\mathbb{R}_{0,1}$, which contradicts the definition of $\mathcal{A}^{\bf\Gamma}$.
\hfill{$\square$}\medskip

 The next result provides the basis part of Theorem \ref{morethree}.

\begin{lem} \label{moretwo} Let ${\bf\Gamma}$ be a non self-dual Borel class of rank at least two, $X$ be an analytic space, and $R$ be a locally countable Borel relation on $X$ whose sections are in $\Delta ({\bf\Gamma })$. Exactly one of the following holds:\smallskip  

(a) the relation $R$ is a $\bf\Gamma$ subset of $X^2$,\smallskip  

(b) there is $(\mathbb{X},\mathbb{R})\!\in\! {\mathcal A}^{\bf\Gamma}$ such that $(\mathbb{X},\mathbb{R})\sqsubseteq_c(X,R)$.\end{lem}

\noindent\bf Proof.\rm\ By Theorem \ref{contBr}, (a) and (b) cannot hold simutaneously. Assume that (a) does not hold. By Theorem \ref{contBr} again, one of the following holds.\medskip

\noindent (1) There is a relation $\mathbb{R}$ on $2^\omega$ such that $\mathbb{R}\cap\Delta (2^\omega )\! =\!\Delta (\mathbb{C})$ and $(2^\omega ,\mathbb{R})\sqsubseteq_c(X,R)$. As $R$ is locally countable Borel, so is $\mathbb{R}$ by Lemma \ref{basicr}, so that we can apply Corollary 3.3 in [L5] if the rank of $\bf\Gamma$ is at least three. This gives $g\! :\! 2^\omega\!\rightarrow\! 2^\omega$ injective continuous such that $g$ preserves $\mathbb{C}$, and 
$\big( g(\alpha ),g(\beta )\big)\!\notin\! \mathbb{R}$ if $\alpha\!\not=\!\beta$. Note that $g$ is a witness for the fact that 
$\big( 2^\omega ,\Delta(\mathbb{C})\big)\sqsubseteq_c(2^\omega ,\mathbb{R})$, so that $\big( 2^\omega ,\Delta(\mathbb{C})\big)\sqsubseteq_c(X,R)$. If the rank of $\bf\Gamma$ is two, then by Corollary \ref{redtwo} we may assume that 
$\mathbb{R}\!\subseteq\!\Delta (2^\omega )\cup\mathbb{P}_f^2$. As the sections of $R$ are countable and $\bormtwo$, so are those of $\mathbb{R}$. In particular, $\mathbb{R}$ has nowhere dense sections. By Corollary \ref{dc}, we may assume that $\mathbb{R}\! =\!\Delta (\mathbb{C})$, so that, here again, 
$\big( 2^\omega ,\Delta(\mathbb{C})\big)\sqsubseteq_c(X,R)$.\medskip

\noindent (2) There is a relation $\mathbb{R}$ on $\mathbb{D}$ such that 
$\mathbb{R}_{0,1}\cap\Delta (2^\omega )\! =\!\Delta (\mathbb{C})$ and $(\mathbb{D},\mathbb{R})\sqsubseteq_c(X,R)$. Note that $\mathbb{R}_{\varepsilon ,\eta}$ is a locally countable Borel relation on $2^\omega$. If the rank of $\bf\Gamma$ is at least three, then Corollary 3.3 in [L5] provides $g\! :\! 2^\omega\!\rightarrow\! 2^\omega$ injective continuous such that $g$ preserves $\mathbb{C}$, and 
$\big( g(\alpha ),g(\beta )\big)\!\notin\!\bigcup_{\varepsilon ,\eta\in 2}~\mathbb{R}_{\varepsilon ,\eta}$ if $\alpha\!\not=\!\beta$. If the rank of $\bf\Gamma$ is two, then Corollary \ref{redtwo} provides $g''\! :\! 2^\omega\!\rightarrow\! 2^\omega$ injective continuous preserving $\mathbb{C}$ such that $\mathbb{R}''\! :=\! 
(g''\!\times\! g'')^{-1}(\bigcup_{\varepsilon ,\eta\in 2}~\mathbb{R}_{\varepsilon ,\eta})\!\subseteq\!\Delta (2^\omega )\cup\mathbb{P}_f^2$. As the sections of $R$ are countable and $\bormtwo$, so are those of $\mathbb{R}''$. In particular, $\mathbb{R}''$ has nowhere dense sections. By Lemma \ref{nd}, we may assume that $\mathbb{R}''\!\subseteq\!\Delta (2^\omega )$. So we may assume that $g$ exists in both cases.\medskip

 We define $h\! :\!\mathbb{D}\!\rightarrow\!\mathbb{D}$ by $h(\varepsilon ,\alpha )\! :=\!\big(\varepsilon ,g(\alpha )\big)$, so that $h$ is injective and continuous. We then set 
 $\mathbb{R}'\! :=\! (h\!\times\! h)^{-1}(\mathbb{R})$. Repeating the notation above, $\mathbb{R}'_{\varepsilon ,\eta}\!\subseteq\!\Delta (2^\omega )$ by the property of $g$, and $h$ is a witness for the fact that $(\mathbb{D},\mathbb{R}')\sqsubseteq_c(\mathbb{D},\mathbb{R})$. This means that we may assume that 
${\mathbb{R}_{\varepsilon ,\eta}\!\subseteq\!\Delta (2^\omega )}$, for $\varepsilon ,\eta\!\in\! 2$, and that $\mathbb{R}_{0,1}\! =\!\Delta (\mathbb{C})\! =\!\Delta (S_0)$. By Corollary \ref{Delta}, we may assume that ${\mathbb{R}_{\varepsilon ,\eta}\!\in\!\{\Delta (S_j)\mid j\!\leq\! 3\}}$. Note that if 
$\mathbb{R}_{\varepsilon, \varepsilon}\! =\!\Delta (S_0)$ for some $\varepsilon\!\in\! 2$, then 
$\big( 2^\omega ,\Delta(\mathbb{C})\big)\sqsubseteq_c(X,R)$. So we may assume that 
$\mathbb{R}_{\varepsilon, \varepsilon}\! =\!\Delta (S_j)$ for some $j\!\in\!\{ 1,2,3\}$. Finally, 
$(\mathbb{D},\mathbb{R}^{\bf\Gamma}_{i,0,0,j})\sqsubseteq_c(\mathbb{D},\mathbb{R}^{\bf\Gamma}_{j,0,0,i})$ if 
$i\! >\! j\! >\! 0$ with witness $(\varepsilon ,\alpha )\!\mapsto\! (1\! -\! \varepsilon ,\alpha )$. So (b) holds.\hfill{$\square$}\medskip

\noindent\bf Proof of Theorem \ref{morethree}.\rm\ For (1) and (3), we apply Lemmas \ref{moretwoant} and \ref{moretwo}. For (2), we use the closure properties of $\sqsubseteq_c$ and the fact that $\mathbb{G}^{\bf\Gamma}_m$ is the only graph in $\mathcal{A}^{\bf\Gamma}$.\hfill{$\square$}

\section{$\!\!\!\!\!\!$ The rank two}\indent

 In this section, $\bf\Gamma$ is a non self-dual Borel class of rank two.
 
\subsection{$\!\!\!\!\!\!$ The diagonally complex relations: minimality and comparability}\indent
 
 Up to restrictions to Cantor sets, two comparable diagonally complex relations are bi-reducible.

\begin{lem} \label{perf} Let $\mathcal{R},\mathcal{R}_*$ be diagonally complex, and $K$ be a finitely dense Cantor set with the property that $(K,\mathcal{R}_*\cap K^2)\leq_c(2^\omega ,\mathcal{R})$. Then 
${(C,\mathcal{R}\cap C^2)\sqsubseteq_c(K,\mathcal{R}_*\cap K^2)}$ for some $C$ finitely dense.\end{lem}

\noindent\bf Proof.\rm\ Let $f$ be a witness for the fact that $(K,\mathcal{R}_*\cap K^2)\leq_c(2^\omega ,\mathcal{R})$. Note that $f$ preserves $\mathbb{C}$ since $\mathcal{R}_*,\mathcal{R}$ agree with $\Delta (\mathbb{C})$ on $\Delta (2^\omega )$. This implies that $f$ is nowhere dense-to-one. Let us prove that there is $g\! :\! 2^\omega\!\rightarrow\! K$ continuous such that $g$ preserves 
$\mathbb{C}$ and $f\big( g(\alpha )\big)\! <_{\text{lex}}\! f\big( g(\beta )\big)$ if $\alpha\!<_{\text{lex}}\!\beta$. We ensure (1)-(6) with 
$C\! =\! K$ and
$$(4)~f(\alpha )\! <_{\text{lex}}\! f(\beta )\mbox{ if }\alpha\!\in\! U_{t0}\mbox{ and }\beta\!\in\! U_{t1}$$
 
 Assume that this is done. If $\alpha\! <_{\text{lex}}\!\beta$, then Condition (4) ensures that $f\big( g(\alpha )\big)\! <_{\text{lex}}\! f\big( g(\beta )\big)$.\medskip
 
 So it is enough to prove that the construction is possible. We set $U_\emptyset\! :=\! K$. Assume that 
$(n_t)_{\vert t\vert\leq l}$ and $(U_t)_{\vert t\vert\leq l}$ satisfying (1)-(6) have been constructed, which is the case for 
$l\! =\! 0$. Fix $t\!\in\! 2^l$. Note that $\alpha_{n_t}$ is in $U_t\cap\mathbb{P}_f$, so that $f(\alpha_{n_t})\!\in\!\mathbb{P}_f$. This gives $s\!\in\! 2^{<\omega}$ such that $f(\alpha_{n_t})\! =\! s0^\infty$. In particular, there is a clopen neighbourhood 
$N\!\subseteq\! U_t$ of $\alpha_{n_t}$ such that $f(\beta )\!\in\! N_s$ if $\beta\!\in\! N$. We choose $n_{t1}\!\geq\! 1$ such that $\alpha_{n_{t1}}\!\in\! N\!\setminus\!\Big(\{\alpha_n\mid n\!\leq\! l\}\cup f^{-1}\big(\{ f(\alpha_{n_t})\}\big)\Big)$, which is possible since $f$ is nowhere dense-to-one. Note that $f(\alpha_{n_t})\! <_{\text{lex}}\! f(\alpha_{n_{t1}})$. It remains to choose a small enough clopen neighbourhood $U_{t\varepsilon}$ of $\alpha_{n_{t\varepsilon}}$ to finish the construction, using the continuity of $f$.\medskip

 Note then that $f_{\vert g[2^\omega ]}$ is a homeomorphism onto $\tilde C\! :=\! f\big[ g[2^\omega ]\big]$. Moreover, 
$$\big( g(\alpha ),g(\beta )\big)\!\in\!\mathcal{R}_*\Leftrightarrow\Big( f\big( g(\alpha )\big) ,f\big( g(\beta )\big)\Big)\!\in\!\mathcal{R}\mbox{,}$$ 
so that $f_{\vert g[2^\omega ]}^{-1}$ is a witness for the fact that 
$(\tilde C,\mathcal{R}\cap\tilde C^2)\sqsubseteq_c(K,\mathcal{R}_*\cap K^2)$, and $\mathcal{R}\cap\tilde C^2$ is not in 
$\bf\Gamma$ since the map $\alpha\!\mapsto\!\Big( f\big( g(\alpha )\big) ,f\big( g(\alpha )\big)\Big)$ reduces $\mathbb{C}$ to 
$\mathcal{R}\cap\tilde C^2$. This implies that $\tilde C\cap\mathbb{C}$ is not separable from 
$\tilde C\!\setminus\!\mathbb{C}$ by a $\bf\Gamma$ set. Using Theorem \ref{Lo-SR}, we get $m\! :\! 2^\omega\!\rightarrow\!\tilde C$ injective continuous such that $\mathbb{C}\! =\! m^{-1}(\mathbb{C})$. It remains to set $C\! :=\! m[2^\omega ]$.\hfill{$\square$}\medskip

 The minimality of diagonally complex relations can be seen on restrictions to Cantor sets.

\begin{lem} \label{G} Let $\mathcal{R}$ be a diagonally complex relation, $X$ be an analytic space, and $R$ be a non-$\bormtwo$ relation on $X$ such that 
${(X,R)\leq_c(2^\omega ,\mathcal{R})}$, with witness $f$. Then there is $g\! :\! 2^\omega\!\rightarrow\! X$ injective continuous such that 
$\mathbb{P}_f\! =\! g^{-1}\big( f^{-1}(\mathbb{P}_f)\big)$.\end{lem}

\noindent\bf Proof.\rm\ Note that $\mathcal{R}\! =\!\Delta (S)\cup (\mathcal{R}\cap\mathbb{P}_f^2)$, where $S\! =\!\mathbb{P}_\infty$ if 
${\bf\Gamma}\! =\!\boratwo$, and $S\! =\!\emptyset$ if ${\bf\Gamma}\! =\!\bormtwo$. As $\mathbb{P}_f^2\!\setminus\!\mathcal{R}$ is countable, it is a $\boratwo$ subset of 
$\mathbb{P}_f^2$, and $\mathcal{R}\cap\mathbb{P}_f^2$ is a $\bormtwo$ subset of $\mathbb{P}_f^2$, which provides $G\!\in\!\bormtwo (2^\omega\!\times\! 2^\omega )$ such that 
$\mathcal{R}\cap\mathbb{P}_f^2\! =\! G\cap\mathbb{P}_f^2$. As 
$R\! =\! (f\!\times\! f)^{-1}(\mathcal{R})\! =\! (f\!\times\! f)^{-1}\big(\Delta (S)\big)\cup\big( (f\!\times\! f)^{-1}(G)\cap f^{-1}(\mathbb{P}_f)^2\big)$, $f^{-1}(\mathbb{P}_f)$ is not 
$\bormtwo$. Theorem \ref{Lo-SR} provides $g$ as desired.\hfill{$\square$}

\begin{cor} \label{crimin} Let $\mathcal{R}$ be a diagonally complex relation, $X$ be an analytic space, and $R$ be a 
non-$\bf\Gamma$ relation on $X$ with $(X,R)\leq_c(2^\omega ,\mathcal{R})$. Then there is a finitely dense Cantor set $C$ such that 
${(C,\mathcal{R}\cap C^2)\sqsubseteq_c(X,R)}$.\end{cor}

\noindent\bf Proof.\rm\ Assume first that ${\bf\Gamma}\! =\!\bormtwo$, and let $f$ be a witness for the fact that 
$(X,R)\leq_c(2^\omega ,\mathcal{R})$. Lemma \ref{G} provides $g\! :\! 2^\omega\!\rightarrow\! X$ injective continuous such that 
$\mathbb{P}_f\! =\! g^{-1}\big( f^{-1}(\mathbb{P}_f)\big)$. We set $\mathcal{R}'\! :=\! (g\!\times\! g)^{-1}(R)$. Note that, for each 
$\alpha\!\in\! 2^\omega$, 
$$(\alpha ,\alpha )\!\in\!\mathcal{R}'\Leftrightarrow\big( g(\alpha ),g(\alpha )\big)\!\in\! R\Leftrightarrow
\Big( f\big( g(\alpha )\big) ,f\big( g(\alpha )\big)\Big)\!\in\!\mathcal{R}\Leftrightarrow f\big( g(\alpha )\big)\!\in\!\mathbb{P}_f\Leftrightarrow\alpha\!\in\!\mathbb{P}_f.$$ 
In other words, $\mathcal{R}'\cap\Delta (2^\omega )\! =\!\Delta (\mathbb{P}_f)$. This argument also shows that 
$\mathcal{R}'\!\subseteq\!\mathbb{P}_f^2$. In other words, $\mathcal{R}'$ is diagonally complex, 
$(2^\omega ,\mathcal{R}')\sqsubseteq_c(X,R)$ and $(2^\omega ,\mathcal{R}')\leq_c(2^\omega ,\mathcal{R})$. Lemma \ref{perf} provides a finitely dense Cantor set $C$ such that $(C,\mathcal{R}\cap C^2)\sqsubseteq_c(2^\omega ,\mathcal{R}')$.

\vfill\eject

 Assume now that ${\bf\Gamma}\! =\!\boratwo$. By Theorem \ref{contBr}, one of the following holds:\medskip

\noindent (1) there is a relation $\mathbb{R}_*$ on $2^\omega$ such that 
$\mathbb{R}_*\cap\Delta (2^\omega )\! =\!\Delta (\mathbb{P}_\infty )$ and $(2^\omega ,\mathbb{R}_*)\sqsubseteq_c(X,R)$.\smallskip

\noindent (2) there is a relation $\mathbb{R}_*$ on $\mathbb{D}$ such that 
$D_\infty\! :=\!\big\{\big( (0,\alpha ),(1,\alpha )\big)\mid\alpha\!\in\!\mathbb{P}_\infty\big\}\!\subseteq\!\mathbb{R}_*$, 
$$D_f\! :=\!\big\{\big( (0,\alpha ),(1,\alpha )\big)\mid\alpha\!\in\!\mathbb{P}_f\big\}\!\subseteq\!\neg\mathbb{R}_*$$ 
and $(\mathbb{D},\mathbb{R}_*)\sqsubseteq_c(X,R)$. We set $D\! :=\! D_\infty\cup D_f$.\medskip

 Assume that (2) holds, which gives $f\! :\!\mathbb{D}\!\rightarrow\! 2^\omega$ continuous such that $\mathbb{R}_*\! =\! (f\!\times\! f)^{-1}(\mathcal{R})$. Note that\medskip
 
\noindent - $(f\!\times\! f)[D_\infty ]\!\subseteq\!\mathcal{R}$ is not separable from $(f\!\times\! f)[D_f]$ by a $\boratwo$ set,\smallskip

\noindent - $(f\!\times\! f)[D_\infty ]\cap\Delta (2^\omega )$ is not separable from $(f\!\times\! f)[D_f]$ by a $\boratwo$ set since 
$\mathcal{R}\!\setminus\!\Delta (2^\omega)\!\subseteq\!\mathbb{P}_f^2$ is $\boratwo$,\smallskip

\noindent - $I_\infty\! :=\! (f\!\times\! f)[D_\infty ]\cap\Delta (2^\omega )$ is not separable from $I_f\! :=\! (f\!\times\! f)[D_f]\cap\Delta (2^\omega )$ by a $\boratwo$ set,\smallskip

\noindent - $R_\infty\! :=\! D\cap (f\!\times\! f)^{-1}( I_\infty )$ is not separable from $R_f\! :=\! D\cap (f\!\times\! f)^{-1}( I_f)$ by a $\boratwo$ set (otherwise 
$R_\infty\!\subseteq\! S\!\subseteq\!\neg R_f$ and $I_\infty$ is separable from $I_f$ by the $K_\sigma$ set $(f\!\times\! f)[S\cap D]$),\smallskip

\noindent - $C_\infty\! :=\!\big\{\alpha\!\in\! 2^\omega\vert\big( (0,\alpha ),(1,\alpha )\big)\!\in\! R_\infty\big\}$ is not separable from 
${C_f\! :=\!\big\{\alpha\!\in\! 2^\omega\vert\big( (0,\alpha ),(1,\alpha )\big)\!\in\! R_f\big\}}$ by a $\boratwo$ set.\medskip

 Theorem \ref{Lo-SR} provides $g\! :\! 2^\omega\!\rightarrow\! 2^\omega$ injective continuous with the properties that 
$\mathbb{P}_\infty\!\subseteq\! g^{-1}(C_\infty )$ and $\mathbb{P}_f\!\subseteq\! g^{-1}(C_f)$. We define a map 
$h\! :\! 2^\omega\!\rightarrow\!\{ 0\}\!\times\! 2^\omega$ by $h(\alpha )\! :=\! (0,\alpha )$, and set $c\! :=\! h\!\circ\! g$, $c'\! :=\! f\!\circ\! c$, $\mathbb{R}'\! :=\! (c'\!\times\! c')^{-1}(\mathcal{R})\! =\! (c\!\times\! c)^{-1}(\mathbb{R}_*)$. As $c$ is injective continuous, $\mathbb{R}'$ is a Borel relation on $2^\omega$ and $(2^\omega ,\mathbb{R}')\sqsubseteq_c(\mathbb{D},\mathbb{R}_*),(X,R)$. If 
$\alpha\!\in\!\mathbb{P}_\infty$, then $\Big(\big( 0,g(\alpha )\big) ,\big( 1,g(\alpha )\big)\Big)\!\in\! R_\infty$, and 
$\Big( f\big( 0,g(\alpha )\big) ,f\big( 1,g(\alpha )\big)\Big)\!\in\! I_\infty\!\subseteq\!\Delta (2^\omega )$. This implies that 
$\Big( f\big( 0,g(\alpha )\big) ,f\big( 0,g(\alpha )\big)\Big)\!\in\!\mathcal{R}$ and $(\alpha ,\alpha )\!\in\!\mathbb{R}'$. Similarly, if 
$\alpha\!\in\!\mathbb{P}_f$, then $(\alpha ,\alpha )\!\notin\!\mathbb{R}'$. Thus 
$\mathbb{R}'\cap\Delta (2^\omega )\! =\!\Delta (\mathbb{P}_\infty )$, and $\mathbb{R}'$ is a witness for the fact that (1) also holds. So (1) holds in any case.\medskip

 This implies that $(2^\omega ,\mathbb{R}_*)\leq_c(2^\omega ,\mathcal{R})$, with witness $f'$. By Lemma \ref{pre}, $f'$ preserves 
$\mathbb{P}_f$. In particular, 
$\mathbb{R}_*\cap\big( (\mathbb{P}_f\!\times\!\mathbb{P}_\infty )\cup (\mathbb{P}_\infty\!\times\!\mathbb{P}_f)\big)\! =\!\emptyset$ since $\mathcal{R}$ is diagonally complex. Note that $\mathbb{R}_*$ has $\boratwo$ sections since $\mathcal{R}$ has. This implies that 
$\mathbb{R}_*\cap\mathbb{P}_\infty^2$ is a Borel relation on $\mathbb{P}_\infty$ whose sections are separable from $\mathbb{P}_f$ by a $\boratwo$ set. Lemmas \ref{bot} and \ref{inf} provide $g'\! :\! 2^\omega\!\rightarrow\! 2^\omega$ injective continuous preserving 
$\mathbb{P}_f$ such that $(g'\!\times\! g')^{-1}(\mathbb{R}_*\cap\mathbb{P}_\infty^2)\!\subseteq\!\Delta (\mathbb{P}_\infty )$. Thus we may assume that $\mathbb{R}_*$ is diagonally complex. We now apply Lemma \ref{perf}.\hfill{$\square$}\medskip

 A consequence of Lemma \ref{perf} is a characterisation of the minimality of diagonally complex relations.

\begin{cor} \label{critmin} Let $\mathcal{R}$ be a diagonally complex relation. The following are equivalent:\smallskip

(a) $\mathcal{R}$ is $\leq_c$ and $\sqsubseteq_c$-minimal among non-$\bf\Gamma$ relations on an analytic space,\smallskip

(b) $(2^\omega ,\mathcal{R})\sqsubseteq_c(C,\mathcal{R}\cap C^2)$ if $C$ is a finitely dense Cantor set.\end{cor}

\noindent\bf Proof.\rm\ (a) $\Rightarrow$ (b) As $\mathbb{P}_f\cap C$ is dense in $C$, $\mathcal{R}\cap\Delta (C)$ and 
$\mathcal{R}\cap C^2$ are not in $\bf\Gamma$. As 
$$(C,\mathcal{R}\cap C^2)\sqsubseteq_c(2^\omega ,\mathcal{R})$$ 
and (a) holds, (b) holds.\medskip

\noindent (b) $\Rightarrow$ (a) Let $X$ be an analytic space and $R$ be a non-$\bf\Gamma$ relation on $X$ with 
$(X,R)\leq_c(2^\omega ,\mathcal{R})$. Corollary \ref{crimin} provides a finitely dense Cantor set $C$ with the property that 
${(C,\mathcal{R}\cap C^2)\sqsubseteq_c(X,R)}$. By (b), $(2^\omega ,\mathcal{R})\sqsubseteq_c(X,R)$.\hfill{$\square$}\medskip

 Another consequence of Lemma \ref{perf} is about the comparison of minimal diagonally complex relations.

\begin{cor} \label{compdiag} Let $\mathcal{R},\mathcal{R}'$ be diagonally complex relations, $\leq_c$ and $\sqsubseteq_c$-minimal among among non-$\bf\Gamma$ relations on an analytic space, with 
$(2^\omega ,\mathcal{R})\leq_c(2^\omega ,\mathcal{R}')$. Then 
$$(2^\omega ,\mathcal{R}')\sqsubseteq_c(2^\omega ,\mathcal{R})$$ 
and $(2^\omega ,\mathcal{R})\sqsubseteq_c(2^\omega ,\mathcal{R}')$.\end{cor}

\noindent\bf Proof.\rm\ Lemma \ref{perf} provides a finitely dense Cantor set $C$ such that 
$(C,\mathcal{R}'\cap C^2)\sqsubseteq_c(2^\omega ,\mathcal{R})$. By Corollary \ref{critmin}, 
$(2^\omega ,\mathcal{R}')\sqsubseteq_c(2^\omega ,\mathcal{R})$. An application of this fact implies that 
$(2^\omega ,\mathcal{R})\sqsubseteq_c(2^\omega ,\mathcal{R}')$.\hfill{$\square$}

\subsection{$\!\!\!\!\!\!$ Proof of Theorems \ref{two} and \ref{baum}}\indent

 We now introduce our antichain of size continuum.\medskip
 
\noindent\bf Notation.\rm\ We define $i\! :\! Q^2\!\rightarrow\!\omega$ as follows. We want to ensure that $i(z,t)\! =\! i(s_z,s_t)$, where $s_t$ is defined before Lemma \ref{inf}. The definition of $i$ is partly inspired by the oscillation map osc defined in [T] after Theorem 6.33 as follows. The elements of $2^{<\omega}$ are identified with finite subsets of $\omega$, through the characteristic function. The oscillation between $z$ and $t$ describes the behavious of the symmetric difference $z\Delta t$. The equivalence relation $\sim_{zt}$ is defined on $z\Delta t$ by 
$$j\sim_{zt}k\Leftrightarrow [\mbox{min}(j,k),\mbox{max}(j,k)]\cap (z\!\setminus\! t)\! =\!\emptyset\vee [\mbox{min}(j,k),\mbox{max}(j,k)]\cap (t\!\setminus\! z)\! =\!\emptyset .$$
Then $\mbox{osc}(z,t)\! :=\!\vert z\Delta t/\!\sim_{zt}\!\!\vert$. For $i$, we work on $z\cup t$ instead of $z\Delta t$, and the definition depends more heavily on the initial segments of $z$ and $t$, in particular on their lexicographic ordering. If 
$z\!\in\! 2^{<\omega}\!\setminus\!\{\emptyset\}$, then we set 
${z^-\! :=\! z\vert\mbox{max}\{ l\! <\!\vert z\vert\mid z\vert l\!\in\! Q\}}$. We also set 
$$\perp_<:=\!\{ (z,t)\!\in\! Q^2\mid\exists i\! <\!\mbox{min}(\vert z\vert ,\vert t\vert )~~z\vert i\! =\! t\vert i\wedge z(i)\! <\! t(i)\}$$ 
and, similarly, $\perp_>:=\!\{ (z,t)\!\in\! Q^2\mid\exists i\! <\!\mbox{min}(\vert z\vert ,\vert t\vert )~~z\vert i\! =\! t\vert i\wedge z(i)\! >\! t(i)\}$. The definition of $i(z,t)$ is by induction on 
$\mbox{max}(\vert z\vert ,\vert t\vert )$. We set 
$$i(z,t)\! :=\!\left\{\!\!\!\!\!\!\!
\begin{array}{ll}
& 0\mbox{ if }z\! =\! t\mbox{,}\cr
& i(z,t^-)\mbox{ if }\vert z\vert\! <\!\vert t^-\vert\vee\big(\vert z\vert\! =\!\vert t^-\vert\wedge (z,t^-)\!\in\perp_<\!\big)\mbox{,}\cr
& i(z,t^-)\! +\! 1\mbox{ if }(\vert z\vert\! <\!\vert t\vert\wedge\vert z\vert\! >\!\vert t^-\vert )\vee
\big(\vert z\vert\! =\!\vert t^-\vert\wedge (z,t^-)\!\notin\perp_<\!\big)\mbox{,}\cr
& i(z^-,t)\mbox{ if }\vert t\vert\! <\!\vert z^-\vert\vee\big(\vert t\vert\! =\!\vert z^-\vert\wedge (t,z^-)\!\in\perp_<\!\big)\mbox{,}\cr
& i(z^-,t)\! +\! 1\mbox{ if }(\vert t\vert\! <\!\vert z\vert\wedge\vert t\vert\! >\!\vert z^-\vert )\vee
\big(\vert t\vert\! =\!\vert z^-\vert\wedge (t,z^-)\!\notin\perp_<\!\big)\mbox{,}\cr
& i(z^-,t^-)\! +\! 1\mbox{ if }\vert z\vert\! =\!\vert t\vert\!\wedge\!\big( (\vert z^-\vert\! <\!\vert t^-\vert\!\wedge\! t^-\! <_{\text{lex}}\! z^-)\!\vee\!
(\vert t^-\vert\! <\!\vert z^-\vert\!\wedge\! z^-\! <_{\text{lex}}\! t^-)\big)\mbox{,}\cr
& i(z^-,t^-)\! +\! 2\mbox{ if }\vert z\vert\! =\!\vert t\vert\!\wedge\!\big( (\vert z^-\vert\! <\!\vert t^-\vert\!\wedge\! z^-\! <_{\text{lex}}\! t^-)\!\vee\!
(\vert t^-\vert\! <\!\vert z^-\vert\!\wedge\! t^-\! <_{\text{lex}}\! z^-)\vee\cr
& \hfill{(\vert z^-\vert\! =\!\vert t^-\vert\!\wedge\! z^-\!\not=\! t^-)\big) .}
\end{array}
\right.$$
Note that $i(z,t)\! =\! i(t,z)$ if $z,t\!\in\! Q$.

\begin{lem} \label{suff} Let $(s_t)_{t\in 2^{<\omega}}$ be a sequence of elements of $Q$ with 
$$\begin{array}{ll}
& (a)~\vert s_z\vert\! <\!\vert s_t\vert\mbox{ if }z\! <_l\! t\mbox{ are in }Q\cr
& (b)~s_t\!\subsetneqq\! s_{t1}\mbox{ if }t\mbox{ is in }Q
\end{array}$$
Then $i(z,t)\! =\! i(s_z,s_t)\mbox{ if }z,t\!\in\! Q$.\end{lem}

\noindent\bf Proof.\rm\ We argue by induction on 
$\mbox{max}(\vert z\vert ,\vert t\vert )$. As $i(z,t)\! =\! i(t,z)$ if $z,t\!\in\! Q$, we may assume that $z\! <_lt$. We go through the cases of the definition of $i$.\medskip

 If $\vert z\vert\! <\!\vert t^-\vert$, then $\vert s_z\vert\! <\!\vert s_{t^-}\vert\!\leq\!\big\vert (s_t)^-\big\vert$ since $s_{t^-}\!\subsetneqq\! s_t$, and
$$i(z,t)\! =\! i(z,t^-)\! =\! i(s_z,s_{t^-})\! =\!\cdots\! =\! i\big( s_z,(s_t)^-\big)\! =\! i(s_z,s_t).$$
If $\vert z\vert\! =\!\vert t^-\vert$ and $(z,t^-)\!\in\perp_<$, then 
$\vert s_z\vert\! <\!\vert s_{t^-}\vert\!\leq\!\vert (s_t)^-\vert$ again and we conclude as above.\medskip

 If $\vert z\vert\! <\!\vert t\vert$ and $\vert z\vert\! >\!\vert t^-\vert$, then $\vert s_z\vert\! <\!\vert s_t\vert$ and $\vert s_z\vert\! >\!\vert s_{t^-}\vert$. Let $s\!\in\! Q$ with 
$s\!\subseteq\!  s_t$ and $\vert s_z\vert\! <\!\vert s\vert$ be of minimal length. Note that $i(z,t)\! =\! i(z,t^-)\! +\! 1\! =\! i(s_z,s_{t^-})\! +\! 1\! =\! i(s_z,s)\! =\! i(s_z,s_t)$.\medskip

 If $\vert z\vert\! =\!\vert t^-\vert$ and $(z,t^-)\!\notin\perp_<$, then either $z\! =\! t^-$, or $(z,t^-)\!\in\perp_>$, so that
 $$i(z,t)\! =\! i(z,t^-)\! +\! 1\! =\! i(s_z,s_{t^-})\! +\! 1\! =\!\cdots\! =\! i\big( s_z,(s_t)^-\big)\! =\! i(s_z,s_t).$$
If $\vert z\vert\! =\!\vert t\vert$, $\vert t^-\vert\! <\!\vert z^-\vert$ and $z^-\! <_{\text{lex}}\! t^-$, then 
$\vert s_{t^-}\vert\! <\!\vert s_{z^-}\vert\! <\!\vert s_z\vert\! <\!\vert s_t\vert$,
$$i(z,t)\! =\! i(z^-,t^-)\! +\! 1\! =\! i(s_{z^-},s_{t^-})\! +\! 1\! =\! i(s_z,s_{t^-})\! +\! 1\! =\! i(s_z,s_t).$$
If $\vert z\vert\! =\!\vert t\vert$, $\vert z^-\vert\! <\!\vert t^-\vert$ and $z^-\! <_{\text{lex}}\! t^-$, then 
$\vert s_{z^-}\vert\! <\!\vert s_{t^-}\vert\! <\!\vert s_z\vert\! <\!\vert s_t\vert$,
$$i(z,t)\! =\! i(z^-,t^-)\! +\! 2\! =\! i(s_{z^-},s_{t^-})\! +\! 2\! =\! i(s_z,s_{t^-})\! +\! 1\! =\! i(s_z,s_t).$$
If $\vert z\vert\! =\!\vert t\vert$, $\vert z^-\vert\! =\!\vert t^-\vert$ and $z^-\!\not=\! t^-$, then $(z^-,t^-)\!\in\perp_<$, 
$\vert s_{z^-}\vert\! <\!\vert s_{t^-}\vert\! <\!\vert s_z\vert\! <\!\vert s_t\vert$, and we conclude as above. Thus $i(z,t)\! =\! i(s_z,s_t)$ in any case, as desired.\hfill{$\square$}\medskip

The next result is the key lemma to prove Theorems \ref{two} and \ref{baum}.

\begin{lem} \label{extract} Let $H\!\subseteq\!\mathbb{P}_f$ be homeomorphic to $\mathbb{P}_f$. Then there is 
$f\! :\! 2^\omega\!\rightarrow\! 2^\omega$ injective continuous satisfying the following properties:\smallskip

(i) $f[\mathbb{P}_f]\!\subseteq\! H$; in particular, for each $t\!\in\! Q$ there is $n_t\!\geq\! 1$ with $f(t0^\infty )\! =\! s_t0^\infty$,\smallskip

(ii) $f[\mathbb{P}_\infty ]\!\subseteq\!\mathbb{P}_\infty$,\smallskip

(iii) $i(z,t)\! =\! i(s_z,s_t)\mbox{ if }z,t\!\in\! Q$.\end{lem}

\noindent\bf Proof.\rm\ We ensure (1)-(6) with $C\! :=\!\overline{H}$, that $U_t$ is of the form $N_{z_t}\cap C$, and 
$$\begin{array}{ll}
& (2)~\alpha_{n_t}\!\in\! N_{z_t}\cap H\mbox{ if }t\mbox{ is in }2^{<\omega}\cr
& (4)~z_{t0}\! <_{\text{lex}}\! z_{t1}\mbox{ and }\vert z_{t0}\vert\! =\!\vert z_{t1}\vert\mbox{ if }t\mbox{ is in }2^{<\omega}\cr
& (7)~\vert s_z\vert\! <\!\vert s_t\vert\mbox{ if }z\! <_l\! t\mbox{ are in }Q\cr
& (8)~s_t\!\subsetneqq\! s_{t1}\mbox{ if }t\mbox{ is in }Q\cr
& (9)~i(z,t)\! =\! i(s_z,s_t)\mbox{ if }z,t\!\in\! Q
\end{array}$$

 It is enough to prove that the construction is possible. We choose $n_\emptyset\!\geq\! 1$ with $s_\emptyset 0^\infty\!\in\! H$, and set $z_\emptyset\! :=\! s_\emptyset$. Assume that $(n_t)_{\vert t\vert\leq l}$ and $(z_t)_{\vert t\vert\leq l}$ satisfying (1)-(9) have been constructed, which is the case for $l\! =\! 0$.\medskip

 Fix $t\!\in\! 2^l$. As $\alpha_{n_t}\!\in\! U_t\! =\! N_{z_t}\cap C$, $z_t\!\subseteq\! s_t0^\infty$. We choose $s\!\subseteq\! s_t0^\infty$ extending $z_t$ and $s_t$, in such a way that 
$N_s\cap (\{\alpha_n\mid n\!\leq\! l\}\!\setminus\!\{\alpha_{n_t}\} )\! =\!\emptyset$ and $\vert s_z\vert\! <\!\vert s\vert$ if $z\! <_l\! t1$. Let $n_{t1}\!\geq\! 1$ such that 
$\alpha_{n_{t1}}\!\in\! N_s\cap H\!\setminus\!\{\alpha_{n_t}\}$, which is possible by density of $H$ in the perfect space $C$. We set $z_{t1}\! =\! s_{t1}$ and 
$z_{t0}\! =\!\alpha_{n_t}\vert\vert z_{t1}\vert$. Note that (1) holds. (2) holds by definition. By our extensions, (3) holds. (4) holds because of the choice of $s$ and $n_{t1}$. (5)-(8) hold by construction. By Lemma \ref{suff}, (9) holds.\hfill{$\square$}\medskip

\noindent\bf Proof of Theorem \ref{baum}.\rm\ We may replace $\mathbb{Q}$ with its topological copy $\mathbb{P}_f$. We set, for 
$z,t\!\in\! Q$, $c(\{ z0^\infty ,t0^\infty\} )\! :=\! i(z,t)$, which is well-defined by symmetry of $i$. As 
$i\big( (10)^k0^\infty ,(01)^k0^\infty\big)\! =\! 2k$ and $i\big( (01)^k0^\infty,(10)^{k+1}0^\infty\big)\! =\! 2k\! +\! 1$, $c$ is onto. If 
$H\!\subseteq\!\mathbb{P}_f$ is homeomorphic to $\mathbb{P}_f$, then Lemma \ref{extract} provides 
$f\! :\! 2^\omega\!\rightarrow\! 2^\omega$. We set $h\! :=\! f_{\vert\mathbb{P}_f}$. If $z,t\!\in\! Q$, then 
$$c\big(\{ h(z0^\infty ),h(t0^\infty )\}\big)\! =\! c(\{ s_z0^\infty ,s_t0^\infty\} )\! =\! i(s_z,s_t)\! =\! i(z,t)\! =\! c(\{ z0^\infty ,t0^\infty\} ) .$$ 
In particular, $c$ takes all the values from $\omega$ on $H^{[2]}$.\hfill{$\square$}\medskip

\noindent\bf Notation.\rm\ We now set, when $(2^\omega ,\mathbb{C})$ is $\bf\Gamma$-good and $\beta\!\in\! 2^\omega$,
$$\mathbb{R}_\beta\! :=\!\Delta (\mathbb{C})\cup\bigcup_{\beta (p)=1}~
\{ (z0^\infty ,t0^\infty )\mid z\not=\! t\!\in\! Q\wedge ~i(z,t)\! =\! p\}\mbox{,}$$
so that $\mathbb{R}_\beta$ is symmetric.

\begin{cor} \label{contmin} Let $\beta\!\in\! 2^\omega$. Then $\mathbb{R}_\beta$ is $\leq_c$ and $\sqsubseteq_c$-minimal among 
non-$\bf\Gamma$ relations on an analytic space.\end{cor}

\noindent\bf Proof.\rm\ Note that $\mathbb{R}_\beta$ is a diagonally complex relation, and therefore not in $\bf\Gamma$. Let $C$ be a finitely dense Cantor set. By Corollary \ref{critmin}, it is enough to see that the inequality $(2^\omega ,\mathbb{R}_\beta)\sqsubseteq_c(C,\mathbb{R}_\beta\cap C^2)$ holds. We apply Lemma \ref{extract} to $H\! :=\! C\cap\mathbb{P}_f$. As $C$ is finitely dense, $H$ is nonempty, countable, metrizable, perfect, and therefore homeomorphic to $\mathbb{P}_f$ (see 7.12 in [K]). Lemma \ref{extract} provides a witness $f\! :\! 2^\omega\!\rightarrow\! 2^\omega$ for the fact that $(2^\omega ,\mathbb{R}_\beta)\sqsubseteq_c(C,\mathbb{R}_\beta\cap C^2)$. 
\hfill{$\square$}

\begin{lem} \label{antitwo} The family $(\mathbb{R}_\beta )_{\beta\in N_0}$ is a $\leq_c$-antichain.\end{lem}

\noindent\bf Proof.\rm\ Assume that $\beta ,\beta'\!\in\! N_0$ and 
$(2^\omega ,\mathbb{R}_\beta )\leq_c(2^\omega ,\mathbb{R}_{\beta'})$, with witness $f$. By Lemma \ref{pre}, $f$ preserves $\mathbb{C}$.\medskip

\noindent\bf Claim.\it\ Let $u_0\!\in\! 2^{<\omega}$. Then we can find $u'_0\!\in\! Q$ and $u_1,u'_1\!\in\! Q\!\setminus\!\{\emptyset\}$ such that $f(u_00^\infty )\! =\! u'_00^\infty$ and $f(u_0u_10^\infty )\! =\! u'_0u'_10^\infty$.\rm\medskip

 Indeed, let $u'_0\!\in\! Q$ with $f(u_00^\infty )\! =\! u'_00^\infty$. We can find a sequence $(\beta_k)_{k\in\omega}$ of points of 
$\mathbb{P}_\infty$ converging to $u_00^\infty$. As $f$ preserves $\mathbb{C}$, $f(u_00^\infty )\!\not=\! f(\beta_k)$, and we can find 
$(n_k)_{k\in\omega}$ stricly increasing such that ${\alpha_{n_k}\vert k\! =\!\beta_k\vert k}$ and 
$f(u_00^\infty )\!\not=\! f(\alpha_{n_k})$, so that $(\alpha_{n_k})_{k\in\omega}$ is a sequence of points different from 
$u_00^\infty$ converging to $u_00^\infty$. Let $v_k,v'_k\!\in\! Q$ with $\alpha_{n_k}\! =\! v_k0^\infty$ and 
$f(\alpha_{n_k})\! =\! v'_k0^\infty$. We may assume that $u_0\!\subsetneqq\! v_k$. As $f$ is continuous, we may assume that $f(\alpha_{n_k})\!\in\! N_{u'_0}$ for each $k$, so that $u'_0\!\subsetneqq\! v'_k$. It remains to choose 
$u_1,u'_1\!\in\! Q\!\setminus\!\{\emptyset\}$ such that $v_1\! =\! u_0u_1$ and $v'_1\! =\! u'_0u'_1$.\hfill{$\diamond$}\medskip

 Assume that $t_0,\cdots ,t_{2k+2}\!\in\! 2^{<\omega}$ are not initial segments of $0^\infty$. Then 
$$(t_0t_10^{\vert t_2\vert}\cdots t_{2k-1}0^{\vert t_{2k}\vert}t_{2k+1}0^\infty ,
t_00^{\vert t_1\vert}t_2\cdots 0^{\vert t_{2k+1}\vert}t_{2k+2}0^\infty )$$ 
is in $\{ (z0^\infty ,t0^\infty )\mid z\not=\! t\!\in\! Q\wedge ~i(z,t)\! =\! 2k\! +\! 2\}$. Similarly, 
$$(t_00^{\vert t_1\vert}t_2\cdots 0^{\vert t_{2k-1}\vert}t_{2k}0^\infty,
t_0t_10^{\vert t_2\vert}\cdots t_{2k-1}0^{\vert t_{2k}\vert}t_{2k+1}0^\infty )$$ 
is in $\{ (z0^\infty ,t0^\infty )\mid z\not=\! t\!\in\! Q\wedge ~i(z,t)\! =\! 2k\! +\! 1\}$. We first apply the claim to $u_0\! :=\! t_0\!\in\! Q$, which gives $t'_0\!\in\! Q$ and $t_1,t'_1\!\in\! Q\!\setminus\!\{\emptyset\}$ with the properties that 
$f(t_00^\infty )\! =\! t'_00^\infty\mbox{ and }f(t_0t_10^\infty )\! =\! t'_0t'_10^\infty$. The continuity of $f$ provides 
$k_1\!\geq\!\vert t_1\vert$ such that $f[N_{t_00^{k_1}}]\!\subseteq\! N_{t'_00^{\vert t'_1\vert}}$. We next apply the claim to 
$u_0\! :=\! t_00^{k_1}$, which gives $\tilde t_2,\tilde t'_2$ in $Q\!\setminus\!\{\emptyset\}$ such that 
$f(t_00^{k_1}\tilde t_20^\infty )\! =\! t'_0\tilde t'_20^\infty$. Note that $0^{\vert t'_1\vert}\!\subseteq\!\tilde t'_2$. We set 
$t_2\! :=\! 0^{k_1-\vert t_1\vert}\tilde t_2$ and $t'_2\! :=\!\tilde t'_2\! -\! 0^{\vert t'_1\vert}$, so that $t_2,t'_2$ are not initial segments of 
$0^\infty$ and 
$$f(t_00^{\vert t_1\vert}t_20^\infty )\! =\! t'_00^{\vert t'_1\vert}t'_20^\infty .$$ 
This argument shows that we can find sequences of finite binary sequences $(t_j)_{j\in\omega}$ and $(t'_j)_{j\in\omega}$ which are not initial segments of $0^\infty$ and satisfy $f(t_00^{\vert t_1\vert}t_2\cdots 0^{\vert t_{2k-1}\vert}t_{2k}0^\infty )\! =\! 
t'_00^{\vert t'_1\vert}t'_2\cdots 0^{\vert t'_{2k-1}\vert}t'_{2k}0^\infty$ and 
$f(t_0t_10^{\vert t_2\vert}\cdots t_{2k-1}0^{\vert t_{2k}\vert}t_{2k+1}0^\infty )\! =\! 
t'_0t'_10^{\vert t'_2\vert}\cdots t'_{2k-1}0^{\vert t'_{2k}\vert}t'_{2k+1}0^\infty$ for each $k\!\in\!\omega$. By the remark after the claim, 
$\beta\! =\!\beta'$.\hfill{$\square$}\medskip

\noindent\bf Proof of Theorem \ref{two}.\rm\ We apply Corollary \ref{contmin} and Lemma \ref{antitwo}.\hfill{$\square$}\medskip

 Theorem \ref{two} shows that if $\leq$ is in $\{\leq_c,\sqsubseteq_c\}$, then the class of non-$\bormtwo$ countable relations on analytic spaces, equipped with $\leq$, contains antichains of size continuum made of minimal relations.

\subsection{$\!\!\!\!\!\!$ Graphs}\indent

 Theorem \ref{two} shows that, for the classes of rank two, any basis must have size continuum. It is natural to ask whether it is also the case for graphs. We will see that it is indeed the case.\medskip

\noindent\bf Notation.\rm\ We define, for $\beta\!\in\! 2^\omega$, 
$\mathbb{G}_\beta\! :=\! s\Big(\big\{\big( (0,\alpha ),(1,\gamma )\big)\!\in\!\mathbb{D}^2\mid (\alpha ,\gamma )\!\in\!\mathbb{R}_\beta\big\}\Big)$. Note that $\mathbb{G}_\beta$ is a locally countable graph.

\begin{lem} \label{both} Let ${\bf\Gamma}\! :=\!\boratwo$, $(2^\omega ,\mathbb{C})$ be $\bf\Gamma$-good, 
$\beta\!\in\! 2^\omega$, $X$ be an analytic space, $R$ be a relation on $X$, and $\mathbb{R}$ be a relation on $\mathbb{D}$ such that $\mathbb{R}_{0,1}\cap\Delta (2^\omega )\! =\!\Delta (\mathbb{C})$ and 
$(\mathbb{D},\mathbb{R})\sqsubseteq_c(X,R)\sqsubseteq_c(\mathbb{D},\mathbb{G}_\beta )$. Then 
$(2^\omega ,\mathbb{R}_\beta )\sqsubseteq_c(2^\omega ,\mathbb{R}_{0,1})$ with witness $h$ having the property that  
$(\mathbb{R}_{0,0}\cup\mathbb{R}_{1,1})\cap h[2^\omega ]^2\! =\!\emptyset$.\end{lem}

\noindent\bf Proof.\rm\ Let $g,f$ be witnesses for the fact that $(\mathbb{D},\mathbb{R})\sqsubseteq_c(X,R)$, $(X,R)\sqsubseteq_c(\mathbb{D},\mathbb{G}_\beta )$ respectively.\medskip

 We set, for $\varepsilon\!\in\! 2$, $X_\varepsilon\! :=\!\{ x\!\in\! X\mid f_0(x)\! =\!\varepsilon\}$, which defines a partition of $X$ into clopen sets. The definition of $\mathbb{G}_\beta$ shows that $R\!\subseteq\! (X_0\!\times\! X_1)\cup (X_1\!\times\! X_0)$. If 
$\alpha_0\!\in\!\mathbb{C}$, then $\big( g(0,\alpha_0),g(1,\alpha_0)\big)\!\in\! R$, which gives $\varepsilon_0\!\in\! 2$ and $l\!\in\!\omega$ such that $g(0,\alpha )\!\in\! X_{\varepsilon_0}$ and $g(1,\alpha )\!\in\! X_{1-\varepsilon_0}$ if $\alpha\!\in\! N_{\alpha_0\vert l}$. This shows that $\mathbb{R}\cap (2\!\times\! N_{\alpha_0\vert l})^2\!\subseteq\! 
s\Big(\big\{\big( (\varepsilon_0,\alpha ),(1\! -\!\varepsilon_0,\gamma )\big)\mid\alpha ,\gamma\!\in\! 2^\omega\big\}\Big)$. Note that 
$\mathbb{R}_{\varepsilon ,\eta}$ is a relation on $2^\omega$. As $\mathbb{R}_{0,1}\cap\Delta (2^\omega )\! =\!\Delta (\mathbb{C})$, 
$\mathbb{R}_{0,1}$ is not in $\bf\Gamma$. Moreover, by symmetry of $\mathbb{R}_\beta$,
$$\begin{array}{ll}
(\alpha ,\gamma )\!\in\!\mathbb{R}_{0,1}\!\!\!\!
& \Leftrightarrow\big( (0,\alpha ),(1,\gamma)\big)\!\in\!\mathbb{R}\Leftrightarrow\big( g(0,\alpha ),g(1,\gamma )\big)\!\in\! R\cr
& \Leftrightarrow\Big( f\big( g(0,\alpha )\big) ,f\big( g(1,\gamma )\big)\Big)\!\in\!\mathbb{G}_\beta\Leftrightarrow
\Big( f_1\big( g(0,\alpha )\big) ,f_1\big( g(1,\gamma )\big)\Big)\!\in\!\mathbb{R}_\beta
\end{array}$$
if $\alpha\!\in\! N_{\alpha_0\vert l}$.

\vfill\eject

\noindent\bf Claim.\it\ $f_1\big( g(0,\alpha )\big)\! =\! f_1\big( g(1,\alpha )\big)$ if $\alpha\!\in\!\mathbb{C}\cap N_{\alpha_0\vert l}$.\rm\medskip

 Indeed, note that $\alpha\!\in\!\mathbb{C}\Leftrightarrow\big( (0,\alpha ),(1,\alpha )\big)\!\in\!\mathbb{R}\Leftrightarrow
\Big( f_1\big( g(0,\alpha )\big) ,f_1\big( g(1,\alpha )\big)\Big)\!\in\!\mathbb{R}_\beta$ if ${\alpha\!\in\! N_{\alpha_0\vert l}}$. We argue by contradiction, which gives 
$\alpha_0\!\in\!\mathbb{C}\cap N_{\alpha_0\vert l}$, $t\!\in\! 2^{<\omega}$, $\eta\!\in\! 2$ and $l'\!\geq\! l$ with the properties that $f_1\big( g(0,\alpha )\big)\!\in\! N_{t\eta}$ and $f_1\big( g(1,\alpha )\big)\!\in\! N_{t(1-\eta )}$ if $\alpha\!\in\! N_{\alpha_0\vert l'}$. In particular, 
$$\mathbb{C}\cap N_{\alpha_0\vert l'}\! =\!\Big\{\alpha\!\in\! N_{\alpha_0\vert l'}\mid
\Big( f_1\big( g(0,\alpha )\big) ,f_1\big( g(1,\alpha )\big)\Big)\!\in\!\mathbb{R}_\beta\cap (N_{t\eta}\!\times\! N_{t(1-\eta )})\Big\}\mbox{,}$$ 
which shows that $\mathbb{C}\cap N_{\alpha_0\vert l'}\!\in\! {\bf\Gamma}\! =\!\boratwo$ since 
$\mathbb{R}_\beta\!\setminus\!\Delta (2^\omega )$ is countable, contradicting the choice of $\mathbb{C}$.\hfill{$\diamond$}\medskip

 By the claim, continuity of $f,g$ and density of $\mathbb{C}$, $f_1\big( g(0,\alpha )\big)\! =\! f_1\big( g(1,\alpha )\big)$ if $\alpha\!\in\! N_{\alpha_0\vert l}$. This shows that 
$(N_{\alpha_0\vert l},\mathbb{R}_{0,1}\cap N_{\alpha_0\vert l}^2)\leq_c(2^\omega ,\mathbb{R}_\beta )$. By the proof of Corollary \ref{critmin} and Corollary \ref{contmin}, 
$$(2^\omega ,\mathbb{R}_\beta )\sqsubseteq_c(N_{\alpha_0\vert l},\mathbb{R}_{0,1}\cap N_{\alpha_0\vert l}^2)\mbox{,}$$
and we are done.\hfill{$\square$}

\begin{lem} \label{contminstruct} Let ${\bf\Gamma}\! :=\!\boratwo$, $(2^\omega ,\mathbb{C})$ be $\bf\Gamma$-good, and 
$\beta\!\in\! 2^\omega$. Then $\mathbb{G}_\beta$ is $\sqsubseteq_c$-minimal among non-$\bf\Gamma$ relations on an analytic space.\end{lem}

\noindent\bf Proof.\rm\ Let $X$ be an analytic space, and $R$ be a non-$\bf\Gamma$ relation on $X$ with 
$(X,R)\sqsubseteq_c(\mathbb{D},\mathbb{G}_\beta )$. As $\mathbb{G}_\beta$ has sections in $\bf\Gamma\! =\!\boratwo$, $R$ too. By Theorem \ref{contBr}, one of the following holds:\medskip  

\noindent (1) there is a relation $\mathbb{R}$ on $2^\omega$ such that 
$\mathbb{R}\cap\Delta (2^\omega )\! =\!\Delta (\mathbb{C})$ and $(2^\omega ,\mathbb{R})\sqsubseteq_c(X,R)$,\smallskip

\noindent (2) there is a relation $\mathbb{R}$ on $\mathbb{D}$ such that 
$\mathbb{R}_{0,1}\cap\Delta (2^\omega )\! =\!\Delta (\mathbb{C})$ and $(\mathbb{D},\mathbb{R})\sqsubseteq_c(X,R)$.\medskip

 As $\mathbb{G}_\beta$ is a digraph, $R$ and $\mathbb{R}$ are digraphs too. So (1) cannot hold. By Lemma \ref{both}, 
$$(2^\omega ,\mathbb{R}_\beta )\sqsubseteq_c(2^\omega ,\mathbb{R}_{0,1})\mbox{,}$$ 
with witness $h$ having the property that $(\mathbb{R}_{0,0}\cup\mathbb{R}_{1,1})\cap h[2^\omega ]^2\! =\!\emptyset$. We define $k\! :\!\mathbb{D}\!\rightarrow\!\mathbb{D}$ by $k(\varepsilon ,\alpha )\! :=\!\big(\varepsilon ,h(\alpha )\big)$. Note that $k$ is injective, continuous, and 
$$\big( (0,\alpha ),(1,\gamma )\big)\!\in\!\mathbb{G}_\beta\Leftrightarrow
\big( h(\alpha ),h(\gamma )\big)\!\in\!\mathbb{R}_{0,1}\Leftrightarrow\big( k(0,\alpha ),k(1,\gamma )\big)\!\in\!\mathbb{R}
\mbox{,}$$
showing that $(\mathbb{D},\mathbb{G}_\beta )\sqsubseteq_c(X,R)$ as desired.\hfill{$\square$}

\begin{lem} \label{perpstruct} Let ${\bf\Gamma}\! :=\!\boratwo$, and $(2^\omega ,\mathbb{C})$ be $\bf\Gamma$-good. Then $(\mathbb{G}_\beta )_{\beta\in N_0}$ is a 
$\sqsubseteq_c$-antichain.\end{lem}

\noindent\bf Proof.\rm\ Let $\beta ,\beta'\!\in\! N_0$. Assume that $(\mathbb{D},\mathbb{G}_\beta )\sqsubseteq_c(\mathbb{D},\mathbb{G}_{\beta'})$. Lemma \ref{both} shows that 
$$(2^\omega ,\mathbb{R}_{\beta'})\sqsubseteq_c\big( 2^\omega ,(\mathbb{G}_\beta )_{0,1}\big)\! =\! (2^\omega ,\mathbb{R}_\beta )\mbox{,}$$ 
so that $\beta\! =\!\beta'$ by Lemma \ref{antitwo}.\hfill{$\square$}

\begin{cor} \label{Fsigmastruct} There is a concrete $\sqsubseteq_c$-antichain of size continuum made of locally countable Borel graphs on $\mathbb{D}$ which are $\sqsubseteq_c$-minimal among non-$\boratwo$ graphs on an analytic space.\end{cor}

\noindent\bf Proof.\rm\ We apply Lemmas \ref{contminstruct} and \ref{perpstruct}.\hfill{$\square$}\medskip

 We now turn to the study of the class $\bormtwo$.

\vfill\eject

\begin{lem} \label{sect} Assume that ${\bf\Gamma}\! =\!\bormtwo$ and $\beta\!\in\! N_0\!\setminus\!\{ 0^\infty\}$. Then 
$\mathbb{R}_\beta$ has a non-$\bormtwo$ section.\end{lem}

\noindent\bf Proof.\rm\ We argue by contradiction, so that 
$\mathbb{R}_\beta$ has nowhere dense sections since it is locally countable. By Corollary \ref{dc}, 
$\big( 2^\omega ,\Delta (\mathbb{C})\big)\sqsubseteq_c(2^\omega ,\mathbb{R}_\beta )$. Note that 
$(2^\omega ,\mathbb{R}_\beta )\sqsubseteq_c\big( 2^\omega ,\Delta (\mathbb{C})\big)$ by Theorem \ref{morethree} and Corollary \ref{compdiag}, which contradicts the injectivity since $\beta\!\in\! N_0\!\setminus\!\{ 0^\infty\}$.\hfill{$\square$}\medskip

\noindent\bf Remark.\rm\ Lemma \ref{sect} shows that a locally countable relation can have a non-$\bormtwo$ section, and therefore be non-$\bormtwo$ because of that. This is not the case for the class $\boratwo$. This is the reason why we cannot argue for the class 
$\bormtwo$ as we did in Corollary \ref{Fsigmastruct} for the class $\boratwo$. More precisely, note that $\mathbb{G}^{{\bormtwo},a}_m$ is a (locally) countable graph with a 
non-$\bormtwo$ section. Moreover, $(\mathbb{S},\mathbb{G}^{{\bormtwo},a}_m)\sqsubseteq_c(\mathbb{D},\mathbb{G}_\beta )$ for any $\beta\!\in\! N_0\!\setminus\!\{ 0^\infty\}$. Indeed, Lemma \ref{sect} gives $\alpha_0\!\in\! 2^\omega$ such that 
$(\mathbb{R}_\beta )_{\alpha_0}$ is not $\bormtwo$ since $\mathbb{R}_\beta$ is symmetric. By Theorem \ref{Lo-SR}, there is 
$g\! :\! 2^\omega\!\rightarrow\! 2^\omega$ injective continuous such that 
$\mathbb{P}_f\! =\! g^{-1}\big( (\mathbb{R}_\beta )_{\alpha_0}\big)$. We define $f\! :\!\mathbb{S}\!\rightarrow\!\mathbb{D}$ by $f(0^\infty )\! :=\! (0,\alpha_0)$ and 
$f(1\alpha )\! :=\!\big( 1,g(\alpha )\big)$, so that $f$ is injective continuous and 
$$\begin{array}{ll}
(\varepsilon\alpha ,\eta\gamma )\!\in\!\mathbb{G}^{{\bormtwo},a}_m\!\!\!\!
& \Leftrightarrow (\varepsilon\! =\! 0\wedge\alpha\! =\! 0^\infty\wedge\eta\! =\! 1\wedge\gamma\!\in\!\mathbb{P}_f)\vee 
(\varepsilon\! =\! 1\wedge\alpha\!\in\!\mathbb{P}_f\wedge\eta\! =\! 0\wedge\gamma\! =\! 0^\infty )\cr
& \Leftrightarrow\bigg(\big( f(\varepsilon\alpha ),f(\eta\gamma )\big)\! =\!\Big( (0,\alpha_0),\big( 1,g(\gamma )\big)\Big)\wedge
\big(\alpha_0,g(\gamma )\big)\!\in\!\mathbb{R}_\beta\bigg)\vee\cr
& \hfill{\bigg(\big( f(\varepsilon\alpha ),f(\eta\gamma )\big)\! =\!\Big( \big( 1,g(\alpha )\big) ,(0,\alpha_0)\Big)\wedge
\big(\alpha_0,g(\alpha )\big)\!\in\!\mathbb{R}_\beta\bigg)}\cr
& \Leftrightarrow\big( f(\varepsilon\alpha ),f(\eta\gamma )\big)\!\in\!\mathbb{G}_\beta .
\end{array}$$
So we need to find other examples for the class $\bormtwo$.

\begin{thm} \label{Gdeltastruct} There is a concrete $\sqsubseteq_c$-antichain of size continuum made of locally countable Borel graphs on $2^\omega$ which are $\sqsubseteq_c$-minimal among non-$\bormtwo$ graphs on an analytic space.\end{thm}

\noindent\bf Proof.\rm\ We set, for $\beta\!\in\! N_0\!\setminus\!\{ 0^\infty\}$, 
$\mathbb{G}_\beta\! :=\!\mathbb{R}_\beta\!\setminus\!\Delta (2^\omega)$. As $\mathbb{R}_\beta$ is symmetric, 
$\mathbb{G}_\beta$ is a graph. In particular, we can speak of $\mbox{proj}[\mathbb{G}_\beta ]$. By Lemma \ref{sect}, 
$\mathbb{R}_\beta$ has a non-$\bormtwo$ section. A consequence of this is that $\mathbb{G}_\beta$ is not $\bormtwo$. Let $X$ be an analytic space and $R$ be a non-$\bormtwo$ relation on $X$ with 
${(X,R)\!\sqsubseteq_c(2^\omega ,\mathbb{G}_\beta )}$, with witness $f$. By injectivity of $f$, we get 
${\Big( X,R\cup (f\!\times\! f)^{-1}\big(\Delta (\mathbb{C})\big)\Big)\sqsubseteq_c(2^\omega ,\mathbb{R}_\beta )}$. As 
${\Delta (\mathbb{C})\!\in\!\boratwo}$ is disjoint from $\mathbb{G}_\beta$, 
$R\cup (f\!\times\! f)^{-1}\big(\Delta (\mathbb{C})\big)$ is not $\bormtwo$. By Corollary \ref{contmin} and injectivity of $f$, 
${(2^\omega ,\mathbb{R}_\beta )\sqsubseteq_c\Big( X,R\cup (f\!\times\! f)^{-1}\big(\Delta (\mathbb{C})\big)\Big)}$. By injectivity of $f$ again, $(2^\omega ,\mathbb{G}_\beta )\sqsubseteq_c(X,R)$, showing the minimality of $\mathbb{G}_\beta$.\medskip

 Assume now that $\beta ,\beta'\!\in\! N_0\!\setminus\!\{ 0^\infty\}$ and 
$(2^\omega ,\mathbb{G}_\beta )\sqsubseteq_c(2^\omega ,\mathbb{G}_{\beta'})$, with witness $f$.\medskip

\noindent\bf Claim.\it\ Let $u_0\!\in\! 2^{<\omega}$ with $u_00^\infty\!\in\!\mbox{proj}[\mathbb{G}_\beta ]$. Then we can find $u'_0\!\in\! Q$ and 
$u_1,u'_1\!\in\! Q\!\setminus\!\{\emptyset\}$ such that $f(u_00^\infty )\! =\! u'_00^\infty$, $u_0u_10^\infty\!\in\!\mbox{proj}[\mathbb{G}_\beta ]$ and 
$f(u_0u_10^\infty )\! =\! u'_0u'_10^\infty$.\rm\medskip

 Indeed, as $u_00^\infty\!\in\!\mbox{proj}[\mathbb{G}_\beta ]$, there is $\gamma$ such that $(u_00^\infty ,\gamma )$ is in 
$\mathbb{G}_\beta$, and $\big( f(u_00^\infty ),f(\gamma )\big)$ is therefore in $\mathbb{G}_{\beta'}$. In particular, 
$f(u_00^\infty )\!\in\!\mathbb{P}_f$ and there is $u'_0\!\in\! Q$ with $f(u_00^\infty )\! =\! u'_00^\infty$. As $f$ is continuous, there is $k_0$ such that $f[N_{u_00^{k_0}}]\!\subseteq\! N_{u'_0}$. \medskip

 Let us prove that $\mbox{proj}[\mathbb{G}_\beta ]$ is dense in $2^\omega$. Let $s\!\in\! 2^{<\omega}$, and $p\! :=\! 2k\! +\!\varepsilon\!\geq\! 1$ with $\beta (p)\! =\! 1$. We set 
$(z,t)\! :=\!\big( s1(10)^{k+\varepsilon}1,s1(01)^k\big)$, so that $z\!\not=\! t\!\in\! Q$, $i(z,t)\! =\! p$,  
$(z0^\infty ,t0^\infty )\!\in\!\mathbb{G}_\beta\cap N_s^2$ and $\mbox{proj}[\mathbb{G}_\beta ]$ meets $N_s$ as desired.\medskip

 Let $\alpha\!\in\!\mbox{proj}[\mathbb{G}_\beta ]\cap N_{u_00^{k_0}1}$, and $\delta$ such that $(\alpha ,\delta )$ is in $\mathbb{G}_\beta$. Then 
$\big( f(\alpha ),f(\delta )\big)$ is in $\mathbb{G}_{\beta'}$, and $f(\alpha )\!\in\!\mathbb{P}_f\cap N_{u'_0}\!\setminus\!\{ u'_00^\infty\}$ by injectivity of $f$. This gives 
$u_1,u'_1\!\in\! Q\!\setminus\!\{\emptyset\}$ with $\alpha\! =\! u_0u_10^\infty$ and $f(\alpha )\! =\! u'_0u'_10^\infty$.\hfill{$\diamond$}\medskip

 The end of the proof is now similar to that of Lemma \ref{antitwo}. Let us indicate the differences. We first apply the claim to 
$u_0\! :=\! t_0\!\in\! Q$ such that $t_00^\infty\!\in\!\mbox{proj}[\mathbb{G}_\beta ]$, which gives $t'_0\!\in\! Q$ and 
$t_1,t'_1\!\in\! Q\!\setminus\!\{\emptyset\}$ such that 
$f(t_00^\infty )\! =\! t'_00^\infty\mbox{, }t_0t_10^\infty\!\in\!\mbox{proj}[\mathbb{G}_\beta ]\mbox{ and }
f(t_0t_10^\infty )\! =\! t'_0t'_10^\infty$. The continuity of $f$ provides $k_1\!\geq\!\vert t_1\vert$ such that 
$f[N_{t_00^{k_1}}]\!\subseteq\! N_{t'_00^{\vert t'_1\vert}}$. We next apply the claim to $u_0\! :=\! t_00^{k_1}$, which gives 
$\tilde t_2,\tilde t'_2$ in $Q\!\setminus\!\{\emptyset\}$ such that $t_00^{k_1}\tilde t_20^\infty\!\in\!\mbox{proj}[\mathbb{G}_\beta ]$ and $f(t_00^{k_1}\tilde t_20^\infty )\! =\! t'_0\tilde t'_20^\infty$. This provides sequences $(t_j)_{j\in\omega}$ and $(t'_j)_{j\in\omega}$ as in the proof of Lemma \ref{antitwo}. By the remark after the claim in the proof of Lemma \ref{antitwo}, $\beta\! =\!\beta'$.\hfill{$\square$}
  
\section{$\!\!\!\!\!\!$ Acyclicity}

\bf Remark.\rm\ Assume that $\beta\!\in\! N_0\!\setminus\!\{ 0^\infty\}$. Then $\mathbb{R}_\beta$ is not s-acyclic. Indeed,
$(0^\infty ,10^\infty,1^20^\infty )$ is a $s(\mathbb{R}_\beta )$-cycle if $\beta (1)\! =\! 1$, $(10^\infty ,010^\infty,0^210^\infty )$ is a $s(\mathbb{R}_\beta )$-cycle if $\beta (2)\! =\! 1$, 
$$(101^{k+1}0^\infty ,1^{k+3}0^\infty,0101^k0^\infty )$$ 
is a $s(\mathbb{R}_\beta )$-cycle if $\beta (2k\! +\! 3)\! =\! 1$, and 
$(0101^20(110)^k0^\infty ,1010^3(10^2)^k0^\infty,010^21^2(101)^k0^\infty )$ is a $s(\mathbb{R}_\beta )$-cycle if 
$\beta (2k\! +\! 4)\! =\! 1$. We will see that Theorem \ref{morethree} can be extended, under a suitable acyclicity assumption. We need to introduce new examples.\medskip

\noindent\bf Notation.\rm\ Let ${\bf\Gamma}$ be a non self-dual Borel class of rank at least two, and $(2^\omega ,\mathbb{C})$ be $\bf\Gamma$-good. We set 
$A\! :=\!\big\{ t\!\in\! 4^{(2^2)}\mid t(0,0)\!\in\! 2\wedge\big( t(0,1)\! =\! 0\vee t(1,0)\! =\! 0\big)\wedge t(1,1)\!\not=\! 0\big\}$. We then set, for $t\!\in\! A$,\medskip

\leftline{$\mathbb{R}^{{\bf\Gamma},a}_t\! :=\!\{ (0^\infty ,0^\infty )\mid t(0,0)\! =\! 1\}\cup
\{ (0^\infty ,1\alpha )\mid\alpha\!\in\! S_{t(0,1)}\}\cup\{ (1\alpha ,0^\infty )\mid\alpha\!\in\! S_{t(1,0)}\} ~\cup$}\smallskip

\rightline{$\{ (1\alpha ,1\alpha )\mid\alpha\!\in\! S_{t(1,1)}\} .$}\medskip
 
\noindent Note that $\mathbb{G}^{{\bf\Gamma},a}_m\! =\!\mathbb{R}^{{\bf\Gamma},a}_{0,0,0,1}$. Finally 
$\mathcal{B}^{\bf\Gamma}\! :=\!\big\{ (\mathbb{S},\mathbb{R}^{{\bf\Gamma},a}_t)\mid t\!\in\! A\big\}$.

\begin{lem} \label{antichac} Let ${\bf\Gamma}$ be a non self-dual Borel class of rank at least two. Then 
${\mathcal A}^{\bf\Gamma}\cup {\mathcal B}^{\bf\Gamma}$ is a 76 elements $\leq_c$-antichain.\end{lem}

\noindent\bf Proof.\rm\ Assume that $t,t'\!\in\! A$ and $\mathbb{R}^{{\bf\Gamma},a}_t$ is $\leq_c$-below 
$\mathbb{R}^{{\bf\Gamma},a}_{t'}$ with witness $f\! :\!\mathbb{S}\!\rightarrow\!\mathbb{S}$. If $\mathbb{R}^{{\bf\Gamma},a}_{t'}$ is symmetric, then $\mathbb{R}^{{\bf\Gamma},a}_t$ is too, and $t'$ is of the form $(\varepsilon ,0,0,\varepsilon')$. As 
$\mathbb{R}^{{\bf\Gamma},a}_t,\mathbb{R}^{{\bf\Gamma},a}_{t'}$ have only one vertical section not in $\bf\Gamma$, 
$f(0^\infty )\! =\! 0^\infty$, $f(1,\alpha )(0)\! =\! 1$ by the choice of $\mathbb{C}$, and the function $f_1(1,.)$ defined by 
$1f_1(1,\alpha )\! =\! f(1\alpha )$ preserves $\mathbb{C}$. Thus $t\! =\! t'$. So we may assume that $\mathbb{R}^{{\bf\Gamma},a}_{t'}$ is not symmetric, i.e., $t'$ is of the form $(\varepsilon ,\varepsilon',\varepsilon'',\varepsilon''')$ with $\varepsilon'\!\not=\!\varepsilon''$, and 
$\varepsilon'\! =\! 0$ or $\varepsilon''\! =\! 0$. If, for example, $\varepsilon'\!\not=\! 0$, then $\mathbb{R}^{{\bf\Gamma},a}_{t'}$ has vertical sections in $\bf\Gamma$, as well as $\mathbb{R}^{{\bf\Gamma},a}_t$, so that $t$ is of the form $(\eta ,\eta',0,\eta''')$ with 
$\eta'\!\not=\! 0$. Note that $\mathbb{R}^{{\bf\Gamma},a}_t,\mathbb{R}^{{\bf\Gamma},a}_{t'}$ have only one horizontal section not in 
$\bf\Gamma$, $f(0^\infty )\! =\! 0^\infty$, $f(1,\alpha )(0)\! =\! 1$ by the choice of $\mathbb{C}$, and $f_1(1,.)$ preserves $\mathbb{C}$. Thus $t\! =\! t'$, even if $\varepsilon''\!\not=\! 0$.\medskip

As the elements of ${\mathcal B}^{\bf\Gamma}$ have a section not in $\bf\Gamma$ and the elements of 
${\mathcal A}^{\bf\Gamma}$ have closed sections, an element of ${\mathcal B}^{\bf\Gamma}$ is not $\leq_c$-below an element of ${\mathcal A}^{\bf\Gamma}$. By Theorem \ref{morethree}, it remains to see that an element $\mathbb{R}$ of 
${\mathcal A}^{\bf\Gamma}$ is not $\leq_c$-below an element $\mathbb{R}^{{\bf\Gamma},a}_t$ of 
${\mathcal B}^{\bf\Gamma}$. We argue by contradiction, which provides $f\! :\! 2^\omega\!\rightarrow\!\mathbb{S}$ or 
$f\! :\!\mathbb{D}\!\rightarrow\!\mathbb{S}$.

\vfill\eject

 In the first case, note first that $f(\alpha )(0)\! =\! 1$, since otherwise $f(\beta )(0)\! =\! 0$ if $\beta$ is in a clopen neighbourhood $C$ of $\alpha$. If $\beta\!\not=\!\gamma\!\in\! C\cap\mathbb{C}$, then $f(\beta )\! =\! f(\gamma )\! =\! 0^\infty$, so that 
$(0^\infty ,0^\infty)\!\in\!\mathbb{R}^{{\bf\Gamma},a}_t$ and $(\beta ,\gamma )\!\in\!\Delta (\mathbb{C})$, which is absurd. This shows that $\Delta (\mathbb{C})\!\in\!\bf\Gamma$ since $(\mathbb{R}^{{\bf\Gamma},a}_t)_{1,1}\!\in\! {\bf\Gamma}$, which is absurd. In the second case, note first that $f(0,\alpha )(0)\! =\! 1$, since otherwise $f(0,\beta )(0)\! =\! 0$ if $\beta$ is in a clopen neighbourhood $C$ of 
$\alpha$. We may assume that there is $\varepsilon_0\!\in\! 2$ with the property that $f(1,\beta )(0)\! =\!\varepsilon_0$ if $\beta\!\in\! C$. If $\beta\!\not=\!\gamma\!\in\! C\cap\mathbb{C}$, then either $\varepsilon_0\! =\! 0$, 
$\big( f(0,\beta ),f(1,\gamma )\big)\! =\! (0^\infty ,0^\infty )\!\in\!\mathbb{R}^{{\bf\Gamma},a}_t$ and 
${\big( (0,\beta ),(1,\gamma )\big)\!\in\!\mathbb{R}}$, or $\varepsilon_0\! =\! 1$, 
$\big( f(0,\beta ),f(1,\beta )\big) ,\big( f(0,\gamma ),f(1,\gamma )\big)\!\in\!\mathbb{R}^{{\bf\Gamma},a}_t$, $f(1,\beta ),f(1,\gamma )$ are in $\{ 1\alpha\mid\alpha\!\in\!\mathbb{C}\}$, $\big( f(0,\beta ),f(1,\gamma )\big)\!\in\!\mathbb{R}^{{\bf\Gamma},a}_t$ and 
$\big( (0,\beta ),(1,\gamma )\big)\!\in\!\mathbb{R}$, which is absurd. Similarly, $f(1,\alpha )(0)\! =\! 1$. This shows that 
$\mathbb{R}_{0,1}\!\in\!\bf\Gamma$ since $(\mathbb{R}_t^{\bf\Gamma})_{1,1}\!\in\!\bf\Gamma$, which is absurd.\hfill{$\square$}\medskip

 We now study the rank two. 
 
\begin{lem} \label{locacy} Let $R$ be a s-acyclic Borel relation on $\mathbb{P}_\infty$. Then we can find a sequence 
$(R_n)_{n\in\omega}$ of relations closed in $\mathbb{P}_\infty\!\times\! 2^\omega$ and $2^\omega\!\times\!\mathbb{P}_\infty$, as well as $f\! :\! 2^\omega\!\rightarrow\! 2^\omega$ injective continuous preserving $\mathbb{P}_f$ such that 
$(f\!\times\! f)^{-1}(R)\!\subseteq\!\bigcup_{n\in\omega}~R_n$.\end{lem}

\noindent\bf Proof.\rm\ Assume first that $R_{\alpha_0}$ is not separable from $\mathbb{P}_f$ by a set in $\boratwo$, for some 
$\alpha_0\!\in\!\mathbb{P}_\infty$. Theorem \ref{Lo-SR} provides $h\! :\! 2^\omega\!\rightarrow\! 2^\omega\!\setminus\!\{\alpha_0\}$ injective continuous such that $\mathbb{P}_\infty\!\subseteq\! h^{-1}(R_{\alpha_0})$ and $\mathbb{P}_f\!\subseteq\! h^{-1} (\mathbb{P}_f)$. We set $R'\! :=\! (h\!\times\! h)^{-1}(R)$. If $\alpha ,\beta ,\gamma\!\in\!\mathbb{P}_\infty$ are pairwise distinct and 
$\beta ,\gamma\!\in\! R'_\alpha$, then $\big( h(\beta ),h(\alpha ),h(\gamma ),\alpha_0\big)$ is a $s(R)$-cycle, which is absurd. Thus $R'$ has vertical sections of cardinality at most two. So, replacing $R$ with $R'$ if necessary, we may assume that $R_\alpha$ is separable from $\mathbb{P}_f$ by a set in $\boratwo$ for each $\alpha\!\in\! 2^\omega$. Similarly, we may assume that $R^\alpha$ is separable from $\mathbb{P}_f$ by a set in $\boratwo$ for each $\alpha\!\in\! 2^\omega$. It remains to apply Lemma \ref{bot}.\hfill{$\square$}
  
\begin{lem} \label{finsq} (a) Let $R$ be a s-acyclic subrelation of $\mathbb{P}_f^2$. Then there is $f\! :\! 2^\omega\!\rightarrow\! 2^\omega$ injective continuous preserving 
$\mathbb{P}_f$ such that $(f\!\times\! f)^{-1}(R)\!\subseteq\!\Delta (\mathbb{P}_f)$.\smallskip

(b) Let $R$ be a s-acyclic subrelation of $(2\!\times\!\mathbb{P}_f)^2$. Then there is $f\! :\! 2^\omega\!\rightarrow\! 2^\omega$ injective continuous preserving 
$\mathbb{P}_f$ such that $(f\!\times\! f)^{-1}(\bigcup_{\varepsilon ,\eta\in 2}~R_{\varepsilon ,\eta})\!\subseteq\!\Delta (\mathbb{P}_f)$.\end{lem}

\noindent\bf Proof.\rm\ (a) Assume that $R_{\alpha_0}$ is not nowhere dense for some $\alpha_0\!\in\!\mathbb{P}_f$. This gives $s\!\in\! 2^{<\omega}$ with the property that  
$N_s\!\subseteq\!\overline{R_{\alpha_0}}$. Note that $N_s\cap R_{\alpha_0}\!\setminus\!\{\alpha_0\}$ is not separable from $N_s\cap\mathbb{P}_\infty$ by a $\bormtwo$ set, by Baire's theorem. Theorem \ref{Lo-SR} provides ${h\! :\! 2^\omega\!\rightarrow\! N_s\!\setminus\!\{\alpha_0\}}$ injective continuous such that 
${\mathbb{P}_f\!\subseteq\! h^{-1}(R_{\alpha_0})}$ and 
$\mathbb{P}_\infty\!\subseteq\! h^{-1}(\mathbb{P}_\infty )$. If $R'\! :=\! (h\!\times\! h)^{-1}(R)$, 
$\alpha ,\beta ,\gamma\!\in\!\mathbb{P}_f$ are pairwise distinct and $\beta ,\gamma\!\in\! R'_\alpha$, then $\big( h(\beta ),h(\alpha ),h(\gamma ),\alpha_0\big)$ is a $s(R)$-cycle, which is absurd. This shows that, replacing $R$ with $R'$ if necessary, we may assume that $R$ has vertical sections of cardinality at most two. In any case, we may assume that $R$ has nowhere dense vertical sections. Similarly, we may assume that $R$ has nowhere dense horizontal sections. It remains to apply Lemma \ref{nd}.\medskip

\noindent (b) We argue similarly. Fix $\varepsilon ,\eta\!\in\! 2$. We replace $R$ with $R_{\varepsilon ,\eta}$. We just have to note that the cycle 
$\big( h(\beta ),h(\alpha ),h(\gamma ),\alpha_0\big)$ becomes 
$\Big(\big(\eta ,h(\beta )\big) ,\big(\varepsilon ,h(\alpha )\big) ,\big(\eta ,h(\gamma )\big) ,\big(\varepsilon ,\alpha_0\big)\Big)$.\hfill{$\square$}

\begin{lem} \label{mix} (a) Let $R$ be a s-acyclic Borel subrelation of $(\mathbb{P}_f\!\times\!\mathbb{P}_\infty )\cup (\mathbb{P}_\infty\!\times\!\mathbb{P}_f)$. Then there is 
$f\! :\! 2^\omega\!\rightarrow\! 2^\omega$ injective continuous preserving $\mathbb{P}_f$ such that $(f\!\times\! f)^{-1}(R)\! =\!\emptyset$.\smallskip

(b) Let $R$ be a s-acyclic Borel subrelation of $\big( (2\!\times\!\mathbb{P}_f)\!\times\! ( 2\!\times\!\mathbb{P}_\infty )\big)\cup
\big( (2\!\times\!\mathbb{P}_\infty )\!\times\! (2\!\times\!\mathbb{P}_f)\big)$. Then there is 
$f\! :\! 2^\omega\!\rightarrow\! 2^\omega$ injective continuous preserving $\mathbb{P}_f$ such that 
$(f\!\times\! f)^{-1}(\bigcup_{\varepsilon ,\eta\in 2}~R_{\varepsilon ,\eta})\! =\!\emptyset$.\end{lem}

\noindent\bf Proof.\rm\ (a) By symmetry, we may assume that $R\!\subseteq\!\mathbb{P}_f\!\times\!\mathbb{P}_\infty$. Assume first that $R_{\alpha_0}$ is not meager for some 
$\alpha_0\!\in\!\mathbb{P}_f$. This gives $s\!\in\! 2^{<\omega}$ with the property that  $N_s\cap R_{\alpha_0}$ is comeager in $N_s$. In particular, 
$N_s\cap\mathbb{P}_f\!\setminus\{\alpha_0\}$ is not separable from $N_s\cap R_{\alpha_0}$ by a $\bormtwo$ set, by Baire's theorem. Theorem \ref{Lo-SR} provides $h\! :\! 2^\omega\!\rightarrow\! N_s\!\setminus\{\alpha_0\}$ injective continuous with 
$\mathbb{P}_f\!\subseteq\! h^{-1}(\mathbb{P}_f)$ and $\mathbb{P}_\infty\!\subseteq\! h^{-1}(R_{\alpha_0})$. If 
$R'\! :=\! (h\!\times\! h)^{-1}(R)$, $\alpha\!\in\!\mathbb{P}_f$, $\beta\not=\!\gamma\!\in\! R'_\alpha$, then 
$\big( h(\beta ),h(\alpha ),h(\gamma ),\alpha_0\big)$ is a $s(R)$-cycle, which is absurd. This shows that, replacing $R$ with $R'$ if necessary, we may assume that $R$ has meager vertical sections. Let $F$ be a meager $\boratwo$ subset of $2^\omega$ containing 
$\bigcup_{\alpha\in\mathbb{P}_f}~R_\alpha$. Note that $\mathbb{P}_f$ is not separable from 
$\mathbb{P}_\infty\!\setminus\! F$ by a $\bormtwo$ set. Theorem \ref{Lo-SR} provides 
$f\! :\! 2^\omega\!\rightarrow\! 2^\omega$ injective continuous with $\mathbb{P}_f\!\subseteq\! f^{-1}(\mathbb{P}_f)$ and 
$\mathbb{P}_\infty\!\subseteq\! f^{-1}(\mathbb{P}_\infty\!\setminus\! F)$. Thus $(f\!\times\! f)^{-1}(R)$ is empty.\medskip

\noindent (b) We argue as in the proof of Lemma \ref{finsq}.(b).\hfill{$\square$}

\begin{cor} \label{ac} (a) Let $R$ be a s-acyclic Borel relation on $2^\omega$. Then there is $f\! :\! 2^\omega\!\rightarrow\! 2^\omega$ injective continuous preserving $\mathbb{P}_f$ such that $(f\!\times\! f)^{-1}(R)\!\subseteq\!\Delta (2^\omega )$.\smallskip

(b) Let $R$ be a s-acyclic Borel relation on $\mathbb{D}$. Then there is $f\! :\! 2^\omega\!\rightarrow\! 2^\omega$ injective continuous preserving $\mathbb{P}_f$ such that 
$(f\!\times\! f)^{-1}(\bigcup_{\varepsilon ,\eta\in 2}~R_{\varepsilon ,\eta})\!\subseteq\!\Delta (2^\omega )$.\end{cor}

\noindent\bf Proof.\rm\ (a) By Lemma \ref{finsq}, we may assume that $R\cap\mathbb{P}_f^2\!\subseteq\!\Delta (2^\omega )$. By Lemma \ref{mix}, we may assume that 
$R\cap\big( (\mathbb{P}_f\!\times\!\mathbb{P}_\infty )\cup (\mathbb{P}_\infty\!\times\!\mathbb{P}_f)\big)\! =\!\emptyset$. By Lemma 
\ref{locacy}, we may assume that $R\!\setminus\!\Delta (2^\omega )$ is contained in the union of a sequence $(R_n)_{n\in\omega}$ of relations on 
$\mathbb{P}_\infty$ which are closed in $\mathbb{P}_\infty\!\times\! 2^\omega$ and in $2^\omega\!\times\!\mathbb{P}_\infty$. By Lemma \ref{inf} applied to $(R_n)_{n\in\omega}$, we may assume that $R\cap\mathbb{P}_\infty^2\!\subseteq\!\Delta (2^\omega )$.\medskip

\noindent (b) By Lemma \ref{finsq}, we may assume that $\bigcup_{\varepsilon ,\eta\in 2}~R_{\varepsilon ,\eta}\cap\mathbb{P}_f^2\!\subseteq\!\Delta (2^\omega )$. By Lemma \ref{mix}, we may assume that $R\cap\Big(\big( (2\!\times\!\mathbb{P}_f)\!\times\! ( 2\!\times\!\mathbb{P}_\infty )\big)\cup
\big( (2\!\times\!\mathbb{P}_\infty )\!\times\! (2\!\times\!\mathbb{P}_f)\big)\Big)\! =\!\emptyset$. By Lemma \ref{locacy}, we may assume that 
$(\bigcup_{\varepsilon ,\eta\in 2}~R_{\varepsilon ,\eta})\!\setminus\!\Delta (2^\omega )$ is contained in the union of a sequence 
$(R_n)_{n\in\omega}$ of relations on $\mathbb{P}_\infty$ which are closed in $\mathbb{P}_\infty\!\times\! 2^\omega$ and in 
$2^\omega\!\times\!\mathbb{P}_\infty$. By Lemma \ref{inf} applied to $(R_n)_{n\in\omega}$, we may assume that 
$\bigcup_{\varepsilon ,\eta\in 2}~R_{\varepsilon ,\eta}\cap\mathbb{P}_\infty^2\!\subseteq\!\Delta (2^\omega )$.\hfill{$\square$}\medskip

 Theorem \ref{genac} for classes of rank at least three is a consequence of Lemma \ref{antichac} and the following result since the elements of 
${\mathcal A}^{\bf\Gamma}$ are contained in either $\Delta (2^\omega )$, or in 
$\big\{\big( (\varepsilon ,x),(\eta ,x)\big)\!\in\!\mathbb{D}^2\mid x\!\in\! 2^\omega\big\}$, and the elements of ${\mathcal B}^{\bf\Gamma}$ are contained in 
$\Delta (\mathbb{S})\cup(\{ 0^\infty\}\!\times\! N_1)\cup (N_1\!\times\! \{ 0^\infty\} )$, which are s-acyclic and closed on the one side, and 
$\mathbb{G}^{\bf\Gamma}_m$ and $\mathbb{G}^{{\bf\Gamma},a}_m$ are the only graphs in 
${\mathcal A}^{\bf\Gamma}\cup {\mathcal B}^{\bf\Gamma}$ on the other side. The next proof is the last one using effective descriptive set theory.\medskip

\noindent\bf Notation.\rm\ We set, for any relation $R$ on $\mathbb{S}$, 
$R_{1,1}\! :=\!\{ (\alpha ,\beta )\!\in\! 2^\omega\!\times\! 2^\omega\mid (1\alpha ,1\beta )\!\in\! R\}$.

\begin{thm} \label{morethreeacy} Let ${\bf\Gamma}$ be a non self-dual Borel class of rank at least three, $X$ be an analytic space, and 
$R$ be a Borel relation on $X$ contained in a s-acyclic Borel relation with $\boratwo$ vertical sections. Exactly one of the following holds:\smallskip  

(a) the relation $R$ is a $\bf\Gamma$ subset of $X^2$,\smallskip  

(b) there is $(\mathbb{X},\mathbb{R})\!\in\! {\mathcal A}^{\bf\Gamma}\cup {\mathcal B}^{\bf\Gamma}$ such that 
$(\mathbb{X},\mathbb{R})\sqsubseteq_c(X,R)$.\end{thm}

\noindent\bf Proof.\rm\ By Theorem \ref{morethree} and since the elements of ${\mathcal B}^{\bf\Gamma}$ have a section not in 
$\bf\Gamma$, (a) and (b) cannot hold simultaneously. Assume that (a) does not hold. By Theorem \ref{contBr}, one of the following holds:\medskip  

\noindent (1) the relation $R$ has at least one section not in $\bf\Gamma$,\smallskip

\noindent (2) there is a relation $\mathcal{R}$ on $2^\omega$ such that 
$\mathcal{R}\cap\Delta (2^\omega )\! =\!\Delta (\mathbb{C})$ and $(2^\omega ,\mathcal{R})\sqsubseteq_c(X,R)$,\smallskip

\noindent (3) there is a relation $\mathcal{R}$ on $\mathbb{D}$ such that 
$\mathbb{R}_{0,1}\cap\Delta (2^\omega )\! =\!\Delta (\mathbb{C})$ and $(\mathbb{D},\mathcal{R})\sqsubseteq_c(X,R)$.

\vfill\eject

\noindent (1) Let $x_0\!\in\! X$ such that, for example, $R_{x_0}$ is not in $\bf\Gamma$, the other case being similar. Note that $R_{x_0}\!\setminus\!\{ x_0\}$ is not separable from $X\!\setminus\! (R_{x_0}\cup\{ x_0\})$ by a set in $\bf\Gamma$. Theorem 
\ref{Lo-SR} provides $h\! :\! 2^\omega\!\rightarrow\! X\!\setminus\!\{ x_0\}$ injective continuous with ${\mathbb{C}\! =\! h^{-1}(R_{x_0})}$. We define ${g\! :\!\mathbb{S}\!\rightarrow\! X}$ by $g(0^\infty )\! :=\! x_0$ and $g(1\alpha )\! :=\! h(\alpha )$, so that $g$ is injective continuous. Considering $(g\!\times\! g)^{-1}(R)$ if necessary, we may assume that $X\! =\!\mathbb{S}$, $x_0\! =\! 0^\infty$ and 
$R_{x_0}\cap N_1\! =\!\{ 1\alpha\mid\alpha\!\in\!\mathbb{C}\}$. Let $\mathcal{A}$ be a s-acyclic Borel relation with $\boratwo$ vertical sections containing $R$.\medskip

 For the simplicity of the notation, we assume that the rank of $\bf\Gamma$ is less than $\omega_1^{\text{CK}}$, and $\mathbb{C},\mathcal{A}$ are $\Borel$. Theorem 3.5 in [Lo1] gives a sequence $(\mathcal{C}_n)$ of $\Borel$ relations with closed vertical sections such that $\mathcal{A}\! =\!\bigcup_{n\in\omega}~\mathcal{C}_n$. By Lemma 2.2.2 in [L5], 
$\Borel\cap 2^\omega$ is countable and $\Ca$, so that $V\! :=\! 2^\omega\!\setminus\! (\Borel\cap 2^\omega )$ is $\Ana$, disjoint from $\Borel\cap 2^\omega$, and 
$V\cap\mathbb{C}$ is not separable from $V\!\setminus\!\mathbb{C}$ by a set in $\bf\Gamma$. We will apply Theorem 3.2 in [L5], where the 
{\bf Gandy-Harrington topology} ${\it\Sigma}_{2^\omega}$ on $2^\omega$ generated by $\Ana (2^\omega )$ is used. Let us prove that 
$\mathcal{A}_{1,1}\cap V^2$ is $({\it\Sigma}_{2^\omega})^2$-meager in $V^2$. It is enough to see that $(\mathcal{C}_n)_{1,1}\cap V^2$ is $({\it\Sigma}_{2^\omega})^2$-nowhere dense in $V^2$ for each $n$. By Lemma 3.1 in [L5], $(\mathcal{C}_n)_{1,1}$ is $({\it\Sigma}_{2^\omega})^2$-closed. We argue by contradiction, which gives $n$ and nonempty 
$\Ana$ subsets $S,T$ of $2^\omega$ with the property that $S\!\times\! T\!\subseteq\! (\mathcal{C}_n)_{1,1}\cap V^2$. By the effective perfect set Theorem (see 4.F1 in [Mo]), $S,T$ are uncountable. So pick $x,y\!\in\! S$ and $z,t\!\in\! T$ pairwise different. Then $(1x,1z,1y,1t)$ is a $s(\mathcal{A})$-cycle, which is absurd. Theorem 3.2 in [L5] provides 
$f\! :\! 2^\omega\!\rightarrow\! 2^\omega$ injective continuous preserving $\mathbb{C}$ with the property that $\big( f(\alpha ),f(\beta )\big)\!\notin\!\mathcal{A}_{1,1}$ if 
$\alpha\!\not=\!\beta$.\medskip

 Considering the set  $(f\!\times\! f)^{-1}(R)$ if necessary, we may assume that $R_{1,1}\!\subseteq\!\Delta (2^\omega )$. We set $E_{1,1}\! :=\!\{\alpha\!\in\! 2^\omega\mid (\alpha ,\alpha )\!\in\! R_{1,1}\}$, so that $E_{1,1}$ is a Borel subset of $2^\omega$ and $R_{1,1}\! =\!\Delta (E_{1,1})$. By Lemma \ref{four}, we may assume that $E_{1,1}\! =\! S_j$ for some $j\!\in\! 4$. If $j\! =\! 0$, then 
$\big( 2^\omega ,\Delta (\mathbb{C})\big)\sqsubseteq_c(X,R)$. So we may assume that $j\!\not=\! 0$. Similarly, 
$\{\alpha\!\in\! 2^\omega\mid (1\alpha ,0^\infty )\!\in\! R\}\! =\! S_j$ for some $j\!\in\! 4$. This provides $t\!\in\! A$ with 
$(\mathbb{S},\mathbb{R}^{{\bf\Gamma},a}_t)\sqsubseteq_c(X,R)$.\medskip

\noindent (2) We partially argue as in (1). Note that $\mathcal{R}$ is Borel and contained in a s-acyclic Borel relation with $\boratwo$ vertical sections $\mathcal{A}$.\medskip

 This time, $\mathcal{C}_n\!\in\!\Borel\big( (2^\omega )^2\big)$. Let us prove that 
$\mathcal{C}_n\cap V^2$ is $({\it\Sigma}_{2^\omega})^2$-nowhere dense in $V^2$ for each $n$. We argue by contradiction, which gives $n$ and nonempty $\Ana$ subsets $S,T$ of $2^\omega$ with $S\!\times\! T\!\subseteq\!\mathcal{C}_n\cap V^2$. Note that $(x,z,y,t)$ is a $s(\mathcal{A})$-cycle, which is absurd. Theorem 3.2 in [L5] provides 
$f\! :\! 2^\omega\!\rightarrow\! 2^\omega$ injective continuous preserving $\mathbb{C}$ such that $\big( f(\alpha ),f(\beta )\big)\!\notin\!\mathcal{A}$ if $\alpha\!\not=\!\beta$.\medskip
  
 Considering $(f\!\times\! f)^{-1}(\mathcal{R})$ if necessary, we may assume that $\mathcal{R}\!\subseteq\!\Delta (2^\omega )$, which means that 
$\big( 2^\omega ,\Delta (\mathbb{C})\big)\sqsubseteq_c(X,R)$.\medskip

\noindent (3) We partially argue as in (2).\medskip

 This time, $\mathcal{C}_n\!\in\!\Borel (\mathbb{D}^2)$. Fix $\varepsilon ,\eta\!\in\! 2$. Let us prove that 
$(\mathcal{C}_n)_{\varepsilon ,\eta}\cap V^2$ is $({\it\Sigma}_{2^\omega})^2$-nowhere dense in $V^2$ for each $n$. We argue by contradiction, which gives $n$ and nonempty 
$\Ana$ subsets $S,T$ of $2^\omega$ with $S\!\times\! T\!\subseteq\! (\mathcal{C}_n)_{\varepsilon ,\eta}\cap V^2$. Note that 
$\big( (\varepsilon ,x),(\eta ,z),(\varepsilon ,y),(\eta ,t)\big)$ is a $s(\mathcal{A})$-cycle, which is absurd. Theorem 3.2 in [L5] provides $f\! :\! 2^\omega\!\rightarrow\! 2^\omega$ injective continuous preserving $\mathbb{C}$ such that $\big( f(\alpha ),f(\beta )\big)\!\notin\!\bigcup_{\varepsilon ,\eta\in 2}~\mathcal{A}_{\varepsilon ,\eta}$ if $\alpha\!\not=\!\beta$.\medskip 

 Considering $(f\!\times\! f)^{-1}(\mathcal{R})$ if necessary, we may assume that $\mathcal{R}_{\varepsilon ,\eta}\!\subseteq\!\Delta (2^\omega )$ and 
$\mathcal{R}_{0,1}\! =\!\Delta (\mathbb{C})$. Theorem \ref{morethree} provides $(\mathbb{X},\mathbb{R})\!\in\! {\mathcal A}^{\bf\Gamma}$ such that 
$(\mathbb{X},\mathbb{R})\sqsubseteq_c(X,R)$.\medskip

 So (b) holds in any case.\hfill{$\square$}\medskip
 
 Theorem \ref{2ac} is an immediate consequence of the following result.

\begin{thm} \label{morethreeacry} Let ${\bf\Gamma}$ be a non self-dual Borel class of rank two, $X$ be an analytic space, and $R$ be a s-acyclic Borel relation on $X$. Exactly one of the following holds:\smallskip

(a) the relation $R$ is a $\bf\Gamma$ subset of $X^2$,\smallskip

(b) there is $(\mathbb{X},\mathbb{R})\!\in\! {\mathcal A}^{\bf\Gamma}\cup {\mathcal B}^{\bf\Gamma}$ such that 
$(\mathbb{X},\mathbb{R})\sqsubseteq_c(X,R)$.\end{thm}

\noindent\bf Proof.\rm\ We partially argue as in the proof of Theorem \ref{morethreeacy}. For the case (1), recall that we may assume that $X\! =\!\mathbb{S}$, $x_0\! =\! 0^\infty$ and $R_{x_0}\cap N_1\! =\!\{ 1\alpha\mid\alpha\!\in\!\mathbb{C}\}$. Corollary \ref{ac} provides $f\! :\! 2^\omega\!\rightarrow\! 2^\omega$ injective continuous preserving $\mathbb{C}$ such that 
$\big( f(\alpha ),f(\beta )\big)\!\notin\! R_{1,1}$ if $\alpha\!\not=\!\beta$. For the case (2), $\mathcal{R}$ is s-acyclic Borel, and by Corollary \ref{ac} we may assume that $\mathcal{R}\!\subseteq\!\Delta (2^\omega )$, which means that 
$\big( 2^\omega ,\Delta (\mathbb{C})\big)\sqsubseteq_c(X,R)$. For the case (3), by Corollary \ref{ac} we may assume that 
$\bigcup_{\varepsilon ,\eta\in 2}~R_{\varepsilon ,\eta}\!\subseteq\!\Delta (2^\omega )$ and $\mathcal{R}_{0,1}\! =\!\Delta (\mathbb{C})$.\hfill{$\square$}\medskip

 Theorem \ref{Sigma02ac} is an immediate consequence of the following result.

\begin{cor} \label{morethree+y} Let $X$ be an analytic space, and $R$ be a s-acyclic Borel relation on $X$ whose sections are in 
$\boratwo$. Exactly one of the following holds:\smallskip  

(a) the relation $R$ is a $\boratwo$ subset of $X^2$,\smallskip  

(b) there is $(\mathbb{X},\mathbb{R})\!\in\! {\mathcal A}^{\boratwo}$ such that $(\mathbb{X},\mathbb{R})\sqsubseteq_c(X,R)$.\smallskip

\noindent In particular, ${\mathcal A}^{\boratwo}$ is a 34 elements $\sqsubseteq_c$ and $\leq_c$-antichain basis.\end{cor}

\noindent\bf Proof.\rm\ We apply Theorem \ref{morethreeacry} and use the fact that the elements of ${\mathcal B}^{\bf\Gamma}$ have a section not in $\bf\Gamma$.\hfill{$\square$}\medskip

\noindent\bf Notation.\rm\ We set\medskip

\leftline{$\mathcal{C}^{\bormtwo}\! :=\!\mathcal{A}^{\bormtwo}~\cup$}\smallskip

\rightline{$\big\{ (\mathbb{S},\mathbb{R}^{\bormtwo ,a}_t)\mid t\!\in\! 4^{(2^2)}\wedge t(0,0),t(0,1),t(1,0)\!\in\! 2\wedge
\big( t(0,1)\! =\! 0\vee t(1,0)\! =\! 0\big)\wedge t(1,1)\!\not=\! 0\big\}$.}\medskip

 Theorem \ref{Pi02ac} is an immediate consequence of the following result.

\begin{cor} \label{morethreeacrypi} Let $X$ be an analytic space, and $R$ be a s-acyclic locally countable Borel relation on $X$. Exactly one of the following holds:\smallskip

(a) the relation $R$ is a $\bormtwo$ subset of $X^2$,\smallskip

(b) there is $(\mathbb{X},\mathbb{R})\!\in\! {\mathcal C}^{\bormtwo}$ such that 
$(\mathbb{X},\mathbb{R})\sqsubseteq_c(X,R)$.\smallskip

\noindent Moreover, ${\mathcal C}^{\bormtwo}$ is a 52 elements $\leq_c$-antichain (and thus a $\sqsubseteq_c$ and a $\leq_c$-antichain basis).\end{cor}

\noindent\bf Proof.\rm\ We apply Theorem \ref{morethreeacry} and use the fact that $\{\alpha\!\in\! 2^\omega\mid (1\alpha ,0^\infty )\!\in\! R\}\! =\! S_j$ for some $j\!\in\! 2$ since $R$ is locally countable. This provides $(\mathbb{S},\mathbb{R}^{\bormtwo ,a}_t)$ in ${\mathcal C}^{\bormtwo}$ below $(X,R)$.\hfill{$\square$}

\section{$\!\!\!\!\!\!$ The rank one} 

\noindent\bf Notation.\rm\ Let $\mathbb{K}\! :=\!\{ 2^{-k}\mid k\!\in\!\omega\}\cup\{ 0\}$, and 
$\mathbb{C}\! :=\!\{ 2^{-k}\mid k\!\in\!\omega\}$. We first set 
$$\begin{array}{ll}
& S_0\! :=\!\{ (x,y)\!\in\!\mathbb{K}^2\mid x,y\!\in\!\mathbb{C}\wedge x\! <\! y\}\mbox{,}\cr  
& S_1\! :=\!\{ (x,y)\!\in\!\mathbb{K}^2\mid x\! =\! y\!\in\!\mathbb{C}\}\mbox{,}\cr 
& S_2\! :=\!\{ (x,y)\!\in\!\mathbb{K}^2\mid x,y\!\in\!\mathbb{C}\wedge x\! >\! y\}\mbox{,}\cr    
& S_3\! :=\!\{ (x,y)\!\in\!\mathbb{K}^2\mid x\!\in\!\mathbb{C}\wedge y\!\notin\!\mathbb{C}\}\mbox{,}\cr   
& S_4\! :=\!\{ (x,y)\!\in\!\mathbb{K}^2\mid x\!\notin\!\mathbb{C}\wedge y\!\in\!\mathbb{C}\}\mbox{,}\cr 
& S_5\! :=\!\{ (x,y)\!\in\!\mathbb{K}^2\mid x,y\!\notin\!\mathbb{C}\}.
\end{array}$$ 

 Note that $(S_j)_{j<6}$ is a partition of $\mathbb{K}^2$, and the vertical sections of $S_2$ and the horizontal sections of $S_0$ are infinite. We set $N\! :=\!\big\{ t\!\in\! 2^6\mid\big( t\!\notin\! 1^6\wedge t(5)\! =\! 0\big)\vee\big( t(2)\! =\! 1\wedge t(3)\! =\! 0\big)\vee\big( t(0)\! =\! 1\wedge t(4)\! =\! 0\big)\big\}$. We first consider the class $\bormone$ and set, for $t\!\in\! 2^6$, 
$\mathbb{R}^{\bormone}_{N,t}\! :=\!\bigcup_{j<6,t(j)=1}~S_j$.\medskip

 We next code relations having just one non-closed vertical section, with just one limit point on this vertical section, out of the diagonal. We set\medskip
  
\leftline{$V\! :=\!\big\{ t\!\in\! (2^6)^{2^2}\mid t(0,0)\!\in\! (\{ 0\}^5\!\times\! 2)\wedge 
t(0,1)\! =\! (0,0,0,0,1,0)\wedge t(1,0)\!\in\! (\{ 0\}^3\!\times\! 2\!\times\!\{ 0\}\!\times\! 2)$}\smallskip

\rightline{$\wedge\ t(1,1)\!\notin\! N\big\} ,$}\medskip

\noindent and $\mathbb{L}\! :=\! (\neg\mathbb{C})\oplus\mathbb{K}$. We set, for $t\!\in\! V$, 
$$\mathbb{R}^{\bormone}_{V,t}\! :=\!\bigcup_{(\varepsilon ,\eta )\in 2^2,j<6,t(\varepsilon ,\eta )(j)=1}~
\big\{\big( (\varepsilon ,x),(\eta ,y)\big)\!\in\!\mathbb{L}^2\mid (x,y)\!\in\! S_j\big\} .$$
Similarly, we code relations having closed vertical sections, and just one non-closed horizontal section, with just one limit point on this horizontal section, out of the diagonal. We set\medskip
  
\leftline{$H\! :=\!\big\{ t\!\in\! (2^6)^{2^2}\mid t(0,0)\!\notin\! N\wedge 
t(0,1)\! =\! (0,0,0,1,0,0)\wedge\ t(1,0)\!\in\! (\{ 0\}^4\!\times\! 2^2)\!\setminus\!\{ (0,0,0,0,1,0)\}$}\smallskip

\rightline{$\wedge\ t(1,1)\!\in\! (\{ 0\}^5\!\times\! 2)\big\} ,$}

\noindent and $\mathbb{M}\! :=\!\mathbb{K}\oplus(\neg\mathbb{C})$. We set, for $t\!\in\! H$, 
$$\mathbb{R}^{\bormone}_{H,t}\! :=\!\bigcup_{(\varepsilon ,\eta )\in 2^2,j<6,t(\varepsilon ,\eta )(j)=1}~
\big\{\big( (\varepsilon ,x),(\eta ,y)\big)\!\in\!\mathbb{M}^2\mid (x,y)\!\in\! S_j\big\} .$$
We define a set of codes for relations on $\mathbb{K}$ with closed sections as follows: 
$$C\! :=\!\{ t\!\in\! 2^6\mid t(0)\! =\! 1\Rightarrow t(4)\! =\! 1\Rightarrow t(5)\! =\! 1\wedge t(2)\! =\! 1\Rightarrow t(3)\! =\! 1\Rightarrow t(5)\! =\! 1\} .$$ 
We define a set of codes for the missing relations in our antichain basis. We set\medskip

\leftline{$S\! :=\!\big\{ t\!\in\! (2^6)^{2^2}\mid t(0,0),t(1,1)\!\notin\! N\wedge t(0,1)\! =\! (0,1,0,0,0,0)\wedge 
t(1,0)\!\in\! C\ \wedge$}\smallskip

\rightline{$t(1,0)\! =\! (0,1,0,0,0,0)\Rightarrow t(0,0)\!\leq_{\text{lex}}\!t(1,1)\big\} .$}\medskip

\noindent We set, for $t\!\in\! S$, 
$\mathbb{R}^{\bormone}_{S,t}\! :=\!\bigcup_{(\varepsilon ,\eta )\in 2^2,j<6,t(\varepsilon ,\eta )(j)=1}~
\big\{\big( (\varepsilon ,x),(\eta ,y)\big)\!\in\!\mathbb{D}^2\mid (x,y)\!\in\! S_j\big\}$. Finally, 
$$\mathcal{A}^{\bormone}\! :=\!\{ (\mathbb{K},\mathbb{R}^{\bormone}_{N,t})\mid t\!\in\! N\}\cup
\{ (\mathbb{L},\mathbb{R}^{\bormone}_{V,t})\mid t\!\in\! V\}\cup\{ (\mathbb{M},\mathbb{R}^{\bormone}_{H,t})\mid t\!\in\! H\}\cup\{ (\mathbb{D},\mathbb{R}^{\bormone}_{S,t})\mid t\!\in\! S\}$$ 
is the 7360 elements $\leq_c$-antichain basis mentioned in the statement of Theorem \ref{one}.\medskip

 We set, for $j\!\in\! 2^2$, $T_j\! :=\!\{ (2^{-2k-j(0)},2^{-2k-j(1)})\mid k\!\in\!\omega\}$. We then define relations on $\mathbb{K}$ by  
$\mathbb{R}^{\bormone}_0\! :=\! T_{(0,1)}$ and $\mathbb{R}^{\bormone}_1\! :=\! T_{(0,1)}\cup T_{(1,0)}$. They describe the 2 elements mentioned at the end of the statement of Theorem \ref{one}. For the class $\boraone$, we simply pass to complements.\medskip

 For graphs, in the case of the class $\bormone$, the 5-elements $\leq_c$-antichain basis is described by the following codes:\medskip
 
\noindent - $(0,0,0,1,1,0)\mbox{ (acyclic) },(1,0,1,0,0,0),(1,0,1,1,1,0)\!\in\! N$,\smallskip

\noindent - $\big( (0,0,0,0,0,0),(0,0,0,0,1,0),(0,0,0,1,0,0),(0,0,0,0,0,0)\big)\!\in\! V$ (acyclic),\smallskip 

\noindent - $\big( (0,0,0,0,0,0),(0,1,0,0,0,0),(0,1,0,0,0,0),(0,0,0,0,0,0)\big)\!\in\! S$ (acyclic).

\vfill\eject

 In order to get the 6-elements $\sqsubseteq_c$-antichain basis, we just have to add $(\mathbb{K},\mathbb{R}^{\bormone}_1)$ (which is acyclic). In the case of the class 
$\boraone$, the 10-elements $\leq_c$ and $\sqsubseteq_c$-antichain basis is described as follows: take\medskip

\noindent - for $(1,1,1,0,0,1)\!\in\! N$, $\cup_{j<6,t(j)=0}~S_j$ (which is acyclic),\smallskip

\noindent - for 
$$t\!\in\!\big\{\big( (0,0,0,0,0,1),(0,0,0,0,1,0),(0,0,0,1,0,0),(\varepsilon_0,1,\varepsilon_0,\varepsilon_1,\varepsilon_1,1)\big)\!\in\! V\mid (\varepsilon_0,\varepsilon_1)\!\not=\! (1,0)\big\}\mbox{,}$$
$\cup_{(\varepsilon ,\eta )\in 2^2,j<6,t(\varepsilon ,\eta )(j)=0}~
\big\{\big( (\varepsilon ,x),(\eta ,y)\big)\!\in\!\mathbb{L}^2\mid (x,y)\!\in\! S_j\big\}$ (which is acyclic when $(\varepsilon_0,\varepsilon_1)\! =\! (1,1)$),\medskip 

\noindent - for\medskip

\leftline{$t\!\in\!\big\{\big( (\varepsilon_0,1,\varepsilon_0,\varepsilon_1,\varepsilon_1,1),(0,1,0,0,0,0),(0,1,0,0,0,0),(\varepsilon_2,1,\varepsilon_2,\varepsilon_3,\varepsilon_3,1)\big)\!\in\! S\mid$}\smallskip

\hfill{$(\varepsilon_0,\varepsilon_1),(\varepsilon_2,\varepsilon_3)\!\not=\! (1,0)\wedge (\varepsilon_0,1,\varepsilon_0,\varepsilon_1,\varepsilon_1,1)\!\leq_{\text{lex}}\!(\varepsilon_2,1,\varepsilon_2,\varepsilon_3,\varepsilon_3,1)\big\}$.}\medskip

\noindent $\cup_{(\varepsilon ,\eta )\in 2^2,j<6,t(\varepsilon ,\eta )(j)=0}~
\big\{\big( (\varepsilon ,x),(\eta ,y)\big)\!\in\!\mathbb{D}^2\mid (x,y)\!\in\! S_j\big\}$.

\section{$\!\!\!\!\!\!$ Uncountable analytic relations}

\bf Notation.\rm\ We set ${{\mathbb H}\! :=\! N_1\!\times\!\{ 0^\infty\}}$, 
${\mathbb V}\! :=\!\{ 0^\infty\}\!\times\! N_1$, and $\mathbb{L}\! :=\! {\mathbb H}\cup {\mathbb V}$. If 
$\mathbb{A}\!\in\!\{ {\mathbb H},{\mathbb V},{\mathbb L}\}$, then $\mathbb{A}^+\! :=\!\mathbb{A}\cup\{ (0^\infty ,0^\infty )\}$. Let $o\! :\! 2^\omega\!\rightarrow\! 2^\omega$ be defined by $o(\alpha )(n)\! :=\!\alpha (n)$ exactly when $n\! >\! 0$. Then $o$ is a homeomorphism and an involution. We set\medskip
 
\leftline{$\mathcal{A}^c\! :=\!\big\{\big( 2^\omega ,(2^\omega )^2\big) ,({\mathbb S},{\mathbb H}),
({\mathbb S},{\mathbb V}),({\mathbb S},{\mathbb L}),({\mathbb S},{\mathbb H}^+),
({\mathbb S},{\mathbb V}^+),({\mathbb S},{\mathbb L}^+),$}\smallskip

\rightline{$(2^\omega ,\not= ),(2^\omega ,<_{\text{lex}}),(2^\omega ,=),(2^\omega ,\leq_{\text{lex}}),
\big( 2^\omega ,\mbox{Graph}(o)\big),\big( 2^\omega ,\mbox{Graph}(o_{\vert N_0})\big)\big\} .$}\medskip

\noindent We enumerate $\mathcal{A}^c\! :=\!\{ {\cal E}_i\mid i\!\leq\! 12\}$ (in the previous order).

\begin{lem} $\mathcal{A}^c$ is an antichain.\end{lem}

\noindent\bf Proof.\rm\ Note that if $(X,A)\sqsubseteq_c(Y,B)$ and $B$ is in some Borel class $\bf\Gamma$, then $A$ is in $\bf\Gamma$ too. The second coordinate of a member of ${\cal E}_0$-${\cal E}_6$ (resp., 
${\cal E}_7$-${\cal E}_8$, ${\cal E}_9$-${\cal E}_{12}$) is clopen (resp., open not closed, closed not open). This proves that no member of ${\cal E}_7$-${\cal E}_{12}$ is reducible to a member of 
${\cal E}_0$-${\cal E}_6$, and that the members of ${\cal E}_7$-${\cal E}_8$ are incomparable to the members of ${\cal E}_9$-${\cal E}_{12}$. Note that the second coordinate of\smallskip

\noindent - ${\cal E}_0$ and ${\cal E}_9$-${\cal E}_{10}$ is reflexive.\smallskip

\noindent - ${\cal E}_1$-${\cal E}_3$, ${\cal E}_7$-${\cal E}_8$ and ${\cal E}_{11}$-${\cal E}_{12}$ is irreflexive.\smallskip

\noindent - ${\cal E}_0$, ${\cal E}_3$, ${\cal E}_6$, ${\cal E}_7$, ${\cal E}_9$ and ${\cal E}_{11}$ is symmetric.\smallskip

\noindent - ${\cal E}_1$-${\cal E}_2$, ${\cal E}_4$-${\cal E}_5$, ${\cal E}_8$-${\cal E}_{10}$ and 
${\cal E}_{12}$ is antisymmetric.\smallskip

\noindent - ${\cal E}_0$-${\cal E}_2$, ${\cal E}_4$-${\cal E}_5$, ${\cal E}_8$-${\cal E}_{10}$ and 
${\cal E}_{12}$ is transitive.\medskip

 Assume that $(X,A)\sqsubseteq_c(Y,B)$, and that $P$ is one of the following properties of relations: reflexive, irreflexive, symmetric, antisymmetric, transitive. We already noticed that $(X,A)$ has $P$ if $(Y,B)$ does. Note also that if $(X,A)$ has $P$, then there is a copy $C$ of $X$ in $Y$ such that $(C,B\cap C^2)$ has $P$. This implies that the only cases to consider are the following. In all these cases, we will prove by contradiction a result of the form $(X,A)\not\sqsubseteq_c(Y,B)$, which gives $i\! :\! X\!\rightarrow\! Y$ injective continuous with $A\! =\! (i\!\times\! i)^{-1}(B)$.\medskip

 ${\cal E}_1\not\sqsubseteq_c{\cal E}_2$: $i(10^\infty )\! =\! i(1^\infty )\! =\! 0^\infty$, which contradicts the injectivity of $i$. Similarly, ${\cal E}_2\not\sqsubseteq_c{\cal E}_1$, ${\cal E}_4$ is incomparable with ${\cal E}_5$, 
${\cal E}_1\not\sqsubseteq_c{\cal E}_5$, ${\cal E}_2\not\sqsubseteq_c{\cal E}_4$.\medskip  

 ${\cal E}_1\not\sqsubseteq_c{\cal E}_4$: $i(0^\infty )\! =\! 0^\infty$, so that $(0^\infty ,0^\infty )\!\in\! {\mathbb H}$, which is absurd. Similarly, ${\cal E}_2\not\sqsubseteq_c{\cal E}_5$ and ${\cal E}_3\not\sqsubseteq_c{\cal E}_6$.\medskip

 ${\cal E}_1\not\sqsubseteq_c{\cal E}_8$: for example $i(10^\infty )\! <_{\text{lex}}\! i(1^\infty )$ and 
 $(10^\infty ,1^\infty )\!\in\!\mathbb{H}$, which is absurd. Similarly, ${\cal E}_2\not\sqsubseteq_c{\cal E}_8$.\medskip

 ${\cal E}_1\not\sqsubseteq_c{\cal E}_{12}$: $\big( i(10^\infty ),i(0^\infty )\big)\! =\! (0\alpha ,1\alpha )$, 
$\big( i(1^\infty ),i(0^\infty )\big)\! =\! (0\beta ,1\beta )$, so that $\alpha\! =\! \beta$, which contradicts the injectivity of $i$. Similarly, ${\cal E}_2\not\sqsubseteq_c{\cal E}_{12}$.\medskip

 ${\cal E}_3\not\sqsubseteq_c{\cal E}_7$: $i(10^\infty )\!\not=\! i(1^\infty )$, so that $(10^\infty ,1^\infty )\!\in\! {\mathbb L}$, which is absurd.\medskip

 ${\cal E}_3\not\sqsubseteq_c{\cal E}_{11}$: 
$\big( i(0^\infty ),i(1^\infty )\big),\big( i(0^\infty ),i(10^\infty )\big)\!\in\!\mbox{Graph}(o)$, so that 
$i(1^\infty )\! =\! i(10^\infty )$, which contradicts the injectivity of $i$.\hfill{$\square$}\medskip

 From now on, $Y$ will be a Hausdorff topological space and $B$ will be an uncountable analytic relation on $Y$. Note that $B\cap\Delta (Y)$ is analytic.

\begin{lem} \label{52} Assume that $B\cap\Delta (Y)$ is uncountable. Then $\big( 2^\omega ,(2^\omega )^2\big)\sqsubseteq_c(Y,B)$, 
${(2^\omega ,=)\sqsubseteq_c(Y,B)}$ or $(2^\omega ,\leq_{\text{lex}})\sqsubseteq_c(Y,B)$.\end{lem}

\noindent\bf Proof.\rm\ The perfect set theorem gives $j\! :\! 2^\omega\!\rightarrow\! B\cap\Delta (Y)$ injective continuous. Note that 
$$\pi\! :=\!\mbox{proj}_0\big[ j[2^\omega ]\big]\! =\!\mbox{proj}_1\big[ j[2^\omega ]\big]$$ 
is a copy of $2^\omega$, $\Delta (\pi )\!\subseteq\! B\cap\pi^2$ and $(\pi ,B\cap\pi^2)\sqsubseteq_c(Y,B)$, so that we may assume that $Y\! =\! 2^\omega$ and $\Delta (2^\omega )\!\subseteq\! B\!\in\!\ana\big( (2^\omega )^2\big)$. By 19.7 in [K], there is a copy $P$ of $2^\omega$ in $2^\omega$ such that $<_{\text{lex}}\cap P^2\!\subseteq\! B$ or 
$<_{\text{lex}}\cap P^2\!\subseteq\!\neg B$. Similarly, there is a copy $Q$ of $2^\omega$ in $P$ such that 
$>_{\text{lex}}\cap Q^2\!\subseteq\! B$ or 
$>_{\text{lex}}\cap Q^2\!\subseteq\!\neg B$.\medskip

\noindent\bf Case 1.\rm\ $<_{\text{lex}}\cap Q^2\!\subseteq\! B$ and $>_{\text{lex}}\cap Q^2\!\subseteq\! B$.\medskip

 Note that $Q^2\! =\! B\cap Q^2$ and $\big( 2^\omega ,(2^\omega )^2\big)\sqsubseteq_c(Y,B)$.\medskip

\noindent\bf Case 2.\rm\ $<_{\text{lex}}\cap Q^2\!\subseteq\!\neg B$ and $>_{\text{lex}}\cap Q^2\!\subseteq\!\neg B$.\medskip

 Note that $\Delta (Q)\! =\! B\cap Q^2$ and $(2^\omega ,=)\sqsubseteq_c(Y,B)$.\medskip

\noindent\bf Case 3.\rm\ $<_{\text{lex}}\cap Q^2\!\subseteq\! B$ and $>_{\text{lex}}\cap Q^2\!\subseteq_c\!\neg B$.\medskip

 Note that $\leq_{\text{lex}}\cap Q\! =\! B\cap Q^2$ and $(2^\omega ,\leq_{\text{lex}})\sqsubseteq_c(Y,B)$.\medskip

\noindent\bf Case 4.\rm\ $<_{\text{lex}}\cap Q^2\!\subseteq\!\neg B$ and $>_{\text{lex}}\cap Q^2\!\subseteq\! B$.\medskip

 Note that $\geq_{\text{lex}}\cap Q^2\! =\! B\cap Q^2$ and $(2^\omega ,\geq_{\text{lex}})\sqsubseteq_c(Y,B)$. But 
$(2^\omega ,\leq_{\text{lex}})\sqsubseteq_c(2^\omega ,\geq_{\text{lex}})$, with witness $i$ defined by 
$i(\alpha )(n)\! :=\! 1\! -\!\alpha (n)$. Thus $(2^\omega ,\leq_{\text{lex}})\sqsubseteq_c(Y,B)$.\hfill{$\square$}

\begin{lem} Assume that there is $C\!\subseteq\! Y$ countable such that 
$B\!\subseteq\! (C\!\times\! Y)\cup (Y\!\times\! C)$. Then there is $1\!\leq\! i\!\leq\! 6$ such that  
${\cal E}_i\sqsubseteq_c(Y,B)$.\end{lem}

\noindent\bf Proof.\rm\ As $B$ is uncountable, there is $y\!\in\! C$ such that $B^y$ or $B_y$ is uncountable. As $B^y$ and $B_y$ are analytic, there is a copy $P$ of $2^\omega$ in $Y$, disjoint from $C$, such that 
$P\!\times\!\{ y\}\!\subseteq\! B$ or $\{ y\}\!\times\! P\!\subseteq\! B$.\medskip

\noindent\bf Case 1.\rm\ $P\!\times\!\{ y\}\!\subseteq\! B$.\medskip

\noindent\bf Case 1.1.\rm\ $(\{ y\}\!\times\! P)\cap B$ is countable.\medskip

 Note that there is a copy $Q$ of $2^\omega$ in $P$ such that $\{ y\}\!\times\! Q\!\subseteq\!\neg B$.

\vfill\eject 
 
\noindent\bf Case 1.1.1.\rm\ $(y,y)\!\notin\! B$.\medskip

 Note that $Q\!\times\!\{ y\}\! =\! B\cap (\{ y\}\cup Q)^2$ and 
$(\mathbb{S},\mathbb{H})\sqsubseteq_c(Y,B)$.\medskip

\noindent\bf Case 1.1.2.\rm\ $(y,y)\!\in\! B$.\medskip

 Note that $Q\!\times\!\{ y\}\cup\{ (y,y)\}\! =\! B\cap (\{ y\}\cup Q)^2$ and $(\mathbb{S},\mathbb{H}^+)\sqsubseteq_c(Y,B)$.\medskip

\noindent\bf Case 1.2.\rm\ $(\{ y\}\!\times\! P)\cap B$ is uncountable.\medskip

 Note that there is a copy $Q$ of $2^\omega$ in $P$ such that $\{ y\}\!\times\! Q\!\subseteq\! B$.\medskip

\noindent\bf Case 1.2.1.\rm\ $(y,y)\!\notin\! B$.\medskip

 Note that $Q\!\times\!\{ y\}\cup\{ y\}\!\times\! Q\! =\! B\cap (\{ y\}\cup Q)^2$ and 
$(\mathbb{S},\mathbb{L})\sqsubseteq_c(Y,B)$.\medskip

\noindent\bf Case 1.2.2.\rm\ $(y,y)\!\in\! B$.\medskip

 Note that $Q\!\times\!\{ y\}\cup\{ y\}\!\times\! Q\cup\{ (y,y)\}\! =\! B\cap (\{ y\}\cup Q)^2$ and 
$(\mathbb{S},\mathbb{L}^+)\sqsubseteq_c(Y,B)$.\medskip

\noindent\bf Case 2.\rm\ $\{ y\}\!\times\! P\!\subseteq\! B$.\medskip

 Similarly, we show that ${\cal E}_2\sqsubseteq _c(Y,B)$, ${\cal E}_3\sqsubseteq_c(Y,B)$, 
${\cal E}_5\sqsubseteq_c(Y,B)$ or ${\cal E}_6\sqsubseteq_c(Y,B)$.\hfill{$\square$}\medskip

 So from now on we will assume that $B\cap\Delta (Y)$ is countable, and that there is no countable subset $C$ of $Y$ such that 
$B\!\subseteq\! (C\!\times\! Y)\cup (Y\!\times\! C)$. In particular, we may assume that $B$ is irreflexive. By Theorem 1 and Remark 2 in [P], there are $\varphi\! :\! 2^\omega\!\rightarrow\! Y$ and $h\! :\!\varphi [2^\omega ]\!\rightarrow\! Y$ injective continuous with 
$\mbox{Graph}(h)\!\subseteq\! B$. As $B$ is irreflexive, we may assume that $h$ has disjoint domain and range. We define 
$i\! :\! 2^\omega\!\rightarrow\! Y$ by $i(0\alpha )\! :=\!\varphi (\alpha )$ and $i(1\alpha )\! :=\! h\big(\varphi (\alpha )\big)$, so that $i$ is injective continuous. We set $A\! :=\! (i\!\times\! i)^{-1}(B)$. Note that $A$ is an analytic digraph on $2^\omega$, which contains 
$\mbox{Graph}(o_{\vert N_0})$, and that $(2^\omega ,A)\sqsubseteq_c(Y,B)$. So from now on we will assume that $Y\! =\! 2^\omega$, 
$B\!\in\!\ana\big( (2^\omega )^2\big)$ is a digraph, and $\mbox{Graph}(o_{\vert N_0})\!\subseteq\! B$. 

\begin{lem} Assume that $B$ is meager. Then there is $11\!\leq\! i\!\leq\! 12$ such that  
${\cal E}_i\sqsubseteq_c(Y,B)$.\end{lem}

\noindent\bf Proof.\rm\ It is enough to find a Cantor subset $P$ of $N_0$ such that 
$B\cap (P\cup o[P])^2\!\subseteq\!\mbox{Graph}(o)$. In order to see this, we distinguish two cases.\medskip

\noindent\bf Case 1.\rm\ $\mbox{Graph}(o_{\vert o[P]})\cap B$ is uncountable.\medskip

 There is a copy $Q$ of $2^\omega$ in $o[P]$ with $\mbox{Graph}(o_{\vert Q})\!\subseteq\! B$, so that 
$\mbox{Graph}(o_{\vert Q\cup o[Q]})\! =\! B\cap (Q\cup o[Q])^2$. Let 
$\psi\! :\! 2^\omega\!\rightarrow\! Q$ be a homeomorphism. We define $j\! :\! 2^\omega\!\rightarrow\! 2^\omega$ by the formulas $j(0\alpha )\! :=\!\psi (\alpha )$ and $j(1\alpha )\! :=\! o\big(\psi (\alpha )\big)$, so that 
$j$ is injective continuous. Note that $\mbox{Graph}(o)\! =\! (j\!\times\! j)^{-1}(B)$, so that 
$\big( 2^\omega ,\mbox{Graph}(o)\big)\sqsubseteq_c(Y,B)$.\medskip

\noindent\bf Case 2.\rm\ $\mbox{Graph}(o_{\vert o[P]})\cap B$ is countable.\medskip

 There is a copy $Q$ of $2^\omega$ in $o[P]$ with $\mbox{Graph}(o_{\vert Q})\!\subseteq\!\neg B$, so that 
$\mbox{Graph}(o_{\vert o[Q]})\! =\! B\cap (Q\cup o[Q])^2$. Let 
$\psi\! :\! 2^\omega\!\rightarrow\! o[Q]$ be a homeomorphism. We define 
$j\! :\! 2^\omega\!\rightarrow\! 2^\omega$ by the formulas $j(0\alpha )\! :=\!\psi (\alpha )$ and 
$j(1\alpha )\! :=\! o\big(\psi (\alpha )\big)$, so that $j$ is injective continuous. Note that 
$\mbox{Graph}(o_{\vert N_0})\! =\! (j\!\times\! j)^{-1}(B)$, so that 
$\big( 2^\omega ,\mbox{Graph}(o_{\vert N_0})\big)\!\sqsubseteq\! (Y,B)$.\medskip

 As $B$ is meager, there is  a sequence $(K_m)_{m\in\omega}$ of meager compact subsets of $(2^\omega )^2$ satisfying the inclusions  
$B\!\setminus\!\mbox{Graph}(o)\!\subseteq\!\bigcup_{m\in\omega}~K_m\!\subseteq\!\neg\mbox{Graph}(o)$.\medskip

 We build a sequence $(U_s)_{s\in 2^{<\omega}}$ of nonempty clopen subsets of 
$N_0$ satisfying the following:
$$\begin{array}{ll}
& (1)~U_{s\varepsilon}\!\subseteq\! U_s\cr
& (2)~\mbox{diam}(U_s)\!\leq\! 2^{-\vert s\vert}\cr
& (3)~U_{s0}\cap U_{s1}\! =\!\emptyset\cr
& (4)~\forall (f,g)\!\in\!\{\mbox{Id}_{2^\omega},o\}^2~~\forall (s,t)\!\in\! (2\!\times\! 2)^{<\omega}~~
(f[U_s]\!\times\! g[U_t])\cap (\bigcup_{m<\vert s\vert}~K_m)\! =\!\emptyset
\end{array}$$
Assume that this is done. We define $\varphi\! :\! 2^\omega\!\rightarrow\! N_0$ by 
$\{\varphi (\alpha )\}\! :=\!\bigcap_{n\in\omega}~U_{\alpha\vert n}$. Note that $\varphi$ is injective and continuous, so that $P\! :=\!\varphi [2^\omega ]$ is a Cantor subset of $N_0$. It remains to check the inclusion 
$B\cap (P\cup o[P])^2\!\subseteq\!\mbox{Graph}(o)$. So let 
$(\alpha ,\beta )\!\in\! B\cap (P\cup o[P])^2\!\setminus\!\mbox{Graph}(o)$, and $m$ with 
$(\alpha ,\beta)\!\in\! K_m$. We get $s,t\!\in\! 2^{m+1}$ such that 
$(\alpha ,\beta )\!\in\!\big( (U_s\cup o[U_s])\!\times\! (U_t\cup o[U_t])\big)\cap K_m$, which is absurd.\medskip

 Let us prove that the construction is possible. We first set $U_\emptyset\! :=\! N_0$. Note that 
$$\{ (\alpha_0,\alpha_1)\!\in\! U_\emptyset^2\mid\alpha_0\!\not=\!\alpha_1\mbox{ and }
\forall (f,g)\!\in\!\{\mbox{Id}_{2^\omega},o\}^2~~\forall (s,t)\!\in\! 2\!\times\! 2~~
\big( f(\alpha_s),g(\alpha_t)\big)\!\notin\! K_0\}$$ 
is a dense open subset of $U_\emptyset^2$. In particular, it is not empty, and we can pick $(\alpha_0,\alpha_1)$ in it. We choose a clopen neighbourhood $U_\varepsilon\!\subseteq\! U_\emptyset$ of $\alpha_\varepsilon$ with diameter at most $2^{-1}$ such that $U_0\cap U_1\! =\!\emptyset$ and 
$\bigcup_{(s,t)\in 2\times 2}~U_s\!\times\! U_t\!\subseteq\!\bigcap_{(f,g)\in\{\mbox{Id}_{2^\omega},o\}^2}~(f\!\times\! g)^{-1}(\neg K_0)$. Assume that $(U_s)_{\vert s\vert\leq l}$ satisfying conditions (1)-(4) have been constructed, which is the case for $l\! =\! 1$, so that from now on $l\!\geq\! 1$. Note that\medskip

\leftline{$\big\{ (\alpha_s)_{s\in 2^{l+1}}\!\in\! (2^\omega )^{2^{l+1}}\mid\forall s\!\in\! 2^{l+1}~
\alpha_s\!\in\! U_{s\vert l}\mbox{ and }\forall u\!\in\! 2^l~\alpha_{u0}\!\not=\!\alpha_{u1}\mbox{ and }$}\smallskip

\rightline{$\forall (f,g)\!\in\!\{\mbox{Id}_{2^\omega},o\}^2~\forall (s,t)\!\in\! 2^{l+1}\!\times\! 2^{l+1}~
\big( f(\alpha_s),g(\alpha_t)\big)\!\notin\!\bigcup_{m\leq l}~K_m\big\}$}\medskip

\noindent is a dense open subset of $\Pi_{s\in 2^{l+1}}~U_{s\vert l}$. We pick $(\alpha_s)_{s\in 2^{l+1}}$ in it, and choose a clopen neighbourhood $U_s\!\subseteq\! U_{s\vert l}$ of $\alpha_s$ with diameter at most 
$2^{-l-1}$ such that $U_{u0}\cap U_{u1}\! =\!\emptyset$ and 
$$\bigcup_{(s,t)\in 2^{l+1}\times 2^{l+1}}~U_s\!\times\! U_t\!\subseteq\!
\bigcap_{(f,g)\in\{\mbox{Id}_{2^\omega},o\}^2,m\leq l}~(f\!\times\! g)^{-1}(\neg K_m).$$
This finishes the proof.\hfill{$\square$}\medskip

 So we may assume that $B$ is not meager. The Baire property of $B$ and 19.6 in [K] give a product of Cantor sets contained in $B$. This means that we may assume that 
$N_0\!\times\! N_1\!\subseteq\! B\!\subseteq\!\neg\Delta (2^\omega )$.\medskip

\noindent\bf Proof of Theorem \ref{cou}.\rm\ (1) We distinguish several cases.\medskip

\noindent\bf Case 1.\rm\ $B\cap (N_1\!\times\! N_0)$ is not meager.\medskip

 By 19.6 in [K], $B\cap (N_1\!\times\! N_0)$ contains a product of Cantor sets, so that we may assume that 
$(N_0\!\times\! N_1)\cup (N_1\!\times\! N_0)\!\subseteq\! B$.\medskip
 
\noindent\bf Case 1.1.\rm\ There is a Cantor subset of $2^\omega$ which is 
$B$-discrete. Then $(\mathbb{S},\mathbb{L})\sqsubseteq_c(Y,B)$.\medskip

 Indeed, assume for example that $Q\!\subseteq\! N_1$ is a Cantor $B$-discrete set. Let $h\! :\! 2^\omega\!\rightarrow\! Q$ be a homeomorphism. We define $i\! :\!\mathbb{S}\!\rightarrow\! Y$ by $i(0^\infty )\! :=\! 0^\infty$ and $i(1\alpha )\! :=\! h(\alpha )$, so that 
$i$ is injective continuous. Clearly $\mathbb{L}\!\subseteq\! (i\!\times\! i)^{-1}(B)$, and the converse holds since $B$ is a digraph and 
$Q$ is $B$-discrete.

\vfill\eject

\noindent\bf Case 1.2.\rm\ No Cantor subset of $2^\omega$ is $B$-discrete. Then 
$(2^\omega ,\not= )\sqsubseteq_c(Y,B)$ or $(2^\omega ,<_{\text{lex}})\sqsubseteq_c(Y,B)$.\medskip

 Indeed, as in the proof of Lemma \ref{52} there is a Cantor subset $Q$ of $2^\omega$ with 
$<_{\text{lex}}\cap Q^2\!\subseteq\! B$ or $<_{\text{lex}}\cap Q^2\!\subseteq\!\neg B$, and 
$>_{\text{lex}}\cap Q^2\!\subseteq\! B$ or $>_{\text{lex}}\cap Q^2\!\subseteq\!\neg B$. As $B$ is irreflexive and no Cantor subset of $2^\omega$ is $B$-discrete, we cannot have $<_{\text{lex}}\cap Q^2\!\subseteq\!\neg B$ and 
$>_{\text{lex}}\cap Q^2\!\subseteq\!\neg B$.\medskip

\noindent\bf Case 1.2.1.\rm\ $<_{\text{lex}}\cap Q^2\!\subseteq\! B$ and $>_{\text{lex}}\cap Q^2\!\subseteq\! B$.\medskip

 Note that $Q^2\!\setminus\!\Delta (Q)\! =\! B\cap Q^2$ and $(2^\omega ,\not= )\sqsubseteq_c(Y,B)$.\medskip

\noindent\bf Case 1.2.2.\rm\ $<_{\text{lex}}\cap Q^2\!\subseteq\! B$ and $>_{\text{lex}}\cap Q^2\!\subseteq\!\neg B$.\medskip

 Note that $<_{\text{lex}}\cap Q^2\! =\! B\cap Q^2$ and $(2^\omega ,<_{\text{lex}})\sqsubseteq_c(Y,B)$.\medskip

\noindent\bf Case 1.2.3.\rm\ $<_{\text{lex}}\cap Q^2\!\subseteq\!\neg B$ and $>_{\text{lex}}\cap Q\!\subseteq\! B$.\medskip

 Note that $>_{\text{lex}}\cap Q\! =\! B\cap Q$ and $(2^\omega ,>_{\text{lex}})\sqsubseteq_c(Y,B)$. But 
$(2^\omega ,<_{\text{lex}})\sqsubseteq_c(2^\omega ,>_{\text{lex}})$, with witness $i$ defined by 
$i(\alpha )(n)\! :=\! 1\! -\!\alpha (n)$. Thus $(2^\omega ,<_{\text{lex}})\sqsubseteq_c(Y,B)$.\medskip

\noindent\bf Case 2.\rm\ $B\cap (N_1\!\times\! N_0)$ is meager.\medskip

 By 19.6 in [K], $(\neg B)\cap (N_1\!\times\! N_0)$ contains a product of Cantor sets, so that we may assume that $N_0\!\times\! N_1\!\subseteq\! B\!\subseteq\!\neg (N_1\!\times\! N_0)$.\medskip

\noindent\bf Case 2.1.\rm\ There is a $B$-discrete Cantor subset of $2^\omega$. Then 
$(\mathbb{S},\mathbb{H})\sqsubseteq_c(Y,B)$ or $(\mathbb{S},\mathbb{V})\sqsubseteq_c(Y,B)$.\medskip

 Indeed, assume for example that $Q\!\subseteq\! N_0$ is a Cantor $B$-discrete set. Then as in Case 1.1 we see that 
$(\mathbb{S},\mathbb{H})\sqsubseteq_c(Y,B)$. Similarly, if $Q\!\subseteq\! N_1$ is a Cantor $B$-discrete set, then 
$(\mathbb{S},\mathbb{V})\sqsubseteq_c(Y,B)$.\medskip

\noindent\bf Case 2.2.\rm\ No Cantor subset of $2^\omega$ is $B$-discrete. Then 
$(2^\omega ,\not= )\sqsubseteq_c(Y,B)$ or $(2^\omega ,<_{\text{lex}})\sqsubseteq_c(Y,B)$.\medskip

 Indeed, we argue as in Case 1.2.\medskip
 
\noindent (2) The indicated elements are the only graphs in $\mathcal{A}^c$, up to the isomorphism $(\varepsilon ,\alpha )\!\mapsto\!\varepsilon\alpha$.\hfill{$\square$}

\section{$\!\!\!\!\!\!$ References}

\noindent [B]\ \ J. E. Baumgartner, Partition relations for countable topological spaces,~\it J. Combin. Theory Ser. A\rm ~43 (1986), 178-195

\noindent [C-L-M]\ \ J. D. Clemens, D. Lecomte and B. D. Miller, Essential countability of treeable equivalence relations,~\it Adv. Math.\rm ~265 (2014), 1-31

\noindent [G]\ \ S. Gao,~\it Invariant Descriptive Set Theory,~\rm Pure and Applied Mathematics, A Series of Monographs and Textbooks, 293, Taylor and Francis Group, 2009

\noindent [Ka]\ \ V. Kanovei,~\it Borel equivalence relations,~\rm Amer. Math. Soc., 2008

\noindent [K1]\ \ A. S. Kechris,~\it Classical Descriptive Set Theory,~\rm Springer-Verlag, 1995

\noindent [K2]\ \ A. S. Kechris, The theory of countable Borel equivalence relations,~\it preprint (see the url 
http://www.math.caltech.edu/~kechris/papers/lectures on CBER02.pdf)\ \rm

\noindent [K-Ma]\ \ A. S. Kechris, and A. S. Marks, Descriptive graph combinatorics,\ \it preprint (see the url http://www.math.caltech.edu/~kechris/papers/combinatorics20.pdf)\ \rm

\noindent [K-S-T]\ \ A. S. Kechris, S. Solecki and S. Todor\v cevi\'c, Borel chromatic numbers,~\it Adv. Math.\ \rm141 (1999), 1-44

\noindent [L1]\ \ D. Lecomte, On minimal non potentially closed subsets of the plane,
\ \it Topology Appl.~\rm 154, 1 (2007) 241-262

\noindent [L2]\ \ D. Lecomte, How can we recognize potentially $\bormxi$ subsets of the plane?,~\it  J. Math. Log.\ \rm  9, 1 (2009), 39-62

\noindent [L3]\ \ D. Lecomte, Potential Wadge classes,~\it\ Mem. Amer. Math. Soc.,\rm ~221, 1038 (2013)

\noindent [L4]\ \ D. Lecomte, Acyclicity and reduction,~\it\ Ann. Pure Appl. Logic\rm ~170, 3 (2019) 383-426

\noindent [L5]\ \ D. Lecomte, On the complexity of Borel equivalence relations with some countability property,~\it\ Trans. Amer. Math. Soc.\rm\ 373, 3 (2020), 1845-1883

\noindent [L-M]\ \ D. Lecomte and B. D. Miller, Basis theorems for non-potentially closed sets and graphs of uncountable Borel chromatic number,\ \it J. Math. Log.\rm\ 8 (2008), 121-162

\noindent [L-Z]\ \ D. Lecomte and R. Zamora, Injective tests of low complexity in the plane,~\it Math. Logic Quart.\rm ~65, 2 (2019), 134-169

\noindent [Lo1]\ \ A. Louveau, A separation theorem for $\Ana$ sets,\ \it Trans. Amer. Math. Soc.\ \rm 260 (1980), 363-378

\noindent [Lo2]\ \ A. Louveau, Ensembles analytiques et bor\'eliens dans les 
espaces produit,~\it Ast\'erisque (S. M. F.)\ \rm 78 (1980)

\noindent [Lo3]\ \ A. Louveau, Two Results on Borel Orders,\ \it J. Symbolic Logic\ \rm 54, 3 (1989), 865-874

\noindent [Lo-SR1]\ \ A. Louveau and J. Saint Raymond, Borel classes and closed games: 
Wadge-type and Hurewicz-type results,\ \it Trans. Amer. Math. Soc.\ \rm 304 (1987), 431-467

\noindent [Lo-SR2]\ \ A. Louveau and J. Saint Raymond, The strength of Borel Wadge determinacy,\ \it 
Cabal Seminar 81-85, Lecture Notes in Math.\ \rm 1333 (1988), 1-30

\noindent [Mo]\ \ Y. N. Moschovakis,~\it Descriptive set theory,~\rm North-Holland, 1980

\noindent [P]\ \ T. C. Przymusi\' nski, On the notion of $n$-cardinality,~\it Proc. Amer. Math. Soc.\rm~69 (1978), 333-338

\noindent [S]\ \ J. H. Silver, Counting the number of equivalence classes of Borel and coanalytic
equivalence relations,~\it Ann. Math. Logic\ \rm 18, 1 (1980), 1-28

\noindent [T]\ \ S. Todor\v cevi\'c,~\it Introduction to Ramsey spaces~\rm Annals of Mathematics Studies, 174. Princeton University Press, Princeton, NJ, 2010. viii+287 pp
 
\noindent [W]\ \ W. W. Wadge,~\it Reducibility and determinateness on the Baire space, Ph. D. Thesis,~\rm University of California, Berkeley, 1983

\end{document}